%% file: splitting.tex
\documentclass[a4paper,DIV=12,11pt]{scrartcl}

\usepackage[utf8]{inputenc}
\usepackage{amsfonts,amssymb,amsmath}
\usepackage{epsfig}

\usepackage{graphics}
\usepackage{subfig}

\usepackage{tikz}
\usetikzlibrary{plotmarks}
\usetikzlibrary{shapes,arrows}
\usetikzlibrary{positioning}
\usetikzlibrary{trees,mindmap}
\usetikzlibrary{fadings,shapes.arrows,shadows}
\usetikzlibrary{decorations.pathreplacing}
\usetikzlibrary{calc}
\usepackage{pgfplots,pgflibraryplotmarks}

\renewcommand{\vec}{\boldsymbol}

\newcommand{\ud}{\,\mathrm{d}}
\newcommand{\R}{\mathbb{R}}
\newcommand{\N}{\mathbb{N}}

\numberwithin{equation}{section}
\numberwithin{figure}{section}

\def\veps{{\vec \varepsilon}}
\def\la{{\langle}}
\def\ra{{\rangle}}

\newtheorem{defi}{Definition}[section]
\newtheorem{theorem}[defi]{Theorem}

\newtheorem{rem}[defi]{Remark}
\newtheorem{cor}[defi]{Corollary}
\newenvironment{mproof}{\paragraph{Proof.}}{\hfill$\blacksquare$}

\begin{document}

\title{Space-time finite element approximation\\ of the 
Biot poroelasticity system\\ with iterative coupling}
\author{M.\ Bause\thanks{bause@hsu-hh.de (corresponding author),
$^\dag$florin.radu@uib.no, $^\ddag$koecher@hsu-hh.de}\mbox{~}, F.\ 
A.\ Radu$^\dag$, U.\ K\"ocher$^\ddag$ \\
{\small $^\ast \, \ddag$ Helmut Schmidt University, Faculty of 
Mechanical Engineering, Holstenhofweg 85,}\\ 
{\small 220433 Hamburg, Germany}\\
{\small $^\dag$ University of Bergen, Department of Mathematics, 
All\'{e}gaten 41,}\\{\small 50520 Bergen, Norway}
}

\maketitle

\begin{abstract}
We analyze an optimized artificial fixed-stress iteration scheme for 
the numerical approximation of the Biot system modelling fluid flow in 
deformable porous media. The iteration is based on a prescribed 
constant artificial volumetric mean total stress in the first half step. 
The optimization comes through the adaptation of a numerical stabilization 
or tuning parameter and aims at an acceleration of the iterations. The 
separated subproblems of fluid flow, written as a mixed first order in space 
system, and mechanical deformation are discretized by space-time 
finite element methods of arbitrary order. Continuous and 
discontinuous discretizations of the time variable are encountered. 
The convergence of the iteration schemes is proved for the continuous and fully discrete 
case. The choice of the optimization parameter is identified in the proofs of convergence 
of the iterations. The analyses are illustrated and confirmed by numerical experiments.   
\end{abstract}


\section{Introduction and mathematical model}
\label{Sec:Introduction}

Many physical and technical problems in mechanical, environmental and petroleum 
engineering as well as in biomechmanics and biomedicine involve interactions between flow 
and mechanical deformation in porous media. Therefore, the ability to simulate coupled 
mechanical deformations and fluid flow in such media is of particular importance from 
the point of view of physical realism. Numerical modelling of such coupled processes is 
complex due to the structure of the equations and continues to remain a challenging 
task. 

We consider modelling flow in deformable porous media by the quasi-static 
Biot system \cite{MW12},
\begin{align}
\label{Eq:B_1}
- \nabla \cdot \vec \sigma (\vec u,p) = \rho_b \vec g\,,\\[0ex]
\label{Eq:B_2}
\partial_t\left(\frac{1}{M}p + \nabla \cdot (b \vec u)\right) + \nabla \cdot 
\vec q  = f \,, \quad \vec q = - \frac{\vec K}{\eta} \left(\nabla 
p -\rho_f \vec g \right)\,, \\[0ex]
\label{Eq:B_3}
p(0) = p_0\,, \quad \vec u(0) = \vec 0\,,
\end{align}
with supplemented boundary conditions and the total stress $\vec \sigma(\vec u,p) = \vec 
\sigma_0 + \vec C:\vec \varepsilon (\vec u) - b(p-p_0)\vec I$, to be satisfied in the 
bounded Lipschitz domain $\Omega\subset \R^{d}$, with $d=2$ or $d=3$, and for the time 
$t\in I=(0,T]$. In \eqref{Eq:B_1}--\eqref{Eq:B_3}, we denote by $\vec u$ the unknown 
displacement field, $p$ the unknown fluid pressure, $\vec \varepsilon(\vec u) = 
(\nabla \vec u +(\nabla \vec u)^\top)/2$ the linearized strain tensor, $\vec C$ the 
Gassmann rank-4 tensor of elasticity, $\vec \sigma_0$ the reference state stress tensor, 
$b$ Biot's coefficient, $\rho_b=\phi \rho_f + (1-\phi)\rho_s$ the bulk density with 
porosity $\phi$ and fluid and solid phase density $\rho_f$ and $\rho_s$, $p_0$ the 
reference state fluid pressure, $M$ Biot's modulus and, finally, by $\vec q$ Darcy's 
velocity or the fluid flux. Eq.\ \eqref{Eq:B_1} models conservation of momentum and the 
first of the equations \eqref{Eq:B_2} describes conservation of mass. The second of the 
equations \eqref{Eq:B_2} is the well-known Darcy law with permeability field $\vec K$ and 
fluid viscosity $\eta$. Further, $\vec g$ denotes gravity or, in general, some body force 
and $f$ is a volumetric source. The quantities $\eta$, $M$, $\rho_{f}$ and $\rho_s$ are 
positive constants. The matrix $\vec K$ is supposed to be symmetric and uniformly 
positive definite. For any symmetric matrix $\vec B$ we assume that $ (\vec C \vec 
B):\vec B \geq a |\vec B |^2 + \lambda \mathrm{tr}(\vec B)^2$ is satisfied with some 
constant $a>0$ and the drained bulk modulus $\lambda$; cf.\ \cite{MWW14}. 
We assume that $\rho_b$ is independent of time and that $\rho_b\vec g = - \nabla \cdot 
\vec \sigma_0$. Here, the quasi-static feature is due to the negligence of the solid's 
acceleration in problem \eqref{Eq:B_1} of mechanical deformation. This prevents the 
applicability of the model \eqref{Eq:B_1}--\eqref{Eq:B_3} to classes of problems for that 
the contrast coefficients, the ratio between the intrinsic characteristic time and the 
characteristic time scale of the domain, are not close to the singular limit of vanishing 
numbers. In \cite{S00} an existence, uniqueness and regularity theory is presented for 
the Biot system \eqref{Eq:B_1}--\eqref{Eq:B_3} in a Hilbert space setting. In 
\cite{S04} the well-posedness is shown for a wider class of diffusion problems in 
poro-elastic media with more general material deformation models. 

Iteratively coupled solution methods for the system of \eqref{Eq:B_1}--\eqref{Eq:B_3} of 
coupled fluid flow and mechanical deformation have recently attracted researchers' 
interest; cf.\ \cite{A16,C15,C16,K11,MW12,MW13,MWW14,PW07,S98} and the references 
therein. Iterative 
coupling is a sequential approach, in that either the problem of flow or the mechanics is 
solved first followed by solving the other system using the already calculated 
information. At each time step this is repeated until a converged solution within a 
prescribed numerical tolerance is obtained. In \cite{K11} it's shown by an analysis that 
some of the splitting approaches may exhibit stability problems. Iterative coupling 
offers the appreciable advantage over the fully coupled method that existing and highly 
developed discretizations and algebraic solver technology, for instance 
preconditioning methods, as well as existing software tools can be reused. The 
construction of efficient preconditioning techniques for solving the arising algebraic 
systems of equations of fully coupled approaches to poroelasticity has not been 
satisfactorily solved yet and continues to remain a field of active research \cite{W16}. 
In particular, this applies to the case in that higher order space and time 
discretization 
techniques are involved. 

In this work we analyze a ''fixed-stress split'' type iterative method; cf.\ 
\cite{K11,MW13}. The fixed-stress split iterative method is based on imposing constant 
volumetric mean total stress $\sigma_v = \sigma_0 + \lambda \nabla \cdot \vec u - 
b (p-p_0)$ in the first half step of fluid flow. In our approach we use some 
optimized fixed-stress split by prescribing a constant artificial volumetric mean total 
stress that is given for $\vec \sigma_0=0$ and $p_0=0$ by $\sigma_v = \lambda 
\nabla \cdot \vec u -L \lambda b^{-1} p$ with some additional numerical parameter 
$L>0$ that has to be tuned to accelerate the iteration procedure and to reduce the 
numbers of iterations that are required for the adherence of a prescribed numerical 
tolerance. In contrast to \cite{MW13}, a mixed formulation of the flow problem 
\eqref{Eq:B_2} is considered here. In this paper we proof the convergence of the 
proposed iteration scheme by a fixed point argument and identify an optimal choice for 
the numerical parameter $L$. This is done for the continuous case of the iteratively 
coupled subproblems of partial differential equations and for the fully discrete case of 
space-time finite approximations of the subproblems. Our analysis yields the same choice 
for the numerical parameter $L$ for the either cases, even though completely different 
techniques of proof are used. Therefore, the acceleration of the iteration's 
convergence is not impacted by the time or space step size or the polynomial degree of 
the finite element methods in time and space. Our numerical tests nicely confirm the 
choice of the acceleration parameter $L$ that is suggested by our numerical analysis of 
the schemes. The numerical results show that the number of required iterations can 
strongly be reduced by using the proposed optimized fixed-stress split iterative method 
along with the suggested choice of the tuning parameter $L$.

For the numerical approximation of the separated subproblems of fluid flow and 
mechanical deformation we use space-time finite element methods. Continuous and 
discontinuous finite element discretizations of the time variable are studied. For the 
spatial discretization of the flow problem mixed finite element methods (cf.\ \cite{C10}) 
ensuring local mass conservation and an inherent approximation of the flux variable are 
used. Due to these properties, mixed finite element methods have 
shown in numerous works their superiority over standard conforming methods for the 
numerical simulation of fluid flow in porous media; cf.\ \cite{JJ07} for its application 
to reservoir geomechanics. For the spatial discretization of the displacement variable a 
standard conforming approach is used in order to simplify the analysis. In the future we 
will use discontinuous Galerkin methods for the discretization of the displacement 
field and the approximation of the subproblem of mechanical deformation since we expect 
from our former works (cf.\ \cite{K14,K15}) on discontinuous Galerkin methods significant 
advantages for future generalizations of the underlying Biot model, for instance, to the 
Biot--Allard system \cite{MW12}. Moreover, the discontinuous Galerkin discretization of 
the displacement variable helps to avoid locking phenomena. For a discussion of locking 
phenomena arising in poroelasticity and remedies we refer to \cite{L13,N16,PW08,PW09,R16} 
and the references therein.  

Since recently, variational time discretization schemes based on continuous or 
discontinuous finite element techniques have been developed to the point
that they can be put into use (cf., e.g., \cite{A17,A11,BK15,DF15,HST12,HST13,K15,S15} 
and the references 
therein) and demonstrated their significant advantages. Higher order methods are 
naturally embedded in these schemes and the uniform variational approach simplifies 
stability and error analyses. Further, goal-oriented error control \cite{BR03} based on 
the dual weighted residual approach relies on variational space-time formulations and the 
concepts of adaptive finite element techniques for changing the polynomial degree as well 
as the length of the time intervals become applicable. However, in the context of 
numerical modelling flow in porous or deformable media higher order space-time finite 
element methods or even only higher order time discretizations have rarely been used in 
practice so far. However, for applications with strong fluctuations of 
physical quantities and involved highly dynamical processes, for instance in 
vibro acoustics and reactive multicomponent and multiphase subsurface flow, as well as 
for future generalizations to more complex models like the Biot--Allard system 
\cite{MW12} the needfulness of developing and analyzing higher order techniques is 
evident.   

The paper is organized as follows. In Sec.\ \ref{Sec:IterScheme} we introduce the 
iterative coupling scheme of subproblems of partial differential equations for fluid flow 
and mechanical deformation and prove its convergence. In Sec.\ \ref{Sec:Discretization} 
the space-time finite element discretization of the subproblems is introduced for a 
continuous and a discontinuous approximation of the time variable. In 
Sec.~\ref{Sec:Convergence} we then prove convergence of the iterations for both families 
of space-time finite element approximations. Sec.\ \ref{Sec:NumRes} illustrates the given 
analyses by numerical computations and confirms our theoretical observations. Sec.\ 
\ref{Sec:Summary} summarizes the results of this work. 

Throughout the paper, our notation is standard. We denote by $H^m(\Omega)$ the Sobolev 
space of $L^2$ functions with derivatives up to order $m$ in $L^2(\Omega)$. By $\langle 
\cdot,\cdot \rangle$ and $ \| \cdot \|$ denote the inner product and norm in 
$L^2(\Omega)$, respectively, where we do not differ in the notation between inner 
products and norms of scalar- and vector-valued functions. For rank-2 tensors $\vec A,\vec 
B\in \R^{d,d}$ we use the notation $\langle \vec A,\vec B\rangle = \int_\Omega 
\sum_{i,j=1}^d A_{ij}B_{ij} \ud \vec x$. Further, let $H^1_0(\Omega)=\{u\in 
H^1(\Omega) \mid u=0 \mbox{ on } \partial \Omega\}$. For the mixed problem formulation of 
the flow problem \eqref{Eq:B_2} we put
\begin{equation*}
 \vec V = \vec H(\mathrm{div};\Omega)\,, \qquad W = L^2(\Omega)\,,
\end{equation*}
where $\vec H(\mathrm{div};\Omega) = \{\vec q \in \vec L^2(\Omega) \mid \nabla
\cdot \vec q \in L^2(\Omega)\}$. Let $X_0 \subset X \subset X_1$ be
three reflexive Banach spaces with continuous embeddings. Then we consider
the following set of spaces,
\begin{align*}
C(\overline{I};X) & =
\{ w : [0,T] \rightarrow X \mid \mbox{$w$ is continuous} \}\,,\\
L^2(I;X) & =
\bigg\{w: (0,T) \rightarrow X \;\; \bigg|\;\;
\int_0^T \| w(t) \|_X^2 \; \ud t < \infty \bigg\}\,,\\
H^1(I;X_0,X_1) & =
\{w \in L^2(I;X_0) \mid \partial_t w \in L^2(I;X_1)\}\,,
\end{align*}
that are equipped with their naturals norms (cf.\ \cite{Ern2010}) and where the
time derivative $\partial_t$ is understood in the sense of distributions on
$(0,T)$. In particular, every function in $H^1(I;X_0,X_1)$ is continuous on
$[0,T]$ with values in $X$; cf.\ \cite{Ern2010}. For $X_0=X=X_1$ we simply
write $H^1(I;X)$.

\section{Iterative coupling scheme and proof of convergence}
\label{Sec:IterScheme}

In this section we introduce our iterative coupling scheme of subproblems of partial 
differential equations and prove its convergence. The fully discrete counterpart of 
the scheme and its convergence is treated below in Sec.\ \ref{Sec:Discretization} and 
\ref{Sec:Convergence}, respectively. In our analysis of the scheme we restrict ourselves 
to homogeneous Dirichlet boundary conditions. In our numerical experiments (cf.\ 
Section \ref{Sec:NumRes}) more general boundary conditions are also encountered. Without 
loss of generality, we assume vanishing initial conditions $p_0 =0$ and $\vec u_0 = \vec 
0$. We put $\vec \sigma_0= \vec 0$ and assume that $\vec g(0) =\vec 0$. Further, we 
prescribe an isotropic material behavior such that the Gassmann rank-4 tensor of 
elasticity $\vec C$ is given by $c_{ijkl} = \lambda \delta_{ij}\delta_{kl}+\mu 
(\delta_{ik}\delta_{jl}+\delta_{il}\delta_{jk})$ and the total stress reads as $\vec 
\sigma(\vec u,p) = 2\mu \vec \varepsilon(\vec u) + \lambda \nabla \cdot \vec u \vec I - 
b\,p\vec I$ with $\mu>0$ and $\lambda$ denoting the Lam\'{e} parameters. We let $\lambda 
>0$ which is satisfied for most of the materials. These assumptions can be relaxed 
further to study more complex and non isotropic material behavior. To simplify the 
notation, we write $\vec K$ instead of $\vec K/\eta$ and add the gravity term of 
\eqref{Eq:B_2} to the right-hand side term $f$. Under these assumptions Eqs.\ 
\eqref{Eq:B_1}--\eqref{Eq:B_3} 
read as 
\begin{align}
\label{Eq:B_4}
- \nabla \cdot \left(2\mu \vec \varepsilon (\vec u) + \lambda \nabla \cdot \vec u \vec I 
- b\,p\vec I\right) = \rho_b \vec g\,,\\[0ex]
\label{Eq:B_5}
\partial_t\left(\frac{1}{M}p + \nabla \cdot (b \vec u)\right) + \nabla \cdot 
\vec q  = f \,, \quad \vec q = - \vec K \nabla 
p \,, \\[0ex]
\label{Eq:B_6}
p(0) = 0\,, \quad \vec u(0) = \vec 0
\end{align}
for $\vec x \in \Omega$ and $t\in I$ with the boundary conditions 
\begin{align}
\label{Eq:B_7}
p = 0 \quad \text{and} \quad  \vec u  = \vec 0\qquad  \mbox{ on } \partial \Omega \times 
I\,.
\end{align}
For the data $\vec g$, $f$ and $\vec K$ we assume at first that the conditions $\vec g 
\in L^2(I;\vec L^2(\Omega))$, $f \in L^2(I;L^2(\Omega))$ and $\vec K\in \vec 
L^\infty(\Omega)$ are satisfied.

For $\Omega=(0,l)^d$ and under periodic boundary conditions for $\vec u$ and $p$ 
with period $l$ and for smooth $l$-periodic functions $p_0$, $f$ and $\vec \sigma_0$ it 
is shown in \cite{MW12} that the system \eqref{Eq:B_1}--\eqref{Eq:B_3} admits a unique 
periodic solution 
\(
\{\vec u,p\} \in C(\overline I;\vec 
H^1_{\mathrm{per}}(\Omega)\cap \vec L^2_0(\Omega))\times H^1(\Omega\times I) \cap 
C(\overline{I};H^1_{\mathrm{per}}(\Omega))
\). Further, for $\vec g\in C_0^\infty(\R^+;\vec L^2_0(\Omega))$ and homogeneous 
initial conditions the solution of the system is smooth in time with 
\(
\{\vec u,p\} \in H^k(I;\vec H^1_{\mathrm{per}}(\Omega)) \times 
H^k(I; H^1_{\mathrm{per}}(\Omega))\,, \mbox{for all $k\in \N$}\,.
\); cf.\ \cite{MW12}.

To solve the equations \eqref{Eq:B_4}--\eqref{Eq:B_7} we use a \textit{fixed-stress 
iterative splitting scheme}; cf.\ \cite{MW13}. This scheme consists in imposing a 
constant artificial volumetric mean total stress $\sigma_{v}= \lambda\nabla \cdot 
\vec u - L\,\lambda\,b^{-1}\,p$ in the first half step. Here, the parameter 
$L>0$ is a free to be chosen constant that is specified below. The supplement 
''artificial'', that is used here, is due to the additional parameter $L$ in contrast to 
the proper definition of the volumetric mean total stress given by $\sigma_{v}= 
\lambda\nabla \cdot \vec u - b\, p$. By adding the parameter $L$ we aim to find 
an iteration scheme with smaller and optimal contraction number compared to the standard 
definition of $\sigma_v$; cf.\ \cite{MW13}. Supposing a constant artificial volumetric 
mean total stress then yields in the first half step of fluid flow  
\begin{equation}
\label{Eq:FSP_1}
\begin{split}
\left(\frac{1}{M}+ L \right)\partial_t p^{k+1} + \nabla 
\cdot \vec q^{k+1} = f - b \nabla \cdot \partial_t \vec u^k + L \partial_t p^k\,,\qquad
\vec q^{k+1}  = - \vec K \nabla p^{k+1} 
\end{split}
\end{equation}
on $\Omega\times I$, $p^{k+1}(0)= 0$ in $\Omega$ and $p^{k+1}=0$ on $\partial\Omega 
\times I$. In each iteration step problem \eqref{Eq:FSP_1} of fluid flow is thus 
decoupled from the mechanical deformation subproblem and can be solved independently. In 
the second half step the effective deformation is then obtained by 
solving 
\begin{equation}
\label{Eq:FSP_2}
- \nabla \cdot \left(2\mu \vec \varepsilon (\vec u^{k+1})+\lambda \nabla \cdot \vec 
u^{k+1} \vec I \right) =  \rho_b \vec g  - b \nabla p^{k+1} 
\end{equation}
on $\Omega\times I$, where $\vec u^{k+1}(0)= \vec 0$ and $\vec u^{k+1}=0$ 
on $\partial\Omega \times I$. 

The weak formulation of problem \eqref{Eq:FSP_1} in the space-time framework then reads 
as follows: \textit{Let $\widetilde f^{\,k} := f - b \, \nabla \cdot \partial_t \vec u^k 
+ L\, \partial_t p^k$  with $\widetilde f^{\,k} \in 
L^2(I;W)$ be given. Find $p^{k+1}\in H^1(I;W)$ and $\vec q^{k+1}\in L^2(I;\vec V)$ such 
that $p^{k+1}(0)=0$ and}
\begin{align}
\label{Eq:WPF1}
\left(\frac{1}{M}+ L \right) \int_I \langle \partial_t p^{k+1}, w\rangle \ud t + 
\int_I \langle \nabla \cdot \vec q^{k+1}, w\rangle \ud t & =  
\int_I \langle \widetilde f^{\,k}, w \rangle \ud t\,, \\[1ex]
\label{Eq:WPF2}
\int_I \langle \vec K^{-1} \vec q^{k+1}, \vec v\rangle \ud t - \int_I \langle 
p^{k+1}, \nabla \cdot \vec v \rangle \ud t & = 0 
\end{align}
\textit{for all $w\in L^2(I;W)$ and $\vec v\in L^2(I;\vec V)$.} 

The weak form of problem \eqref{Eq:FSP_2} reads as follows: \textit{Let $p^{k+1}\in 
H^1(I;W)$ be given. Find $\vec u^{k+1} \in H^1(I;\vec H^1(\Omega))\cap L^2(I;\vec 
H^1_0(\Omega))$ such that $\vec 
u(0)=\vec 0$ and}
\begin{equation}
\label{Eq:WPM}
\begin{aligned}
\int_I 2 \mu \langle \vec \varepsilon(\vec u^{k+1}), \vec \varepsilon(\vec z)\rangle \ud 
t + \int_I \lambda \langle & \nabla \cdot \vec u^{k+1}, \nabla \cdot \vec z \rangle \ud 
t\\
& = \rho_b \int_I \langle \vec g , \vec z \rangle \ud t +  b \int_I \langle p^{k+1}\vec 
I, \vec \varepsilon(\vec z)\rangle \ud t
\end{aligned}
\end{equation}
\textit{for all $z \in L^2(I;\vec H^1_0(\Omega))$.}

To simplify the notation, we put 
\begin{align*}
\mathcal W & =\{w \in H^1(I;H^1(\Omega))\mid w\in C(\overline 
I;H^1_{0}(\Omega))\}\,,\\[1ex]
\vec{\mathcal V} & = \{\vec v\in  L^2(I;\vec V) \mid \vec v\in  C(\overline 
I;\vec L^2(\Omega))\}\,,\\[1ex]
\vec{\mathcal Z} & = \{ \vec z \in H^1(I;\vec H^1(\Omega)) \mid \vec  z\in  
C(\overline 
I;\vec H^1_{0}(\Omega))\,, \; \partial_t \vec u \in L^2(I;\vec H^2(\Omega))\}\,.
\end{align*}

The following theorem shows the convergence of the iteration scheme \eqref{Eq:WPF1} to 
\eqref{Eq:WPM}. In contrast to \cite{MW13} our proof is based on a mixed formulation of 
the flow problem. Moreover, the proof is presented explicitly here in order to show that 
the convergence proofs for the iteration scheme on the continuous and discrete level lead 
to the same optimal parameter $L$, even though completely different techniques of proof 
are used.     

\begin{theorem}
\label{Thm:ConvCont}
Suppose that $\partial \Omega$ and the permeability field $\vec K$  are sufficiently 
regular. Let $f\in L^2(I;H^1(\Omega))$ and $\vec g\in H^1(I;\vec L^2(\Omega))$ be 
satisfied. Then, for any $L \geq b^2/(2\lambda)$ the operator $\vec{\mathcal S}: 
(p^k,\vec q^k,\vec u^k) \mapsto (p^{k+1},\vec 
q^{k+1},\vec u^{k+1})$ maps $\vec{\mathcal{D}} = \{\{p,\vec q,\vec u\} \in \mathcal W 
\times \vec{\mathcal V} \times \vec{\mathcal Z}  \mid p(0) = 0\,,\; \vec u(0) = \vec 0 
\}$ into itself and is a contraction mapping on $\vec{\mathcal D}$. Therefore, the 
operator $\vec{\mathcal S}$ has a unique fixed point in $\vec{\mathcal D}$. The 
contraction constant is smallest for $L = b^2/(2\lambda)$ with value $L\,M/(L\, M+1)$.
\end{theorem}

\begin{mproof}
Firstly, we show that $\vec{\mathcal S}$ maps $\vec{\mathcal{D}}$ into itself. For 
this, let 
$\{p^k,\vec q^k,\vec u^k\} \in \mathcal W \times \vec{\mathcal V} \times 
\vec{\mathcal Z}$ be given. Under the assumptions of the theorem it follows that 
\[
\widetilde f^{\,k} \in L^2(I;H^1(\Omega)) \quad \text{for}\quad \widetilde f^{\,k} = f - 
b \, \nabla \cdot \partial_t \vec u^k 
+ L\, \partial_t p^k\,.
\]
The variational problem \eqref{Eq:WPF1}, \eqref{Eq:WPF2} then admits a unique solution 
$p^{k+1}\in \mathcal{W}$ and $\vec q^{k+1}\in \vec{\mathcal V}$. This directly follows 
from parabolic regularity theory;  cf., e.g., \cite{E10}. For $p^{k+1}\in 
H^1(I;H^1(\Omega))$ the second of the right-hand side terms in Eq.\ \eqref{Eq:WPM} can be 
rewritten as 
\[
\int_I 
\left\langle p^{k+1}\vec I, \vec \varepsilon(\vec z)\right\rangle \ud \tau  = \int_I 
\left\langle \nabla p^{k+1}, \vec z\right\rangle \ud \tau \,.
\]
By means of elliptic regularity theory the variational problem \eqref{Eq:WPM} then 
admits a unique solution $\vec u^{k+1}\in \vec{\mathcal Z}$; cf., e.g., \cite{CW98}.  We 
note that $p^{k+1}\in \mathcal{W}$, $\vec q^{k+1}\in \vec{\mathcal V}$ and $\vec 
u^{k+1}\in \vec{\mathcal Z}$ are even strong solutions of the problems \eqref{Eq:FSP_1} 
and \eqref{Eq:FSP_2}, respectively.

Secondly, we now show that the operator $\vec{\mathcal S}$ is a contraction mapping on 
$\vec{\mathcal{D}}$. With 
\begin{equation}
\label{Def:VMS}
\sigma_v =  \lambda\nabla \cdot \vec u - \dfrac{L\,\lambda}{b}\,p
\end{equation}
and $S_p^{k+1}=p^{k+1}-p^k$, $\vec S_{\vec q}^{k+1}=\vec q^{k+1}-\vec 
q^k$, $\vec S_{\vec u}^{k+1}=\vec u^{k+1}-\vec u^k$, $S_{\sigma_v}^{k}= \sigma_v^k - 
\sigma_v^{k-1}$ for the differences of the iterates we get from the first of 
the equations \eqref{Eq:WPF1} that 
\begin{equation}
 \label{Eq:Con_01}
\left(\frac{1}{M}+ L \right)\int_0^t \left\langle \partial_t S_p^{k+1},w\right\rangle \ud 
\tau   + 
\int_0^t \left\langle \nabla 
\cdot \vec S_{\vec q}^{k+1},w\right\rangle \ud \tau = - \int_0^t \frac{b}{\lambda} 
\left\langle \partial_t S_{\sigma_v}^k,w \right\rangle \ud \tau  
\end{equation}
for all $w\in L^2(I;W)$. Choosing $w=S_p^{k+1}$ in Eq.\ \eqref{Eq:Con_01} and using the 
inequalities of Cauchy--Schwarz and Cauchy--Young we obtain that  
\begin{equation*}
\begin{aligned}
& \left(\frac{1}{M}+ L \right) \dfrac{b^2}{L^2 \lambda^2} \int_0^t \int_\Omega 
\bigg| 
\frac{L \lambda}{b} \partial_t S_p^{k+1}\bigg|^2 \ud \vec x\ud \tau + \int_0^t 
\left\langle \nabla \cdot \vec S_{\vec q}^{k+1}, \partial_t 
S_p^{k+1}\right\rangle  \ud \tau \\[2ex]
& = - \int_0^t \frac{b}{\lambda} \left\langle \partial_t S_{\sigma_v}^k , 
\partial_t S_p^{k+1} \right\rangle \ud \tau\\[2ex]
& \leq \frac{\varepsilon}{2} \, \dfrac{b^2}{L^2 \lambda^2} \int_0^t \int_\Omega 
\bigg| 
\frac{L \lambda}{b} \partial_t S_p^{k+1}\bigg|^2 \ud \vec x \ud \tau + 
\frac{b^2}{2\varepsilon \lambda^2}\int_0^t \int_\Omega  \left|\partial_t 
S_{\sigma_v}^k\right|^2 
\ud \vec x \ud \tau\,.
\end{aligned}
\end{equation*}
Choosing $\varepsilon = L + \frac{1}{M}$, we then get that 
\begin{equation}
\label{Eq:Con_02}
\begin{aligned}
\int_0^t \int_\Omega \bigg| \frac{L \lambda}{b} \partial_t S_p^{k+1}\bigg|^2 \ud 
\vec x \ud \tau &  + \gamma \int_0^t \left\langle \nabla \cdot \vec S_{\vec q}^{k+1}, 
\partial_t S_p^{k+1}\right\rangle \ud \tau \\[2ex]
& \leq \left(\frac{L}{L+1/M}\right)^2 \int_0^t \int_\Omega  \left|\partial_t 
S_{\sigma_v}^k\right|^2 
\ud \vec x \ud \tau
\end{aligned}
\end{equation}
with $\gamma = \frac{2L^2 \lambda^2}{b^2} \cdot \frac{1}{L+\frac{1}{M}} >0$.

Next, taking the time derivative of the second of the equations \eqref{Eq:FSP_1} and 
testing the resulting identity with $\vec v = \vec S_{\vec q}^{k+1}$, we 
have that 
\begin{equation*}
\int_0^t \left\langle \vec K^{-1} \partial_t \vec S_{\vec q}^{k+1}, \vec S_{\vec 
q}^{k+1}\right\rangle \ud \tau - \int_0^t \left\langle \nabla \cdot \vec 
S_{\vec q}^{k+1}, 
\partial_t S_p^{k+1}\right\rangle \ud \tau  = 0\,.
\end{equation*}
By means of $\frac{1}{2}\frac{\ud }{\ud \tau} \langle a, a\rangle = \langle \partial_t 
a,a\rangle$ we conclude from the previous equation that 
\begin{align}
\nonumber
&\int_0^t \left\langle\nabla \cdot \vec S_{\vec q}^{k+1}, 
\partial_t S_p^{k+1}\right\rangle \ud \tau = \int_0^t \dfrac{1}{2}\frac{\ud}{\ud \tau} 
\left\langle \vec K^{-1} \vec S_{\vec q}^{k+1}, \vec S_{\vec q}^{k+1}\right\rangle \ud 
\tau \\[2ex]
\nonumber
& = \frac{1}{2} \left\langle \vec K^{-1}\vec S_{\vec q}^{k+1}(t), \vec S_{\vec 
q}^{k+1}(t)\right\rangle -  \frac{1}{2} \left\langle \vec K^{-1}\vec S_{\vec q}^{k+1}(0), 
\vec 
S_{\vec q}^{k+1}(0)\right\rangle\\[1ex] 
\label{Eq:Con_03}
& = \frac{1}{2} \left\|\vec K^{-1/2} \vec S_{\vec q}^{k+1}(t)\right\|^2\,.
\end{align}
We note that by definition and Eq.\ \eqref{Eq:FSP_1} along with the regularity 
conditions of $\vec{\mathcal{D}}$ it holds that $\vec S_{\vec q}^{k+1}(0) = -\vec K 
\nabla p_0 + \vec K \nabla p_0 = \vec 0$.

Finally,  taking the time derivative of Eq.\ \eqref{Eq:FSP_2} and testing the 
resulting equation with $\vec z = \partial_t \vec S_{\vec u}^{k+1}$ we get that 
\begin{equation}
\label{Eq:Con_04}
\begin{aligned}
\frac{2\, L\, \lambda^2}{b^2} \int_0^t \int_\Omega 2\mu  \left| \vec 
\varepsilon(\partial_t \vec S_{\vec u}^{k+1})\right|^2& \ud \vec x \ud \tau  + 
\frac{2\,\lambda \, L\, \lambda^2}{b^2} \int_0^t \int_\Omega \left|\nabla \cdot 
\partial_t \vec S_{\vec u}^{k+1}\right|^2 \ud \vec x \ud \tau \\[2ex]
& = 2 \int_0^t \left\langle \frac{L\,\lambda }{b}\partial_t 
S_p^{k+1},\lambda\nabla 
\cdot \left(\partial_t \vec S_{\vec u}^{k+1}\right)\right\rangle \ud \tau\,.
\end{aligned}
\end{equation}
Applying the algebraic identity 
\begin{equation*}
2 \langle a,b\rangle = \langle a,a\rangle + \langle b,b\rangle - \langle a-b,a-b\rangle
\end{equation*}
to the right-hand side of Eq.\ \eqref{Eq:Con_04} and recalling definition \eqref{Def:VMS} 
we find that 
\begin{equation}
\label{Eq:Con_05}
\begin{aligned}
& \frac{2\, L\, \lambda^2}{b^2} \int_0^t \int_\Omega 2\mu \left| \vec 
\varepsilon(\partial_t \vec S_{\vec u}^{k+1})\right|^2 \ud \vec x \ud \tau  + 
\frac{2\,\lambda \, L}{b^2} \int_0^t \int_\Omega \left|\lambda 
\nabla \cdot 
\partial_t \vec S_{\vec u}^{k+1}\right|^2 \ud \vec x \ud \tau \\[2ex]
& = \int_0^t \int_\Omega \left|\frac{L\, \lambda}{b} \partial_t 
S_p^{k+1}\right|^2 \ud \vec x\ud \tau + \int_0^t \int_\Omega \left|\lambda \nabla 
\cdot \partial_t \vec S_{\vec u}^{k+1}\right|^2 \ud \vec x \ud \tau\\[2ex]
& \quad - \int_0^t \int_\Omega \left|\partial_t S_{\sigma_v}^{k+1}\right|^2 \ud \vec x 
\ud \tau \,.
\end{aligned}
\end{equation}

Finally, summing up the relations \eqref{Eq:Con_02}, \eqref{Eq:Con_03} and 
\eqref{Eq:Con_05} yields that 
\begin{equation}
\label{Eq:Con_06}
\begin{aligned}
&\int_0^t \int_\Omega \left|\partial_t S_{\sigma_v}^{k+1}\right|^2 \ud \vec x \ud 
t  +  \frac{2\, L\, \lambda^2}{b^2} \int_0^t \int_\Omega 2\mu \left| \vec 
\varepsilon(\partial_t \vec S_{\vec u}^{k+1})\right|^2 \ud \vec x \ud \tau  \\[2ex]
& \quad + \left(\frac{2\lambda\, L}{b^2}-1\right) \int_0^t \int_\Omega 
\left|\lambda \nabla \cdot 
\partial_t \vec S_{\vec u}^{k+1}\right|^2 \ud \vec x \ud \tau + \frac{\gamma}{2} 
\left\|\vec 
K^{-1/2} \vec S_{\vec q}^{k+1}(t)\right\|^2\\[2ex]
& \leq \left(\frac{L}{L+1/M}\right)^2 \int_0^t \int_\Omega  \left|\partial_t 
S_{\sigma_v}^k\right|^2 
\ud \vec x \ud \tau\,.
\end{aligned}
\end{equation}
Inequality \eqref{Eq:Con_06} yields a contraction map only if $L \geq b^2/(2\lambda)$. 
The contraction constant is smallest for $L=b^2/(2 \lambda)$.  

On the space $\vec{\mathcal{D}}$ the expression on the left-hand side of 
\eqref{Eq:Con_06} defines a metric by 
\begin{equation*}
\label{Eq:Con_07}
\begin{aligned}
& d_{\vec{\mathcal{D}}} \Big((\vec u,p,\vec q),(\vec 0,0,\vec 0)\Big)\\[1ex]
& = \frac{4\, L\, \lambda^2\, \mu}{b^2} \|\vec 
\varepsilon(\partial_t {\vec u}) \|_{\vec L^2(\Omega\times I)}^2 + 
\left(\frac{2\lambda\, L}{b^2}-1\right)\|\lambda 
\nabla \cdot \partial_t {\vec u} \|_{L^2(\Omega\times I)}^2 \\[1ex]
& \qquad + \left\|\partial_t \left(\lambda \nabla \cdot \vec u - \frac{L 
\lambda}{b} p \right)\right\|_{L^2(\Omega\times I)}^2 + \frac{\gamma}{2} 
\max_{0\leq 
t \leq T} \|\vec K^{-1/2} \vec q(t)\|^2_{\vec L^2(\Omega)}\,.
\end{aligned}
\end{equation*}
Summarizing the previous steps, we note that the operator $\vec{\mathcal S}$ mapsto 
$\vec{\mathcal{D}}$ into itself and, by inequality \eqref{Eq:Con_06}, satisfies
\begin{equation*}
d_{\vec{\mathcal{D}}} \Big((\vec u^{k+1},p^{k+1},\vec q^{k+1})-(\vec u^{k},p^{k},\vec 
q^{k})\Big) \leq \delta \, d_{\vec{\mathcal{D}}}\Big((\vec 
u^{k},p^{k},\vec q^{k})-(\vec u^{k-1},p^{k-1},\vec q^{k-1})\Big)
\end{equation*}
with $\delta = \frac{L}{L+\frac{1}{M}}$. Therefore, the operator $\vec{\mathcal S}$ is a 
contraction mapping and by the contraction mapping principle, it has a unique fixed point.
\end{mproof}

\section{Space-time discretization}
\label{Sec:Discretization}

In this section we introduce our space time finite approximation of the subproblems 
\eqref{Eq:WPF1}, \eqref{Eq:WPF2} and \eqref{Eq:WPM} of fluid flow and mechanical 
deformation by space-time finite element techniques. For the discretization of the time 
variable we consider using continuous and discontinuous finite element methods. For the 
spatial discretization of the subproblem of fluid flow mixed finite element techniques 
are applied. Standard conforming finite element methods are used for the 
spatial discretization of the subproblem of mechanical deformation. The derivation of the 
discrete systems is done briefly here. For the application of space-time finite element 
methods to the subproblems of our iteration scheme and the derivation of their algebraic 
formulations as well as for the construction of appropriate iterative linear solvers and 
preconditioning techniques we refer to \cite{BK15,K15}. 

We decompose the time interval $(0,T]$ into $N$ subintervals $I_n=(t_{n-1},t_n]$, where 
$n\in \{1,\ldots ,N\}$ and $0=t_0<t_1< \cdots < t_{n-1} < t_n = T$ and $\tau = 
\max_{n=1,\ldots N} (t_n-t_{n-1})$. Further we denote by $\mathcal{T}_h=\{K\}$ a finite 
element decomposition of mesh size $h$ of the polyhedral domain $\overline{\Omega}$ into 
closed subsets $K$, quadrilaterals in two dimensions and hexahedrals in three dimensions. 
For the spatial discretization of \eqref{Eq:FSP_1} we use a mixed finite element 
approach. We choose the class of Raviart–Thomas elements for the two-dimensional case and 
the class of Raviart--Thomas--N\'ed\'elec elements in three space dimensions, where $W_h^s 
\subset L^2(\Omega)$ with $\vec W_{h}^s=\{ w_h\in L^2(\Omega) \mid w_h{}_{|_{K}}\circ T_K \in \mathbb{Q}_s\}$ and $\vec V_h^s \subset \vec H(\mathrm{div};\Omega)$ denote the 
corresponding inf-sup stable pair of finite element spaces; cf.\ \cite{BK15,C10,Q94} for the exact definition of $\vec V_h^s$.  Here, $\mathbb Q_s$ is the space of polynomials that are of degree less than or equal to 
$s$ with respect to each variable $x_1,\ldots, x_d$ and $T_K$ is a suitable invertible 
mapping of the reference cube $\widehat K$ to the element $K$ of the triangulation 
$\mathcal T_h$.  For 
the spatial approximation of the displacement field $\vec u$ of \eqref{Eq:FSP_2} we 
discretize the space variables by means of a conforming Galerkin method with finite 
element space $\vec H_{h}^l=\{\vec z_h\in C(\overline \Omega) \mid \vec 
z_h{}_{|_{K}}\circ T_K \in \mathbb{Q}_l^d\,, \; \vec z_h{}_{|\partial \Omega} = \vec 0\}$. The fully discrete space-time finite element spaces of functions that 
are continuous in time are 
then given by 
\begin{align}
\label{Def:FE1}
\mathcal{W}_{\tau,h}^{r,s} & = \{w_{\tau,h}\in C(\overline I;L^2(\Omega))\mid 
w_{\tau,h}{}_{|I_n}\in \mathcal P_r(I_n;W_h^s)\} \,,\\[1ex]
\label{Def:FE2}
\vec{\mathcal{V}}_{\tau,h}^{r,s} & = \{\vec v_{\tau,h}\in C(\overline I;\vec 
H(\mathrm{div};\Omega))\mid 
\vec v_{\tau,h}{}_{|I_n}\in \mathcal P_r(I_n;\vec V_{h}^s)\} \,,\\[1ex]
\label{Def:FE3}
\vec{\mathcal{Z}}_{\tau,h}^{r,l} & = \{\vec z_{\tau,h}\in 
C(\overline I;\vec{H}^1_0(\Omega))\mid \vec z_{\tau,h}{}_{|I_n}\in \mathcal P_r(I_n;\vec 
H_{h}^l)\}\,,
\end{align}
where $\mathcal P_r(I_n;X)$ denotes the space of all polynomials in time up to degree 
$r\geq 0$ on $I_n$ with values in $X$. We choose $l=s+1$ to equilibrate the convergence 
rates of the spatial discretization for the three unknowns $p,\vec q$ and 
$\vec u$; cf.\ \cite[Part I, Thm.\ 5.2]{PW07}. For short, we will also use the abbreviations 
$W_h= W_h^s$, $\vec V_h = \vec V_{h}^s$ and $\vec H_h = \vec H_{h}^{s+1}$ in the sequel. 

Discontinuous counterparts $\widetilde{\mathcal{W}}_{\tau,h}^{r,s}$, 
$\vec{\widetilde{\mathcal{V}}}_{\tau,h}^{r,s}$ and $ 
\vec{\widetilde{\mathcal{H}}}_{\tau,h}^{r,l}$ of the spaces 
\eqref{Def:FE1}--\eqref{Def:FE3}, consisting of functions not necessarily being 
continuous 
in time, are then defined by 
\begin{align}
\label{Def:FE4}
\widetilde{\mathcal{W}}_{\tau,h}^{r,s} & = \{w_{\tau,h}\in L^2( 
I;L^2(\Omega))\mid 
w_{\tau,h}{}_{|I_n}\in \mathcal P_r(I_n;W_h^s)\,, \; w_{\tau,h}(0)\in W_h^s\} \,,\\[1ex]
\label{Def:FE5}
\vec{\widetilde{\mathcal{V}}}_{\tau,h}^{r,s} & = \{\vec v_{\tau,h}\in L^2( 
I;\vec H(\mathrm{div};\Omega))\mid 
\vec v_{\tau,h}{}_{|I_n}\in \mathcal P_r(I_n;\vec V_{h}^s)\,,\; \vec v_{\tau,h}(0)\in 
\vec V_h^s\} \,,\\[1ex]
\label{Def:FE6}
\vec{\widetilde{\mathcal{Z}}}_{\tau,h}^{r,l} & = \{\vec z_{\tau,h}\in 
L^2(I;\vec{H}^1_0(\Omega))\mid \vec z_{\tau,h}{}_{|I_n}\in \mathcal P_r(I_n;\vec 
H_{h}^l)\,,\; \vec z_{\tau,h}(0)\in \vec H_h^l\}\,.
\end{align}

\subsection{The cGP($\vec r$)--MFEM($\vec s$)cG($\vec s$+1) approach.} 
\label{Sec:cGPMFEM}

The space-time finite element approximation of the flow problem \eqref{Eq:WPF1}, 
\eqref{Eq:WPF2} by a continuous finite element approach in time reads as follows: 
{\em Let $\vec u^k_{\tau,h} \in \vec{\mathcal{Z}}_{\tau,h}^{r,s+1}$, $p^k_{\tau,h} \in 
\mathcal{W}_{\tau,h}^{r,s}$ be given and}
\begin{align*}
l^k_{p}(\vec w_{\tau,h}) & = \left\langle f - b \nabla \cdot \partial_t 
\vec u^k_{\tau,h}+ L \partial_t p^{k}_{\tau,h}, w_{\tau,h} \right\rangle \,,
\end{align*}
{\em for $w_{\tau,h} \in  \widetilde{\mathcal{W}}_{\tau,h}^{r-1,s}$. Find 
$p^{k+1}_{\tau,h} \in \mathcal{W}_{\tau,h}^{r,s}$ and $\vec q^{k+1}_{\tau,h}\in 
\vec{\mathcal{V}}_{\tau,h}^{r,s}$ with $p^{k+1}_{\tau,h}(0) = 0$ such that}
\begin{align}
\nonumber 
&\sum_{n=1}^N \Bigg\{\int_{t_{n-1}}^{t_n} \left(\frac{1}{M}+ L \right) \langle 
\partial_t 
p_{\tau,h}^{k+1},w_{\tau,h}\rangle \ud t\\[0.5ex] 
\label{Eq:FDM}
&\qquad \qquad  \qquad + \int_{t_{n-1}}^{t_n}  \langle 
\nabla \cdot \vec q_{\tau,h}^{k+1},w_{\tau,h} 
\rangle \ud t \Bigg\}  = \sum_{n=1}^N \int_{t_{n-1}}^{t_n} l^k_{p}(w_{\tau,h}) \ud 
t\,, \\[1ex]
\label{Eq:FDD}
& \sum_{n=1}^N \Bigg\{\int_{t_{n-1}}^{t_n} \langle \vec{K}^{-1} \vec 
q_{\tau,h}^{k+1}, \vec v_{\tau,h}\rangle \ud t - 
\int_{t_{n-1}}^{t_n}  \langle p^{k+1}_{\tau,h}, \nabla \cdot \vec 
v_{\tau,h}\rangle \ud t\Bigg\} = 0 
\end{align}
{\em for all $w_{\tau,h} \in \widetilde{\mathcal{W}}_{\tau,h}^{r-1,s}$ and $\vec 
s_{\tau,h} \in \vec{\widetilde{\mathcal{V}}}_{\tau,h}^{r-1,s}$.}

The corresponding space-time finite element approximation of the problem \eqref{Eq:WPM} 
of mechanical deformation reads as follows: {\em Let 
$p^{k+1}_{\tau,h} \in \mathcal{W}_{\tau,h}^{r,s}$ be given and}  
\begin{equation*}
\label{Eq:FDMD_0}
l^{k+1}_{\vec u}(\vec z_{\tau,h})  =  b \langle p_{\tau,h}^{k+1}, \nabla \cdot \vec 
z_{\tau,h}\rangle 
\end{equation*}
{\em for $\vec z_{\tau,h} \in \vec{\widetilde{\mathcal{Z}}}_{\tau,h}^{r,s+1}$. Find $\vec 
 u^{k+1}_{\tau,h} \in \vec{\mathcal{Z}}_{\tau,h}^{r,s+1}$ with 
$\vec u^{k+1}_{\tau,h}(0)=\vec 0$ such that}
\begin{equation}
\label{Eq:FDMD}
\begin{split}
\sum_{n=1}^N \Bigg\{\int_{t_{n-1}}^{t_n} 2 \mu \langle \vec \varepsilon(\vec 
u^{k+1}_{\tau,h}) , \vec \varepsilon(\vec z_{\tau,h}) \rangle \ud t & 
+\int_{t_{n-1}}^{t_n} \lambda \langle \nabla \cdot \vec u^{k+1}, \nabla \cdot \vec 
z^{k+1}\rangle \ud t  \Bigg\}\\[1ex] & = \sum_{n=1}^N 
\int_{t_{n-1}}^{t_n} l^{k+1}_{\vec u}(\vec z_{\tau,h}) \ud t
\end{split}
\end{equation}
{\em for all $\vec z_{\tau,h} \in \vec{\widetilde{\mathcal{Z}}}_{\tau,h}^{r-1,s+1}$.}

On the subinterval $\overline I_n$ we expand the discrete functions $p^{k}_{\tau,h} 
\in \mathcal{W}_{\tau,h}^{r,s}$, $\vec q^{k}_{\tau,h}\in 
\vec{\mathcal{V}}_{\tau,h}^{r,s}$ 
and $\vec u^{k}_{\tau,h} \in \vec{\mathcal{H}}_{\tau,h}^{r,s+1}$ in terms of Lagrangian 
basis functions $\varphi_{n,j}$ with respect to $r+1$ nodal points 
$t_{n,j}\in \overline I_n$, $j=0,\ldots , r$, for the time variable such that they admit 
the representations 
\begin{equation}
\label{Eq:RepSolBasis}
p_{\tau,h}^{k}{}_{|I_n} (t) = \sum_{j=0}^r P^{j,k}_{n,h} 
\varphi_{n,j}(t)\,, \; \vec q_{\tau,h}^{k}{}_{|I_n} (t) = \sum_{j=0}^r \vec 
Q^{j,k}_{n,h} \varphi_{n,j}(t)\,, \;\vec u^{k}_{\tau,h}{}_{|I_n} (t) = \sum_{j=0}^r 
\vec U^{j,k}_{n,h} 
\varphi_{n,j}(t)
\end{equation}
for $t\in \overline I_n$ with coefficient functions $P^{j,k}_{n,h} \in W_h$, $\vec 
Q^{j,k}_{n,h} \in \vec V_h$ and $\vec U^{j,k}_{n,h} \in \vec H_h$ for $j=0,\ldots ,r$. 
Then we replace the variational problems \eqref{Eq:FDM}, \eqref{Eq:FDD} and 
\eqref{Eq:FDMD} by the following system of equations: {\em Let $n\in 
\{1,\ldots ,N\}$. Find coefficient functions $P_{n,h}^{i, k+1}\in W_h$ for $i=0,\ldots, 
r$ 
 and $\vec U_{n,h}^{i, k+1} \in \vec H_h$, $\vec Q_{n,h}^{i, k+1}\in \vec V_h$ for 
$i=1,\ldots, r$ such that}
\begin{align}
\nonumber
\dfrac{1}{M} \sum_{j =0}^r \alpha_{ij} \langle P_{n,h}^{j, k+1},  w_h \rangle + L 
\sum_{j =0}^r \alpha_{ij} \langle P_{n,h}^{j, k+1} - P_{n,h}^{j, k},   w_h \rangle  + 
\tau_n \beta_{ii} \langle \nabla \cdot \vec Q_{n,h}^{i, k+1} & , w_h \rangle \\ 
\label{Eq:AlgebProb1}
  = \tau_n \beta_{ii}\langle f(t_{n,i}), w_h\rangle -b \sum_{j =0}^r \alpha_{ij} 
\langle \nabla 
\cdot \vec U_{n,h}^{j, k} & , w_h \rangle\,,
\\[1ex]
\label{Eq:AlgebProb2}
\langle \vec K^{-1} \vec Q_{n,h}^{i, k+1}, \vec v_h \rangle - \langle P_{n,h}^{i, k+1}, 
\nabla \cdot \vec v_h \rangle  & = 0\,,\\[1ex]
\label{Eq:AlgebProb3}
2 \mu \langle \vec \varepsilon(\vec U_{n,h}^{i, k+1}),\vec \varepsilon(\vec z_h)\rangle + 
\lambda \langle \nabla \cdot \vec U_{n,h}^{i, k+1}, \nabla \cdot \vec z_h \rangle - b 
\langle P_{n,h}^{i, k+1}, \nabla \cdot \vec z_h \rangle & = 0
\end{align}
\emph{for all $w_h \in W_h$, $\vec v_h \in \vec V_h$, $\vec z_h \in \vec H_h$ and $i= 
1,\ldots, r$, where $P_{n,h}^{0, k+1}$ is defined by the continuity constraint in 
time of the discrete solution $p^{k+1}_{\tau,h}\in \mathcal{W}_{\tau,h}^{r,s}$, i.e.\ 
$P_{n,h}^{0, 
k+1}=\lim_{l\rightarrow \infty} p^l_{\tau,h}{}_{|I_{n-1}}(t_{n-1})$ for $n>1$ and 
$P_{n,h}^{0,k+1}=0$ for $n=1$.}

The coefficients $\alpha_{ij}$ and $\beta_{ii}$ in 
\eqref{Eq:AlgebProb1}--\eqref{Eq:AlgebProb3} are defined by
\begin{equation*}
 \label{Def:AlphaBeta}
 \alpha_{ij}  = \int_{I_n} \varphi_{n,j}'(t) \cdot \varphi_{n,i}(t)\ud t\,, \quad 
\beta_{ii} = \int_{I_n} \varphi_{n,i} (t) \cdot \varphi_{n,i}(t)\ud t\,, \quad  
i=1,\ldots, r\,, \; j=0,\ldots, r\,.
\end{equation*}
 
\begin{rem}
\begin{itemize}
\item The scheme \eqref{Eq:FDM}--\eqref{Eq:FDMD} defines a Galerkin--Petrov method, 
since the trial spaces \eqref{Def:FE1}--\eqref{Def:FE3} and test spaces 
\eqref{Def:FE4}--\eqref{Def:FE6} differ. 
\item For all technical details of the derivation of the semi-algebraic equations 
\eqref{Eq:AlgebProb1}--\eqref{Eq:AlgebProb3} we refer to, e.g., \cite{BK15,BRK15,K15,S10}.  

\item We note that \eqref{Eq:AlgebProb1}--\eqref{Eq:AlgebProb3} is not the local counterpart of  \eqref{Eq:FDM}--\eqref{Eq:FDMD} on $I_n$, i.e.\ the formulation of \eqref{Eq:AlgebProb1}--\eqref{Eq:AlgebProb3} on the subinterval $I_n$ by a suitable choice of a test basis in time with support in $\overline{I}_n$ (cf.\ \cite{BK15,BRK15,K15,S10}), since in \eqref{Eq:FDM}--\eqref{Eq:FDMD} the iteration process is performed globally on $\overline{I}$. In contrast to this, the scheme \eqref{Eq:AlgebProb1}--\eqref{Eq:AlgebProb3} is based on iterating on each of the subintervals $I_n$ before proceeding to the next one. 

\item For the treatment of the continuity constraint in time we put $t_{n,0}=t_{n-1}$ 
for the nodal points of the Lagrangian basis functions. The other points $t_{n,1},\ldots, 
t_{n,r}$ are chosen as the quadrature points of the $r$-point Gauss quadrature formula on 
$I_n$ which is exact if the function to be integrated is a polynomial of degree less or 
equal to $2r-1$. In particular, there holds that $\varphi_{n,j} (t_{n,i}) = 
\delta_{i,j}$ for $i,j=0,\ldots , r$.

\item  The variational formulations \eqref{Eq:AlgebProb1}--\eqref{Eq:AlgebProb3} solely 
depend on the values of the flux and the displacement variable in the Gauss quadrature 
points as Eqs.\ \eqref{Eq:AlgebProb1} and \eqref{Eq:AlgebProb3} show, i.e.\ they 
depend on $\vec Q_{n,h}^{i, k+1}$ and $\vec U_{n,h}^{i,k+1}$ for $i=1,\ldots ,r$. We then 
define the flux and the displacement variable in the grid points by extrapolation, in 
this way also ensuring the continuity in time, i.e.\ $\vec Q_{n,h}^{0,k+1}=\vec 
q_{\tau,h}{}_{|I_{n-1}}(t_{n-1})$ and $\vec U_{n,h}^{0,k+1}=\vec 
u_{\tau,h}{}_{|I_{n-1}}(t_{n-1})$; cf.\ 
\cite{BK15,BRK15,K15,S10}.

\item We define the discrete initial flux as a suitable finite element approximation in 
$\vec V_h$ of $\vec q(0) = - \vec K\nabla p_0$, if $p_0$ is sufficiently regular. If 
this is not the case we take a regular approximation. The discrete initial flux is only 
needed for having a consistent notation and the extrapolation argument of the previous 
item in the first subinterval $I_1$. The discrete initial flux is of no relevance for the 
analysis of the scheme. 
\end{itemize} 
\end{rem}

\subsection{The dG($\vec r$)--MFEM($\vec s$)cG($\vec s$+1) approach.} 
\label{Sec:dGMFEM}

The space-time finite element approximation of the flow problem \eqref{Eq:WPF1}, 
\eqref{Eq:WPF2} by a discontinuous finite element approach in time (cf.\ 
\cite{DF15,T06,BK15,K14}) reads as follows: 
{\em Let $\vec u^k_{\tau,h} \in \vec{\widetilde{\mathcal{Z}}}_{\tau,h}^{r,s+1}$, 
$p^k_{\tau,h} \in \widetilde{\mathcal{W}}_{\tau,h}^{r,s}$ be given and}
\begin{align*}
l^k_{p}(\vec w_{\tau,h}) & = \left\langle f - b \nabla \cdot \partial_t 
\vec u^k_{\tau,h}+ L \partial_t p^{k}_{\tau,h}, w_{\tau,h} \right\rangle \,,
\end{align*}
{\em for $w_{\tau,h} \in  \widetilde{\mathcal{W}}_{\tau,h}^{r,s}$. Find 
$p^{k+1}_{\tau,h} \in \widetilde{\mathcal{W}}_{\tau,h}^{r,s}$ and 
$\vec q^{k+1}_{\tau,h}\in 
\vec{\widetilde{\mathcal{V}}}_{\tau,h}^{r,s}$ with $p^{k+1}_{\tau,h} (0) = 0$ such that}

\begin{align}
\nonumber 
&\sum_{n=1}^N \Bigg\{\int_{t_{n-1}}^{t_n} \left(\frac{1}{M}+ L \right) \langle 
\partial_t 
p_{\tau,h}^{k+1},w_{\tau,h}\rangle \ud t + \int_{t_{n-1}}^{t_n}  \langle 
\nabla \cdot q_{\tau,h}^{k+1},w_{\tau,h} 
\rangle \ud t \Bigg\}\\[0.5ex] 
\label{Eq:D_FDM}
&\qquad +\Bigg(\frac{1}{M}+L\Bigg)\left\langle 
\left[p^{k+1}_{\tau,h}\right]_{n-1},w_{\tau,h}(t_{n-1} ^+)\right\rangle = \sum_{n=1}^N 
\int_{t_{n-1}}^{t_n} l^k_{p}(w_{\tau,h}) \ud 
t \\[1ex]
\nonumber
& \qquad + L \left\langle 
\left[p^{k}_{\tau,h}\right]_{n-1},w_{\tau,h}(t_{n-1} ^+)\right\rangle - \left\langle 
\left[\nabla \cdot \vec u_{\tau,h}^k\right]_{n-1},w_{\tau,h}(t_{n-1} ^+)\right\rangle\,,
\\[1ex]
\label{Eq:D_FDD}
& \sum_{n=1}^N \Bigg\{\int_{t_{n-1}}^{t_n} \langle \vec{K}^{-1} \vec 
q_{\tau,h}^{k+1}, \vec v_{\tau,h}\rangle \ud t - 
\int_{t_{n-1}}^{t_n}  \langle p^{k+1}_{\tau,h}, \nabla \cdot \vec 
v_{\tau,h}\rangle \ud t\Bigg\} = 0 
\end{align}
{\em for all $w_{\tau,h} \in \widetilde{\mathcal{W}}_{\tau,h}^{r,s}$ and $
\vec v_{\tau,h} \in \vec{\widetilde{\mathcal{V}}}_{\tau,h}^{r,s}$.}

Here we use the notation 
\begin{equation*}
p_{\tau,h}^k (t_n^-) = \lim_{t\rightarrow t_n-0} p^k_{\tau,h}(t) \,, \quad 
p_{\tau,h}^k (t_n^+) = \lim_{t\rightarrow t_n+0} p^k_{\tau,h}(t) \,, \quad 
[p_{\tau,h}^k]_n = p^k_{\tau,h} (t_n^+) - p^k_{\tau,h} (t_n^-)\,,
\end{equation*}
and analogously for the displacement field $\vec u^k_{\tau,h}$.

The corresponding space-time finite element approximation of the problem \eqref{Eq:WPM} 
of mechanical deformation reads as follows: {\em Let 
$p^{k+1}_{\tau,h} \in \widetilde{\mathcal{W}}_{\tau,h}^{r,s}$ be given and}  
\begin{equation*}
\label{Eq:D_FDMD_0}
l^{k+1}_{\vec u}(\vec z_{\tau,h})  =  b \langle p_{\tau,h}^{k+1} , \nabla \cdot 
\vec z_{\tau,h}\rangle 
\end{equation*}
{\em for $\vec z_{\tau,h} \in \vec{\widetilde{\mathcal{Z}}}_{\tau,h}^{r,s+1}$. Find $\vec 
 u^{k+1}_{\tau,h} \in \vec{\widetilde{\mathcal{Z}}}_{\tau,h}^{r,s+1}$ with 
$\vec u^{k+1}_{\tau,h} (0) = \vec 0$  such that}
\begin{equation}
\label{Eq:D_FDMD}
\begin{split}
\sum_{n=1}^N \Bigg\{\int_{t_{n-1}}^{t_n} 2 \mu \langle \vec \varepsilon(\vec 
u^{k+1}_{\tau,h}) , \vec \varepsilon(\vec z_{\tau,h}) \rangle \ud t & 
+\int_{t_{n-1}}^{t_n} \lambda \langle \nabla \cdot \vec u_{\tau,h}^{k+1}, \nabla \cdot 
\vec z_{\tau,h}^{k+1} \rangle
\ud t  \Bigg\}\\[1ex] & = \sum_{n=1}^N 
\int_{t_{n-1}}^{t_n} l^{k+1}_{\vec u}(\vec z_{\tau,h}) \ud t
\end{split}
\end{equation}
{\em for all $\vec z_{\tau,h} \in \vec{\widetilde{\mathcal{Z}}}_{\tau,h}^{r,s+1}$.}

On $ I_n$ we expand the discrete functions $p^{k}_{\tau,h} \in 
\widetilde{\mathcal{W}}_{\tau,h}^{r,s}$, $\vec q^{k}_{\tau,h}\in 
\vec{\widetilde{\mathcal{V}}}_{\tau,h}^{r,s}$ 
and $\vec u^{k}_{\tau,h} \in \vec{\widetilde{\mathcal{Z}}}_{\tau,h}^{r,s+1}$ in time in 
terms of 
Lagrangian basis functions $\varphi_{n,j}$ with respect to $r+1$ nodal points 
$t_{n,j}\in I_n$, \begin{equation}
\label{Eq:dGRepSolBasis}
p_{\tau,h}^{k}{}_{|I_n} (t) = \sum_{j=0}^r P^{j,k}_{n,h} 
\varphi_{n,j}(t)\,, \; \vec q_{\tau,h}^{k}{}_{|I_n} (t) = \sum_{j=0}^r \vec 
Q^{j,k}_{n,h} \varphi_{n,j}(t)\,, \;\vec u^{k}_{\tau,h}{}_{|I_n} (t) = \sum_{j=0}^r 
\vec U^{j,k}_{n,h} 
\varphi_{n,j}(t)
\end{equation}
for $t\in I_n$ with coefficient functions $P^{j,k}_{n,h} \in W_h$, $\vec 
Q^{j,k}_{n,h} \in \vec V_h$ and $\vec U^{j,k}_{n,h} \in \vec H_h$ for $j=0,\ldots ,r$. 
The nodal points $t_{n,j}$, with $j=0,\ldots,r$,  are chosen as the quadrature points of the $r$+1-point Gauss 
quadrature formula on $I_n$ which is exact for polynomials of degree less or equal to 
$2r+1$.

Then we replace the variational problems \eqref{Eq:D_FDM}, \eqref{Eq:D_FDD} and 
\eqref{Eq:D_FDMD} by the following system of equations: {\em Let $n\in 
\{1,\ldots ,N\}$. Find coefficient functions $P_{n,h}^{i, k+1}\in W_h$, $\vec U_{n,h}^{i, 
k+1} \in \vec H_h$ and $\vec Q_{n,h}^{i, k+1}\in \vec V_h$ for $i=0,\ldots, r$ such that}
\begin{align}
\nonumber
\dfrac{1}{M} \sum_{j =0}^r \widetilde\alpha_{ij} \langle P_{n,h}^{j, k+1},  w_h \rangle + 
L 
\sum_{j =0}^r \widetilde\alpha_{ij} \langle P_{n,h}^{j, k+1} - P_{n,h}^{j, k},   w_h 
\rangle  + 
\tau_n \widetilde\beta_{ii} \langle \nabla \cdot \vec Q_{n,h}^{i, k+1} & , w_h \rangle \\ 
\label{Eq:D_AlgebProb1}
= \tau_n \widetilde\beta_{ii}\langle f(t_{n,i}), w_h\rangle -b \sum_{j =0}^r 
\alpha_{ij} \langle \nabla 
\cdot \vec U_{n,h}^{j, k} & , w_h \rangle \\
\nonumber
+ \gamma_i \, \frac{1}{M}\langle p^\infty_{\tau,h}(t_{n-1}^-),w_h\rangle + \gamma_i 
\langle b \nabla \cdot \vec u^\infty_{\tau,h}(t_{n-1}^-)& ,w_h\rangle\,,\\[1ex]
\label{Eq:D_AlgebProb2}
\langle \vec K^{-1} \vec Q_{n,h}^{i, k+1}, \vec v_h \rangle - \langle P_{n,h}^{i, k+1}, 
\nabla \cdot \vec v_h \rangle  & = 0\,,\\[1ex]
\label{Eq:D_AlgebProb3}
2 \mu \langle \vec \varepsilon(\vec U_{n,h}^{i, k+1}),\vec \varepsilon(\vec z_h)\rangle + 
\lambda \langle \nabla \cdot \vec U_{n,h}^{i, k+1}, \nabla \cdot \vec z_h \rangle - b 
\langle P_{n,h}^{i, k+1}, \nabla \cdot \vec z_h \rangle & = 0
\end{align}
\emph{for all $w_h \in W_h$, $\vec v_h \in \vec V_h$, $\vec z_h \in \vec H_h$ and $i= 
0,\ldots, r$, where $p^\infty_{\tau,h}(t_{n-1}^-) = \lim_{l\rightarrow 
\infty}p^l_{\tau,h}(t_{n-1}^-)$ and $\vec u^\infty_{\tau,h}(t_{n-1}^-) = 
\lim_{l\rightarrow 
\infty} \vec u^l_{\tau,h}(t_{n-1}^-)$ for $n>1$ as well as  $p_{\tau,h}(t_{n-1}^-)=0$ and 
$\vec u_{\tau,h}(t_{n-1}^-)=\vec 0$ for $n=1$.}

The coefficients $\widetilde\alpha_{ij}$, $\widetilde\beta_{ii}$, $\alpha_{ij}$ and 
$\gamma_i$ are defined by 
\begin{equation*}
\widetilde\alpha_{ij} = \alpha_{ij}+ \gamma_i\cdot\gamma_j\,, \qquad 
\widetilde\beta_{ii} = \beta_{ii}\,,\qquad \gamma_{i} = \varphi_{n,i}(t_{n-1}^+)
\end{equation*}
with 
\begin{equation*}
\alpha_{ij} = \int_{I_n} \varphi_{n,j}'(t) \cdot \varphi_{n,i}(t)\ud t\,, \quad 
\beta_{ii} = \int_{I_n} \varphi_{n,i} (t) \cdot \varphi_{n,i}(t)\ud t\
\end{equation*}
for $i,j=0,\ldots,r$.

\section{Convergence of the iteration schemes}
\label{Sec:Convergence}

Now we prove the convergence of the iterative splitting schemes that we introduced in 
Sec.\ 
\ref{Sec:Discretization}.

\subsection{The cGP($\vec r$)--MFEM($\vec s$)cG($\vec s$+1) approach.} 
\label{Subsec::cGConv}

In this subsection we prove the (linear) convergence of the splitting schemes 
\eqref{Eq:AlgebProb1}--\eqref{Eq:AlgebProb3} based on a continuous Galerkin 
discretization of the time variable. For this we show that the scheme is 
subject to a contraction principle such that a unique fixed point is obtained. This  
convergence is proved in strong energy norms.

In the sequel, we denote by $p_{\tau,h} \in \mathcal{W}_{\tau,h}^{r,s}$, 
$\vec q_{\tau,h}\in \vec{\mathcal{V}}_{\tau,h}^{r,s}$ and $\vec u_{\tau,h} \in 
\vec{\mathcal{Z}}_{\tau,h}^{r,s+1}$, with 
\begin{equation}
\label{EQ:RepFDS}
p_{\tau,h}{}_{|I_n} (t) = \sum_{j=0}^r P^{j}_{n,h} \varphi_{n,j}(t)\,, \; \vec 
q_{\tau,h}{}_{|I_n} (t) = \sum_{j=0}^r \vec Q^{j}_{n,h} \varphi_{n,j}(t)\,, \;\vec 
u_{\tau,h}{}_{|I_n} (t) = \sum_{j=0}^r 
\vec U^{j}_{n,h} 
\varphi_{n,j}(t)
\end{equation}
for $t\in \overline I_n$, the space-time finite element approximation of the Biot system 
\eqref{Eq:B_4}--\eqref{Eq:B_7} that is defined by skipping the upper indices in the 
problems \eqref{Eq:FDM}, \eqref{Eq:FDD} and \eqref{Eq:FDMD}, respectively. Thus we 
tacitly suppose that the coupled system that is obtained by discretizing the Biot 
model \eqref{Eq:B_4}--\eqref{Eq:B_7} in the space-time finite element spaces 
\eqref{Def:FE1}--\eqref{Def:FE3} admits a unique solution. By means of our variational 
framework for the time discretization the existence and uniqueness of the solution can be 
shown along the lines of \cite[Part I, Sec.\ 4]{PW07}, where the proof is given 
for the spatially semidiscretized problem.  

For the sake of brevity, we define the following variables quantifying the errors between 
this  space-time finite element approximation of the Biot system 
\eqref{Eq:B_4}--\eqref{Eq:B_7} and its approximation after $k$ iterations of the 
proposed scheme \eqref{Eq:AlgebProb1}--\eqref{Eq:AlgebProb3}. For fixed $n\in 
\{1,\ldots ,N\}$ we put  
\begin{equation*}
\begin{aligned}
E_p^{j, k} & = P_{n, h}^{j, k} -  P_{n, h}^{j}\,, && j \in \{0,\ldots,r\}\,, &&& 
e_p^{k}(t) 
&= 
\sum_{j=0}^r E_p^{j, k}\varphi_{n,j}(t) \,, && t\in \overline I_n\,, \\[0ex]
\vec E_{\vec q}^{j, k} &= \vec Q_{n, h}^{j, k} -  \vec Q_{n, h}^{j}\,, && j \in 
\{0,\ldots,r\} \,, &&& \vec e_{\vec q}^{k}(t) &= \sum_{j=0}^r \vec E_{\vec q}^{j, 
k}\varphi_{n,j}(t)\,, && t\in \overline I_n\,,\\[0ex]
\vec E_{\vec u}^{j, k} &= \vec U_{n, h}^{j, k} -  \vec U_{n, h}^{j}\,, && j \in 
\{0,\ldots,r\} \,, &&& \vec e_{\vec u}^{k}(t) &= \sum_{j=0}^r \vec E_{\vec u}^{j, 
k}\varphi_{n,j}(t)\,, && t\in \overline I_n\,.
\end{aligned}
\end{equation*}

In order to simplify the notation below, we further introduce the abbreviations
\begin{equation}
 \label{Eq:DefS}
S_p^{i, k+1} = \sum_{j =0}^r \alpha_{ij} E_p^{j, k+1}\,, \quad
\vec S_{\vec q}^{i, k+1} =  \sum_{j =0}^r \alpha_{ij} \vec E_{\vec q}^{j, k+1}\,,\quad 
\vec S_{\vec u}^{i, k+1} =  \sum_{j =0}^r \alpha_{ij} \vec E_{\vec u}^{j, k+1}
\end{equation}
with $S_p^{i, k+1} \in W_h$, $\vec S_{\vec q}^{i, k+1} \in \vec V_h$ and $\vec 
S_{\vec u}^{i, k+1} \in \vec  H_h$ for $i=1,\ldots,r$. 

\begin{rem} Due to the continuity constraint in time that is imposed 
by the definition of the space-time finite element spaces 
\eqref{Def:FE1}--\eqref{Def:FE3} and incorporated into the scheme \eqref{Eq:AlgebProb1}--\eqref{Eq:AlgebProb3}  there holds that 
\begin{equation}
\label{Eq:Ek_eq_0}
E_p^{0,k} = 0\,, \qquad \vec E_{\vec q}^{0,k} = \vec {0}\,, \qquad \vec E_{\vec u}^{0,k} 
= \vec {0} 
\end{equation}
for any iteration index $k \in \N$.
\end{rem}

\begin{theorem} 
\label{Thm:ConvDisc}
Let $p_{\tau,h} \in \mathcal{W}_{\tau,h}^{r,s}$, $\vec q_{\tau,h}\in 
\vec{\mathcal{V}}_{\tau,h}^{r,s}$ and $\vec u_{\tau,h} \in 
\vec{\mathcal{Z}}_{\tau,h}^{r,s+1}$ denote the 
fully discrete space-time finite element approximation of the Biot system 
\eqref{Eq:B_4}--\eqref{Eq:B_7}. On $I_n$ let $\{p_{\tau,h} ,\vec q_{\tau,h},\vec u_{\tau,h}\}$ be represented by \eqref{EQ:RepFDS} and let $\{p^k_{\tau,h} ,\vec 
q^k_{\tau,h}, \vec u^k_{\tau,h} \}$ be defined  by \eqref{Eq:RepSolBasis} with coefficient functions being given by the scheme \eqref{Eq:AlgebProb1}--\eqref{Eq:AlgebProb3}. 
Then, for any $L \geq b^2/(2\lambda)$ the sequence 
$\{S_{p}^{i,k}\}_k$, for $i=1,\ldots,r$, converges geometrically in $W_h$. For 
$n=1,\ldots, N$ this implies the convergence of 
$\{p^k_{\tau,h}(t_n),\vec q^k_{\tau,h}(t_n),\vec u^k_{\tau,h}(t_n)\}$ to 
$ \{p_{\tau,h}(t_n),\vec q_{\tau,h}(t_n),\vec u_{\tau,h}(t_n)\}$ in $W_h\times \vec 
V_h \times \vec H_h$ for $k\rightarrow \infty$. 
\end{theorem}

\begin{mproof}
We split the proof into several steps.

\medskip
\noindent \textbf{1.\ Step (Error equations).} By substracting equations 
\eqref{Eq:AlgebProb1}--\eqref{Eq:AlgebProb3} from the system that is obtained by 
discretizing the coupled Biot model \eqref{Eq:B_4}--\eqref{Eq:B_7} in the space-time 
finite element spaces \eqref{Def:FE1}--\eqref{Def:FE3}, respectively, we obtain the error 
equations  
\begin{align}
\dfrac{1}{M} \sum_{j =0}^r \alpha_{ij} \la E_p^{j, k+1}, w_h \ra + L \sum_{j =0}^r 
\alpha_{ij} \la E_p^{j, k+1} - E_p^{j, k}, w_h \ra \qquad \qquad \qquad \nonumber \\
  + \tau_n \beta_{ii} \la \nabla \cdot \vec E_{\vec q}^{i, k+1}, w_h \ra = -b \sum_{j 
=0}^r \alpha_{ij} \la \nabla \cdot \vec E_{\vec u}^{j, k}, w_h \ra \,,
\label{proof_eq_7}\\[2ex]
 \la \vec K^{-1} \vec E_{\vec q}^{i, k+1}, \vec v_h \ra - \la E_p^{i, k+1}, \nabla \cdot 
\vec v_h \ra = 0\,,
\label{proof_eq_8} \\[2ex]
 2 \mu \la \vec \varepsilon(\vec E_{\vec u}^{i, k+1}), \vec \veps(\vec z_h)\ra + \lambda 
\la \nabla \cdot \vec E_{\vec u}^{i, k+1},   \nabla \cdot \vec z_h\ra - b \la E_p^{i, 
k+1}, \nabla \cdot \vec z_h \ra= 0 \phantom{\,,}
\label{proof_eq_9}
\end{align}
for all $w_h \in W_h$, $\vec v_h \in \vec V_h$, $\vec z_h \in \vec H_h$ and $i=1, 
\ldots, r$.

In the next steps we choose appropriate test functions in the Eqs.\ 
\eqref{proof_eq_7}--\eqref{proof_eq_9}, respectively, and sum up resulting identities.   

\medskip
\noindent\textbf{2.\ Step (Choice of test function in Eq.\ \eqref{proof_eq_7}).} 
We test Eq.~\eqref{proof_eq_7} with $w_h = \sum_{j =0}^r \alpha_{ij} E_p^{j, k+1}$ 
to get 
that 
\begin{equation}
\label{proof_eq_10}
\begin{aligned}
\dfrac{1}{M} \bigg\| \sum_{j =0}^r & \alpha_{ij} E_p^{j, k+1} \bigg\|^2  + L 
\bigg\langle \sum_{j =0}^r \alpha_{ij} (E_p^{j, k+1}-E_p^{j, k}), \sum_{j =0}^r 
\alpha_{ij} E_p^{j, k+1} \bigg\rangle \\[1ex]
& +\tau_n \beta_{ii} \bigg\langle \nabla \cdot \vec E_{\vec q}^{i, k+1}, \sum_{j =0}^r 
\alpha_{ij} E_p^{j, k+1}\bigg\rangle = - b \bigg\langle \sum_{j =0}^r \alpha_{ij} \nabla 
\cdot \vec E_{\vec u}^{j, k}, \sum_{j =0}^r \alpha_{ij} E_p^{j, k+1} \bigg\rangle
\end{aligned}
\end{equation}
for any $i \in \{1, \ldots, r\}$. Using the notation \eqref{Eq:DefS}, we can rewrite 
Eq.\ \eqref{proof_eq_10} as 
\begin{equation}
\label{Eq:proof_eq_10.1}
\begin{split}
\dfrac{1}{M} \| S_p^{i, k+1} \|^2 + L \la S_p^{i, k+1} - S_p^{i, k}, S_p^{i, k+1} \ra  + 
\tau_n \beta_{ii} & \la \nabla \cdot \vec E_{\vec q}^{i, k+1}, S_p^{i, k+1} \ra\\[1ex]
& = - b \la \nabla \cdot \vec S_{\vec u}^{i, k}, S_p^{i, k+1} \ra\,.
\end{split}
\end{equation}
We note that $\beta_{ii} > 0 $ for $i=0,\ldots, r$; cf.\ \cite[Lemma 2.2]{BRK15}. Now, 
dividing Eq.\ \eqref{Eq:proof_eq_10.1} by $\beta_{ii}$ and using the algebraic 
identity 
\begin{equation*}
\la x-y, x \ra = \dfrac{1}{2} \| x \|^2 + \dfrac{1}{2} \| x -y \|^2 - \dfrac{1}{2} \| y 
\|^2\, 
\end{equation*}
we recover Eq.\ \eqref{Eq:proof_eq_10.1} in the equivalent form that 
\begin{equation}
\label{proof_eq_11}
\begin{aligned}
\bigg(\dfrac{1}{M \beta_{ii}}  + \dfrac{L}{ 2 \beta_{ii}}\bigg)  \| S_p^{i, k+1} \|^2 & + 
\dfrac{L}{ 2 \beta_{ii}}  \| S_p^{i, k+1} - S_p^{i, k} \|^2 - \dfrac{L}{ 2 \beta_{ii}} 
\| 
S_p^{i, k} \|^2 \\[1ex]
& + \tau_n \la \nabla \cdot \vec E_{\vec q}^{i, k+1},  S_p^{i, k+1} \ra = - 
\dfrac{b}{\beta_{ii}} \la \nabla \cdot \vec S_{\vec u}^{i, k}, S_p^{i, k+1} \ra
\end{aligned}
\end{equation}
for $i=1,\ldots,r$. 

\medskip
\noindent\textbf{3.\ Step (Summation of Eq.\ \eqref{proof_eq_8} and choice of test 
function).} Firstly, we note that Eq.\ \eqref{proof_eq_8} is also satisfied for $i=0$ by 
means of the observation \eqref{Eq:Ek_eq_0}. Changing the index $i$ in Eq.\  
\eqref{proof_eq_8} to $j$, multiplying the resulting equation with $\alpha_{ij}$ and, 
then, summing up from $j=0$ to $r$ and recalling Eq.\ \eqref{Eq:DefS} yields that 
\begin{equation}
\la \vec K^{-1} \vec S_{\vec q}^{i, k+1}, \vec v_h \ra - \la S_p^{i, k+1}, \nabla \cdot 
\vec v_h \ra = 0
\label{proof_eq_12}
\end{equation}
for all $\vec v_h \in \vec V_h$ and $i=1,\ldots, r$. Testing Eq.\ \eqref{proof_eq_12} 
with $\vec v_h =\tau_n \vec E_{\vec q}^{i, k+1} \in \vec V_h$ we get that 
\begin{equation}\label{proof_eq_14}
\tau_n \la \vec K^{-1} \vec S_{\vec q}^{i, k+1}, \vec E_{\vec q}^{i, k+1} \ra - \tau_n 
\la S_p^{i, k+1}, \nabla \cdot \vec E_{\vec q}^{i, k+1}\ra = 0\,.
\end{equation}
Adding Eq.\ 
\eqref{proof_eq_14} to Eq.\ \eqref{proof_eq_11} then gives that 
\begin{equation}
\label{proof_eq_15}
\begin{aligned}
\bigg(\dfrac{1}{M \beta_{ii}} + \dfrac{L}{ 2 \beta_{ii}}\bigg)\| S_p^{i, k+1} \|^2 & + 
\dfrac{L}{ 2 \beta_{ii}} \| S_p^{i, k+1} - S_p^{i, k} \|^2 - \dfrac{L}{ 2 \beta_{ii}} \| 
S_p^{i, k} 
\|^2 \\[1ex]
& + \tau_n \la \vec K^{-1} \vec S_{\vec q}^{i, k+1}, \vec E_{\vec q}^{i, k+1}\ra = - 
\dfrac{b}{\beta_{ii}} \la \nabla \cdot \vec S_{\vec u}^{i, k}, S_p^{i, k+1} \ra 
\end{aligned}
\end{equation}
for all $i = 1, \ldots, r$.

\medskip
\noindent\textbf{4.\ Step (Summation of Eq.\ \eqref{proof_eq_9} and choice of test 
function).} Similarly, we note that Eq.\ \eqref{proof_eq_9} is also satisfied for $i=0$ 
by means of the observation \eqref{Eq:Ek_eq_0}. Changing the index $i$ in Eq.\ 
\eqref{proof_eq_9} to $j$, multiplying the resulting equation with $\alpha_{ij}$ and, 
then, summing up from $j=0$ to $r$ and recalling Eq.\ \eqref{Eq:DefS} yields that 
\begin{equation}
2 \mu \la \vec \veps(\vec S_{\vec u}^{i, k+1}), \vec \veps(\vec z_h)\ra + \lambda \la 
\nabla \cdot \vec S_{\vec u}^{i, k+1},   \nabla \cdot \vec z_h\ra - b \la S_p^{i, 
k+1}, \nabla \cdot \vec z_h \ra = 0
\label{proof_eq_13}
\end{equation}
for all $\vec z_h \in \vec H_h$ and $i=1,\ldots, r$. Testing Eq.\ \eqref{proof_eq_13} 
with $\vec z_h = \dfrac{1}{\beta_{ii}} \vec S_{\vec u}^{i, k} \in \vec H_h$ yields that
\begin{equation}
\label{proof_eq_16}
\dfrac{2 \mu}{\beta_{ii}} \la \vec \veps(\vec S_{\vec u}^{i, k+1}), \vec 
\veps(\vec S_{\vec u}^{i, k})\ra + \dfrac{\lambda}{\beta_{ii}} \la \nabla \cdot 
\vec S_{\vec u}^{i, k+1},   \nabla \cdot  \vec S_{\vec u}^{i, 
k}\ra - \dfrac{b}{\beta_{ii}}  \la S_p^{i, k+1}, \nabla \cdot \vec S_{\vec u}^{i, k} \ra 
= 0
\end{equation}
for $i = 1, \ldots, r$, where we again used that $\beta_{ii}>0$; cf.\ \cite[Lemma 
2.2]{BRK15}. Adding Eq.\ \eqref{proof_eq_16} to Eq.\ \eqref{proof_eq_15} 
leads to
\begin{equation}\label{proof_eq_17}
\begin{split}
\bigg(\dfrac{1}{M \beta_{ii}}  + \dfrac{L}{ 2 \beta_{ii}} & \bigg)\| S_p^{i, k+1} \|^2 + 
\dfrac{L}{ 2 \beta_{ii}} \| S_p^{i, k+1} - S_p^{i, k} \|^2 + \tau_n \la \vec K^{-1} 
\vec S_{\vec q}^{i, k+1}, \vec E_{\vec q}^{i, k+1}\ra \\[1ex]
& +\dfrac{2 \mu}{\beta_{ii}} \la \vec \veps(\vec S_{\vec u}^{i, k+1}), 
\vec \veps(\vec S_{\vec u}^{i, k})\ra 
+ \dfrac{\lambda}{\beta_{ii}} \la \nabla \cdot \vec S_{\vec u}^{i, k+1},   \nabla \cdot  
\vec S_{\vec u}^{i, k}\ra =  \dfrac{L}{ 2 \beta_{ii}} \| S_p^{i, k} \|^2
\end{split}
\end{equation}
for $i = 1, \ldots, r$. 

In the next step we consider the resulting incremental equation that is obtained by 
substracting Eq.\ \eqref{proof_eq_13} written for two consecutive iteration indices from 
each other.

\medskip
\noindent\textbf{5.\ Step (Formation of incremental equation for \eqref{proof_eq_13}, 
choice of test function and summation).} 
We return to Eq.\ \eqref{proof_eq_13}, write it for two consecutive iterations, $k$ and 
$k +1$, and substract the resulting equations from each other to obtain that 
\begin{equation}\label{proof_eq_17b}
2 \mu \la \vec \veps(\vec S_{\vec u}^{i, k+1} - \vec S_{\vec u}^{i, k}), 
\vec \veps(\vec z_h)\ra + \lambda \la \nabla 
\cdot (\vec S_{\vec u}^{i, k+1} - \vec S_{\vec u}^{i, k}),   \nabla \cdot \vec z_h \ra - 
b 
\la S_p^{i, k+1} - S_p^{i, k}, \nabla \cdot \vec z_h  \ra = 0
\end{equation}
for all $\vec z_h \in \vec H_h$ and $i = 1,\ldots ,r$. Choosing $ \vec z_h = \vec 
S_{\vec u}^{i, k+1} - \vec S_{\vec u}^{i, k} \in \vec H_h $ in Eq.\ \eqref{proof_eq_17b}, 
we find that 
\begin{equation}\label{proof_eq_18}
2 \mu \| \vec \veps(\vec S_{\vec u}^{i, k+1} - \vec S_{\vec u}^{i, k}) \|^2 + \lambda \| 
\nabla \cdot (\vec S_{\vec u}^{i,k+1} - \vec S_{\vec u}^{i, k})\|^2 =  b \la S_p^{i, k+1} 
- S_p^{i, k}, \nabla \cdot (\vec S_{\vec u}^{i, k+1} 
- \vec S_{\vec u}^{i, k}) \ra 
\end{equation}
for $i = 1,\ldots ,r$. By dividing Eq.\ \eqref{proof_eq_18} by $ \beta_{ii} > 0$ and 
summing up the resulting identity from $i = 1 $ to $r$ we obtain that 
\begin{equation}\label{proof_eq_20}
\begin{split}
\sum_{i =1}^r \dfrac{2 \mu}{\beta_{ii}} \| \vec \veps(\vec S_{\vec u}^{i, k+1} - \vec 
S_{\vec u}^{i, k}) \|^2 & + \sum_{i =1}^r \dfrac{\lambda}{\beta_{ii}}  \| \nabla \cdot 
(\vec S_{\vec u}^{i, k+1} - \vec S_{\vec u}^{i, 
k})\|^2\\[1ex]
& =  \sum_{i =1}^r \dfrac{b}{\beta_{ii}} \la S_p^{i, k+1} - S_p^{i, k}, 
\nabla \cdot (\vec S_{\vec u}^{i, k+1} - \vec S_{\vec u}^{i, k}) \ra\,.
\end{split}
\end{equation}
Further, from Eq.\ \eqref{proof_eq_18} we get by means of the inequality 
of Cauchy-Schwarz that 
\begin{equation}\label{proof_eq_19}
\lambda \| \nabla \cdot (\vec S_{\vec u}^{i, k+1} - \vec S_{\vec u}^{i, k})\| \le b \| 
S_p^{i, k+1} - S_p^{i,k} \|
\end{equation}
for $i = 1,\ldots ,r$. 

Next, we combine the derived relations.

\medskip
\textbf{6.\ Step (Summation of Eq.\ \eqref{proof_eq_17} over $\vec i$ and combination 
with derived relations).} 
Using the algebraic identity
\begin{displaymath}
\la x, y \ra = \dfrac{1}{4} \| x + y \|^2 -  \dfrac{1}{4} \| x - y \|^2,
\end{displaymath}
we get from Eq.\ \eqref{proof_eq_17} that
\begin{equation}\label{proof_eq_21}
\begin{split}
\bigg(\dfrac{1}{M \beta_{ii}} & + \dfrac{L}{ 2 \beta_{ii}}\bigg)\| S_p^{i, k+1} \|^2 + 
\dfrac{L}{ 2 \beta_{ii}} \| S_p^{i, k+1} - S_p^{i, k} \|^2 + \tau_n \la \vec K^{-1} 
\vec S_{\vec q}^{i, k+1}, \vec E_{\vec q}^{i, k+1}\ra \\[1ex]
& +\dfrac{\mu}{ 2 \beta_{ii}} \|  \veps(\vec S_{\vec u}^{i, k+1}+ \vec S_{\vec u}^{i, 
k})\|^2 + \dfrac{\lambda}{ 4 \beta_{ii}}\|  \nabla \cdot (\vec S_{\vec u}^{i, k+1} + 
\vec S_{\vec u}^{i, k})\|^2 \\[1ex]
& - \dfrac{\mu}{ 2 \beta_{ii}} \|  \veps(\vec S_{\vec u}^{i, k+1}- \vec S_{\vec u}^{i, 
k})\|^2 - \dfrac{\lambda}{ 4 \beta_{ii}}\|  \nabla \cdot (\vec S_{\vec u}^{i, k+1} - 
\vec S_{\vec u}^{i, k})\|^2 =  \dfrac{L}{ 2 \beta_{ii}} \| S_p^{i, k} \|^2
\end{split}
\end{equation}
for $i = 1, \ldots, r$. Summing up Eq.\ \eqref{proof_eq_21} from $i = 1$ to $r$ and  
using the relations \eqref{proof_eq_20} and  \eqref{proof_eq_19}, we find that 
\begin{equation}\label{proof_eq_22}
\begin{split}
\sum_{i =1}^r & \bigg\{ \bigg(\dfrac{1}{M \beta_{ii}} + \dfrac{L}{ 
2 \beta_{ii}}\bigg) \| S_p^{i, k+1} \|^2  \\[1ex]
& \quad + \dfrac{L}{ 2 \beta_{ii}} \| S_p^{i, k+1} - 
S_p^{i, k} \|^2 + \tau_n 
\la \vec K^{-1} \vec S_{\vec q}^{i, k+1}, \vec E_{\vec q}^{i, k+1}\ra\bigg\} \\[1ex]
& \quad + \sum_{i =1}^r \bigg \{ \dfrac{\mu}{ 2 \beta_{ii}} \|  \vec \veps(\vec S_{\vec 
u}^{i, 
k+1}+ \vec S_{\vec u}^{i, k})\|^2 + \dfrac{\lambda}{ 4 \beta_{ii}}\|  \nabla \cdot 
(\vec S_{\vec u}^{i, k+1} + \vec S_{\vec u}^{i, k})\|^2\bigg\} \\[1ex]
& \le  \sum_{i =1}^r  \dfrac{L}{ 2 \beta_{ii}} \| S_p^{i, k} \|^2  + \sum_{i =1}^r 
\dfrac{b}{4 \beta_{ii}} \la S_p^{i, k+1} - S_p^{i, k}, 
\nabla \cdot (\vec S_{\vec u}^{i, k+1} - \vec S_{\vec u}^{i, k}) \ra\\[1ex]
& \le  \sum_{i 
=1}^r  \dfrac{L}{ 2 \beta_{ii}} \| S_p^{i, k} \|^2 + \sum_{i 
=1}^r  \dfrac{b^2}{4 \lambda \beta_{ii}} \| S_p^{i, k+1} - S_p^{i, k} \|^2\,.
\end{split}
\end{equation}
From \cite[Lemma 2.3]{BRK15} we conclude that
\begin{align}
\nonumber
\sum_{i =1}^r &  \langle \vec K^{-1} \vec S_{\vec q}^{i, k+1}, \vec E_{\vec q}^{i, 
k+1}\ra \\[1ex]
\nonumber
& = \dfrac{1}{2} \langle \vec K^{-1} \vec e_{\vec q}^{k+1} (t_{n}), \vec 
e_{\vec q}^{k+1} (t_{n}) \rangle - \dfrac{1}{2}\langle \vec K^{-1} \vec e_{\vec q}^{k+1} 
(t_{n-1}), \vec e_{\vec q}^{k+1} (t_{n-1}) \rangle\\[1ex]
& = \dfrac{1}{2} \langle \vec K^{-1} \vec e_{\vec q}^{k+1} (t_{n}), \vec 
e_{\vec q}^{k+1} (t_{n}) \rangle  \geq \dfrac{k}{2} \| \vec e_{\vec q}^{k+1} 
(t_{n}) \|^2\,,
\label{proof_eq_23}
\end{align}
due to $\vec e_{\vec q}^{k+1} (t_{n-1}) = \vec 0$ by means of Eq.\ \eqref{Eq:Ek_eq_0}. In 
Eq.\ \eqref{proof_eq_23} the constant $k$ denotes the lower bound of the uniformly 
positive definite matrix $\vec K^{-1}$.

We are now in a position to perform our final contraction argument.

\medskip
\noindent\textbf{7.\ Step (Contraction argument).} 
Combining Eq.\ \eqref{proof_eq_22} with Eq.\ \eqref{proof_eq_23} shows that 
\begin{equation}\label{proof_eq_25}
\begin{split}
\sum_{i =1}^r \bigg(\dfrac{1}{M \beta_{ii}}  + \dfrac{L}{ 2 \beta_{ii}}\bigg)\| S_p^{i, 
k+1} \|^2 & + \sum_{i =1}^r  \dfrac{L}{ 2 \beta_{ii}} \| S_p^{i, k+1} - S_p^{i, k} \|^2 + 
k\, \dfrac{\tau_n}{2} \| \vec e_{\vec q}^{k+1} (t_{n}) \|^2 \\
& \le   \sum_{i =1}^r  \dfrac{L}{ 2 \beta_{ii}} \| S_p^{i, k} \|^2 + \sum_{i 
=1}^r  \dfrac{b^2}{4 \lambda \beta_{ii}} \| S_p^{i, k+1} - S_p^{i, k} \|^2.
\end{split}
\end{equation}
The inequality \eqref{proof_eq_25} shows the geometric convergence of the 
iterates $S_p^{i, k}$ in $W_h$, for $i=1,\ldots, r$, for any parameter $L \ge 
b^2/(2 \lambda)$. The optimal choice of $L$, still ensuring the geometric 
convergence, is thus given by $L = {b^2}/(2 \lambda)$. The geometric convergence of 
$S_p^{i, k}$ along with Eq.\ \eqref{proof_eq_25} then implies the convergence of $\vec 
e_{\vec q}^{k} (t_{n})$ to $\vec{0}$ for 
$k\rightarrow \infty$. 

By using error equation \eqref{proof_eq_8} for Darcy's law together with the 
convergence of $\vec e_{\vec q}^{k} (t_{n})$ to $\vec{0}$ for $k\rightarrow \infty$ we 
directly get the convergence of $e_p^{k} (t_{n})$ to $0$ for $k \rightarrow \infty$. 
Moreover, the error equation \eqref{proof_eq_9} for the subproblem of mechanics 
deformation along with the previous convergence results then implies the convergence of 
$\vec e_{\vec u}^{k} (t_{n})$ to $\vec {0}$ for $k \rightarrow \infty$. This proves the 
assertion of the theorem.
\end{mproof}

\begin{rem}
The optimal constant  $L = b^2/(2 \lambda)$ identified in the previous proof is the same 
as the one that is obtained in Thm.\ \ref{Thm:ConvCont} for the iteration 
scheme \eqref{Eq:WPF1}, \eqref{Eq:WPF2} and \eqref{Eq:WPM} on the level of the partial 
differential equations, even though different different techniques of proof are applied 
in Thm.\ \ref{Thm:ConvCont} and Thm.\ \ref{Thm:ConvDisc}, respectively. We note that the 
optimal choice of $L$ does not depend on the time stepping scheme, i.e.\ on the 
particular choice of the parameter $r$. Moreover, our result is consistent to the 
analysis given in \cite{MW13} where a convergence proof is given for the continuous case 
of subproblems of partial differential equations with the flow problem being written 
in a non-mixed setting.   
\end{rem}

\begin{cor}
\label{Cor:ContGaussConv}
For $j=0,\ldots, r$, the iterates $\{P^{j,k}_{n,h},\vec Q^{j,k}_{n,h},\vec 
U^{j,k}_{n,h}\}$ converge to $\{P^{j}_{n,h},\vec Q^{j}_{n,h},\vec 
U^{j}_{n,h}\}$ for $k\rightarrow \infty$ in $W_h\times \vec V_h \times \vec H_h$. This 
implies the convergence of $p^k_{\tau,h}$, $\vec q^k_{\tau,h}$ and $\vec u^k_{\tau,h}$ in 
$L^2(I_n;W)$ and 
$L^2(I_n;\vec L^2(\Omega))$, respectively. 
\end{cor}

\begin{mproof}
From \cite[Lemma 2.3]{BRK15} along with the first of the identities \eqref{Eq:Ek_eq_0} we 
conclude that 
\begin{equation}
\label{Eq:Cor1}
\frac{1}{2}\, \|E_p^k(t_n)\|^2 = \int_{t_{n-1}}^{t_n} \langle \partial_t e_p^k, e_p^k 
\rangle \ud t  =  \sum_{i = 1}^r  \sum_{j = 1}^r \alpha_{ij} \langle
E_{p}^{j,k}, E_{p}^{i,k} \rangle.
\end{equation}
Since the matrix $(\alpha_{ij})_{i,j = 1, \dots, r}$ in Eq.\ \eqref{Eq:Cor1} 
is  positive definite (cf.\ \cite[p.\ 1784]{Karakashian1999}) it follows that 
\begin{equation}
\label{Eq:Cor2}
\frac{1}{2}\, \|E_p^k(t_n)\|^2 \geq \alpha_0 \sum_{j=1}^r \|E_p^{i,k}\|^2
\end{equation}
with some constant $\alpha_0>0$. From \eqref{Eq:Cor2} along with Thm.\ \ref{Thm:ConvDisc} 
we conclude the convergence of $P_{n,h}^{j,k}$ in $\vec W_h$. The convergence of $E_{\vec 
q}^{i,k}$ and $E_{\vec u}^{i,k}$ to $\vec 0$ for $k\rightarrow \infty$ is then a direct 
consequence of \eqref{proof_eq_8} and \eqref{proof_eq_9}, respectively. Finally, the 
convergence of $p^k_{\tau,h}$ in $L^2(I_n;W)$ and of $\vec q^k_{\tau,h}$, $\vec 
u^k_{\tau,h}$ in $L^2(I_n;\vec L^2(\Omega))$ follows from the second result in 
\cite[Lemma 2.3]{BRK15}.
\end{mproof}

\subsection{The dG($\vec r$)--MFEM($\vec s$)cG($\vec s$+1) approach.}

In this subsection we prove the (linear) convergence of the splitting schemes 
\eqref{Eq:D_AlgebProb1}--\eqref{Eq:D_AlgebProb3} based on a discontinuous Galerkin 
discretization of the time variable. Again, we show that the scheme is 
subject to a contraction principle such that a unique a fixed point is obtained. 

In the sequel, we denote by $p_{\tau,h} \in \widetilde{\mathcal{W}}_{\tau,h}^{r,s}$, 
$\vec q_{\tau,h}\in \vec{\widetilde{\mathcal{V}}}_{\tau,h}^{r,s}$ and $\vec u_{\tau,h} \in 
\vec{\widetilde{\mathcal{Z}}}_{\tau,h}^{r,s+1}$, with 
\begin{equation}
\label{EQ:dGRepFDS}
p_{\tau,h}{}_{|I_n} (t) = \sum_{j=0}^r P^{j}_{n,h} \varphi_{n,j}(t)\,, \; \vec 
q_{\tau,h}{}_{|I_n} (t) = \sum_{j=0}^r \vec Q^{j}_{n,h} \varphi_{n,j}(t)\,, \;\vec 
u_{\tau,h}{}_{|I_n} (t) = \sum_{j=0}^r 
\vec U^{j}_{n,h} 
\varphi_{n,j}(t)
\end{equation}
for $t\in I_n$, the space-time finite element approximation of the Biot system 
\eqref{Eq:B_4}--\eqref{Eq:B_7} that is defined by skipping the upper indices in the 
problems \eqref{Eq:D_FDM}, \eqref{Eq:D_FDD} and \eqref{Eq:D_FDMD}, respectively. Thus we 
tacitly suppose that the coupled system that is obtained by discretizing the Biot 
model \eqref{Eq:B_4}--\eqref{Eq:B_7} in the space-time finite element spaces 
\eqref{Def:FE4}--\eqref{Def:FE6} admits a unique solution.

\begin{theorem} 
\label{Thm:DGConvDisc}
Let $p_{\tau,h} \in \widetilde{\mathcal{W}}_{\tau,h}^{r,s}$, $\vec q_{\tau,h}\in 
\vec{\widetilde{\mathcal{V}}}_{\tau,h}^{r,s}$ and $\vec u_{\tau,h} \in 
\vec{\widetilde{\mathcal{Z}}}_{\tau,h}^{r,s+1}$  denote the 
fully  discrete space-time finite element approximation of the Biot system 
\eqref{Eq:B_4}--\eqref{Eq:B_7}.  On $I_n$ let $\{p_{\tau,h},\vec q_{\tau,h},\vec u_{\tau,h}\}$ be represented by \eqref{EQ:dGRepFDS} and let $\{p^k_{\tau,h} , \vec 
q^k_{\tau,h},\vec u^k_{\tau,h} \}$ be defined by \eqref{Eq:dGRepSolBasis}  with coefficient functions being given by the scheme \eqref{Eq:D_AlgebProb1}--\eqref{Eq:D_AlgebProb3}. Then, for any $L \geq 
b^2/(2\lambda)$ the sequence $\{S_{p}^{i,k}\}_k$, for $i=1,\ldots,r$, converges 
geometrically in $W_h$. For 
$n=1,\ldots, N$ this implies the convergence of 
$\{p^k_{\tau,h}(t_n^\pm),\vec q^k_{\tau,h}(t_n^\pm),\vec u^k_{\tau,h}(t_n^\pm)\}$ to 
$ \{p_{\tau,h}(t_n^\pm),\vec q_{\tau,h}(t_n^\pm),\vec u_{\tau,h}(t_n^\pm)\}$ in 
$W_h\times 
\vec V_h \times \vec H_h$ for $k\rightarrow \infty$. 
\end{theorem}

\begin{mproof}
The proof follows the lines of the proof of Thm.\ \ref{Thm:ConvDisc}. 
Therefore, we restrict ourselves to presenting the differences only. We use the notation and abbreviations of Subsec.\ \ref{Subsec::cGConv}.

We put
\begin{equation*}
S_p^{i, k+1} = \sum_{j =0}^r \widetilde\alpha_{ij} E_p^{j, k+1}\,, 
\quad
\vec S_{\vec q}^{i, k+1} =  \sum_{j =0}^r \widetilde\alpha_{ij} \vec 
E_{\vec q}^{j, k+1}\,,\quad 
\vec S_{\vec u}^{i, k+1} =  \sum_{j =0}^r \widetilde\alpha_{ij} \vec 
E_{\vec u}^{j, k+1}
\end{equation*}
with $S_p^{i, k+1} \in W_h$, $\vec S_{\vec q}^{i, k+1} \in \vec V_h$ and $\vec 
S_{\vec u}^{i, k+1} \in \vec  H_h$ for $i=0.,\ldots,r$. 
By the same arguments as in the proof of Thm.\ \ref{Thm:ConvDisc} we find that 
\begin{equation}
\label{Eq:DG0}
\begin{split}
\sum_{i =0}^r & \bigg\{ \bigg(\dfrac{1}{M \beta_{ii}} + \dfrac{L}{ 
2 \beta_{ii}}\bigg) \| S_p^{i, k+1} \|^2   + \dfrac{L}{ 2 \beta_{ii}} \| S_p^{i, k+1} - 
S_p^{i, k} \|^2 + \tau_n 
\la \vec K^{-1} \vec S_{\vec q}^{i, k+1}, \vec E_{\vec q}^{i, k+1}\ra\bigg\} \\[1ex]
& \quad + \sum_{i =0}^r \bigg \{ \dfrac{\mu}{ 2 \beta_{ii}} \|  \vec \veps(\vec S_{\vec 
u}^{i, 
k+1}+ \vec S_{\vec u}^{i, k})\|^2 + \dfrac{\lambda}{ 4 \beta_{ii}}\|  \nabla \cdot 
(\vec S_{\vec u}^{i, k+1} + \vec S_{\vec u}^{i, k})\|^2\} \\[1ex]
& \le  \sum_{i =0}^r  \dfrac{L}{ 2 \beta_{ii}} \| S_p^{i, k} \|^2 + \sum_{i 
=0}^r  \dfrac{b^2}{4 \lambda \beta_{ii}} \| S_p^{i, k+1} - S_p^{i, k} \|^2\,.
\end{split}
\end{equation}

From \cite[Lemma 2.3]{BRK15} we get that
\begin{align}
\nonumber
\sum_{i =0}^r &  \langle \vec K^{-1} \vec S_{\vec q}^{i, k+1}, \vec E_{\vec q}^{i, 
k+1}\ra \\[0.5ex]
\nonumber
& = \sum_{i =0}^r \bigg\langle \vec K^{-1} \sum_{j=0}^r \alpha_{ij} \vec E_{\vec q}^{j, 
k+1}, \vec E_{\vec q}^{i,k+1}\bigg\rangle + \sum_{i =0}^r \bigg\langle \vec K^{-1} 
\sum_{j=0}^r \gamma_i\gamma_j \vec E_{\vec q}^{j, 
k+1}, \vec E_{\vec q}^{i,k+1}\bigg\rangle \\[1ex] 
\nonumber
& = \dfrac{1}{2} \langle \vec K^{-1} \vec E_{\vec q}^{k+1} (t_{n}^-), 
\vec 
e_{\vec q}^{k+1} (t_{n}^-) \rangle - \dfrac{1}{2}\langle \vec K^{-1} 
\vec e_{\vec q}^{k+1} 
(t_{n-1}^+), \vec e_{\vec q}^{k+1} (t_{n-1}^+) \rangle\\[0.5ex]
\nonumber
& \quad + \bigg\langle \vec K^{-1} 
\sum_{j=0}^r \gamma_j \vec E_{\vec q}^{j, 
k+1}, \sum_{i =0}^r  \gamma_i \vec E_{\vec q}^{i,k+1}\bigg\rangle\\[1ex]
& = \dfrac{1}{2} \langle \vec K^{-1} \vec e_{\vec q}^{k+1} (t_{n}^-), 
\vec e_{\vec q}^{k+1} (t_{n}^-) \rangle + \dfrac{1}{2}\langle \vec K^{-1} 
\vec e_{\vec q}^{k+1}(t_{n-1}^+), \vec e_{\vec q}^{k+1} (t_{n-1}^+) \rangle\,,
\label{Eq:DG1}
\end{align}
where we used that 
\begin{equation*}
 \bigg\langle \vec K^{-1} 
\sum_{j=0}^r \gamma_j \vec E_{\vec q}^{j, 
k+1}, \sum_{i =0}^r  \gamma_i \vec E_{\vec q}^{i,k+1}\bigg\rangle= \langle\vec K^{-1} 
\vec e_{\vec q}^{k+1} (t_{n-1}^+), \vec e_{\vec q}^{k+1} (t_{n-1}^+)  \rangle\,.
\end{equation*}

We are now in a position to perform our final contraction argument. Combining Eq.\ 
\eqref{Eq:DG0} with Eq.\ \eqref{Eq:DG1} shows that 
\begin{equation}
\label{Eq:DG2}
\begin{split}
& \sum_{i =0}^r \bigg(\dfrac{1}{M \beta_{ii}}  + \dfrac{L}{ 2 
\beta_{ii}}\bigg)\| S_p^{i, 
k+1} \|^2 + \sum_{i =0}^r  \dfrac{L}{ 2 \beta_{ii}} \| S_p^{i, k+1} 
- S_p^{i, k} \|^2\\[1ex]
& \qquad \qquad + k \dfrac{\tau_n}{2} \| \vec e_{\vec q}^{k+1} 
(t_{n}^-)\|^2 + k \dfrac{\tau_n}{2} \| \vec e_{\vec q}^{k+1} 
(t_{n-1}^+)\|^2\\[1ex]
& \le   \sum_{i =0}^r  \dfrac{L}{ 2 \beta_{ii}} \| S_p^{i, k} \|^2 + 
\sum_{i 
=0}^r  \dfrac{b^2}{4 \lambda \beta_{ii}} \| S_p^{i, k+1} - S_p^{i, k} 
\|^2\,.
\end{split}
\end{equation}
The inequality \eqref{Eq:DG2} shows the geometric convergence of the iterates $S_p^{i, 
k}$ in $W_h$, for $i=0,\ldots, r$, for any parameter $L \ge b^2/(2 \lambda)$. Again, the 
optimal choice of $L$ is given by $L=b^2/(2 \lambda)$. The geometric convergence of 
$S_p^{i, k}$ along with Eq.\ \eqref{Eq:DG2} then implies the convergence of $\vec e_{\vec 
q}^{k} (t_{n-1}^+)$ and $\vec e_{\vec q}^{k} (t_{n}^-)$ to $\vec{0}$ for $k\rightarrow 
\infty$. 

\bigskip
By using error equation \eqref{proof_eq_8} for Darcy's law together with the 
convergence of $\vec e_{\vec q}^{k} (t_{n}^-)$ to $\vec{0}$ for $k\rightarrow \infty$ we 
directly get the convergence of $e_p^{k} (t_{n}^-)$ to $0$ for $k \rightarrow \infty$. 
Moreover, the error equation \eqref{proof_eq_9} for the subproblem of mechanics 
deformation along with the previous convergence results then implies the convergence of 
$\vec e_{\vec u}^{k} (t_{n}^-)$ to $\vec 
{0}$ for $k \rightarrow \infty$. The convergence of $e_p^{k} (t_{n-1}^+)$ to ${0}$ 
and $\vec e_{\vec u}^{k} (t_{n-1}^+)$ to $\vec {0}$ for $k \rightarrow \infty$ follow 
similarly.

\end{mproof}

\begin{cor}
\label{Cor:DGDisGaussConv}
For $j=0,\ldots, r$, the iterates $\{P^{j,k}_{n,h},\vec Q^{j,k}_{n,h},\vec 
U^{j,k}_{n,h}\}$ converge to $\{P^{j}_{n,h},\vec Q^{j}_{n,h},\vec 
U^{j}_{n,h}\}$ for $k\rightarrow \infty$ in $W_h\times \vec V_h \times \vec H_h$. This 
implies the convergence of $p^k_{\tau,h}$, $\vec q^k_{\tau,h}$ and $\vec u^k_{\tau,h}$ in 
$L^2(I_n;W)$ and $L^2(I_n;\vec L^2(\Omega))$, respectively. 
\end{cor}

\begin{mproof}
From \cite[Lemma 2.3]{BRK15} we conclude that 
\begin{equation}
\label{Eq:Cor1_2}
\frac{1}{2}\, \|e_p^k(t_n^-)\|^2 - \frac{1}{2}\, \|e_p^k(t_{n-1}^+)\|^2 = 
\int_{t_{n-1}}^{t_n} \langle \partial_t e_p^k, e_p^k 
\rangle \ud t  =  \sum_{i = 0}^r  \sum_{j = 0}^r \alpha_{ij} \langle
E_{p}^{j,k}, E_{p}^{i,k} \rangle.
\end{equation}
Since the matrix $(\alpha_{ij})_{i,j = 0, \dots, r}$ in Eq.\ \eqref{Eq:Cor1_2} 
is  positive definite (cf.\ \cite[p.\ 1784]{Karakashian1999}) it follows that 
\[
\frac{1}{2}\, \|e_p^k(t_n^-)\|^2 - \frac{1}{2}\, \|e_p^k(t_{n-1}^+)\|^2 \geq \alpha_0 
\sum_{j=0}^r \|E_p^{i,k}\|^2
\]
with some constant $\alpha_0>0$. Thm.\ \ref{Thm:DGConvDisc} then implies the convergence 
of $P^{j,k}$ to $P^{j}$ for $k\rightarrow \infty$ and $j=0,\ldots 
,r$. The convergence of $\{\vec Q^{j,k}_{n,h},\vec U^{j,k}_{n,h}\}$ to 
$\{\vec Q^{j}_{n,h},\vec U^{j}_{n,h}\}$ is now a direct consequence of 
\eqref{Eq:D_AlgebProb2} and \eqref{Eq:D_AlgebProb3}, respectively. Finally, the 
convergence of $e_p^k$, $\vec e_{\vec q}^k$ and $\vec e_{\vec u}^k$ to zero in 
$L^2(I_n;W)$ and $L^2(I_n;\vec L^2(\Omega))$, respectively, follows from the exactness of 
the $r$+1-point Gauss quadrature formula on $I_n$ all for polynomials of maximum degree 
$2r+1$. 
\end{mproof}

\section{Numerical experiments} 
\label{Sec:NumRes}

In this section we study the numerical performance properties of the fixed-stress  
iteration schemes \eqref{Eq:FDM}--\eqref{Eq:FDMD} 
and \eqref{Eq:D_FDM}--\eqref{Eq:D_FDD}, respectively, along with with the proposed 
choice of our analyses $L = b/(2\lambda)$ for the numerical parameter $L$. For the time 
discretization we consider a continuous approximation with piecewise linear and quadratic 
polynomials, i.e.\ a cGP(1) and cGP(2) approach (cf.\ Sec.\ \ref{Sec:cGPMFEM}), as well 
as a discontinuous approximation with piecewise constant and linear polynomials, i.e.\ a 
dG(0) and dG(1) approach (cf.\ Sec.\ \ref{Sec:dGMFEM}). In our computations we shall 
study numerically the sharpness of our theoretical result of Sec.\ \ref{Sec:Convergence} 
that $L=b^2/(2\lambda)$ provides an optimal choice of $L$ with respect to an 
acceleration of the convergence behaviour of the fixed point iterations. The 
implementation of the schemes is done in our front-end simulation tool for the latest 
\texttt{deal.II} version 8 library and allows distributed-parallel numerical simulations; 
cf.\ \cite{BK16,DEAL,K15,BK15} for details.

\begin{figure}
\centering
\subfloat[Test setting]{
\includegraphics[width=.30\linewidth]{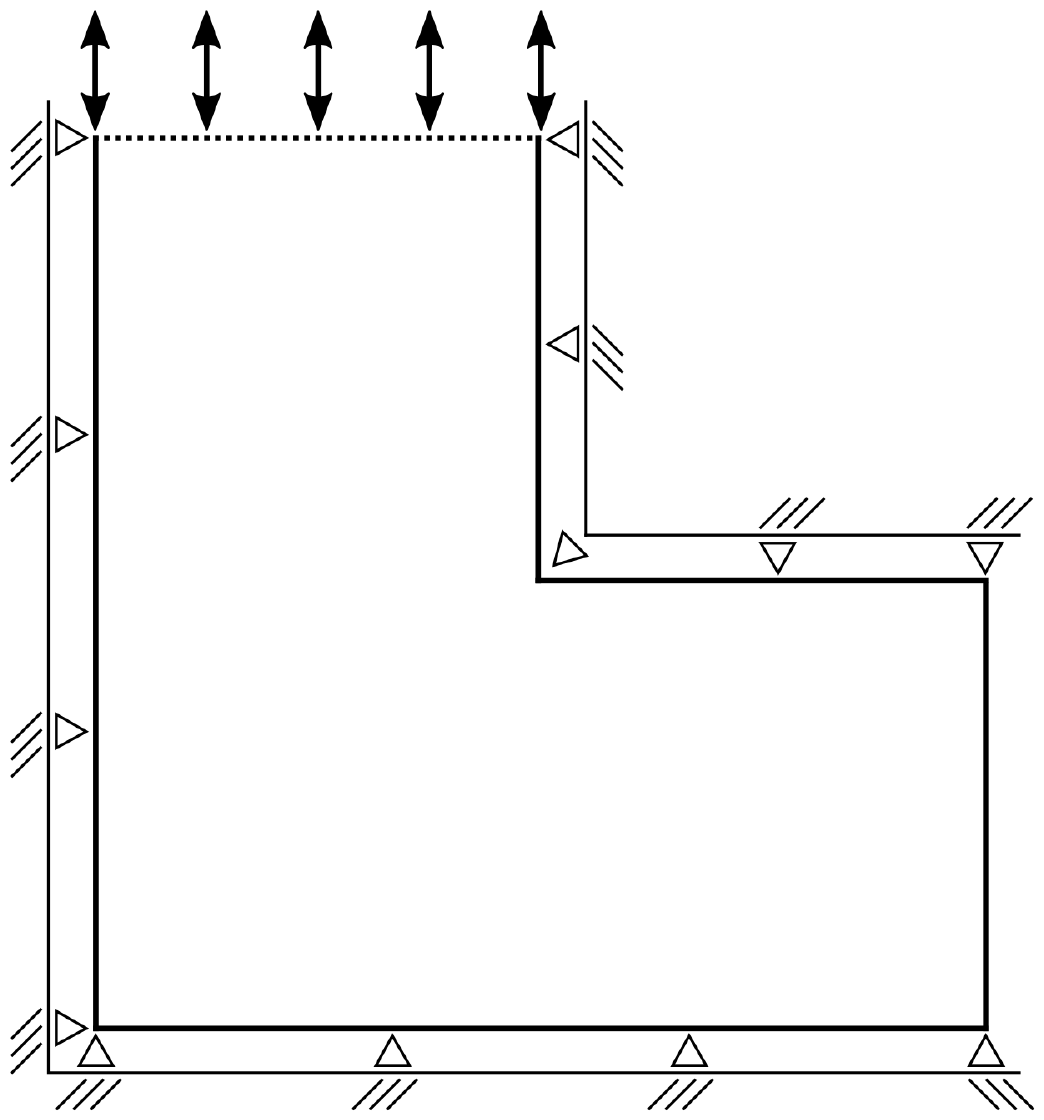}
}%
\subfloat[Pressure $p$ for $t=0.26$]{
\includegraphics[width=.30\linewidth]{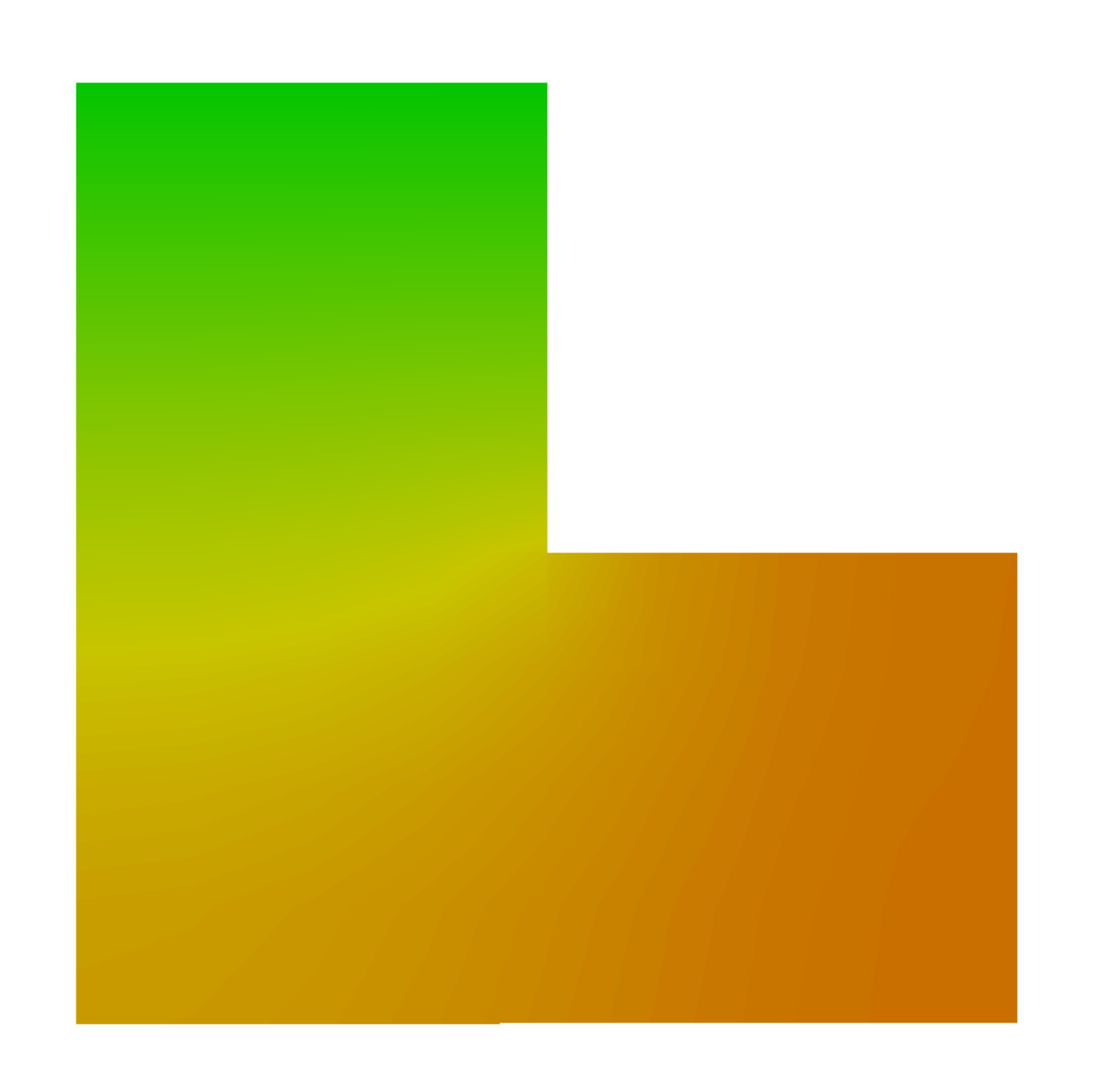}
}%
\subfloat[Displacement $\boldsymbol u$ for $t=0.26$]{
\includegraphics[width=.30\linewidth]{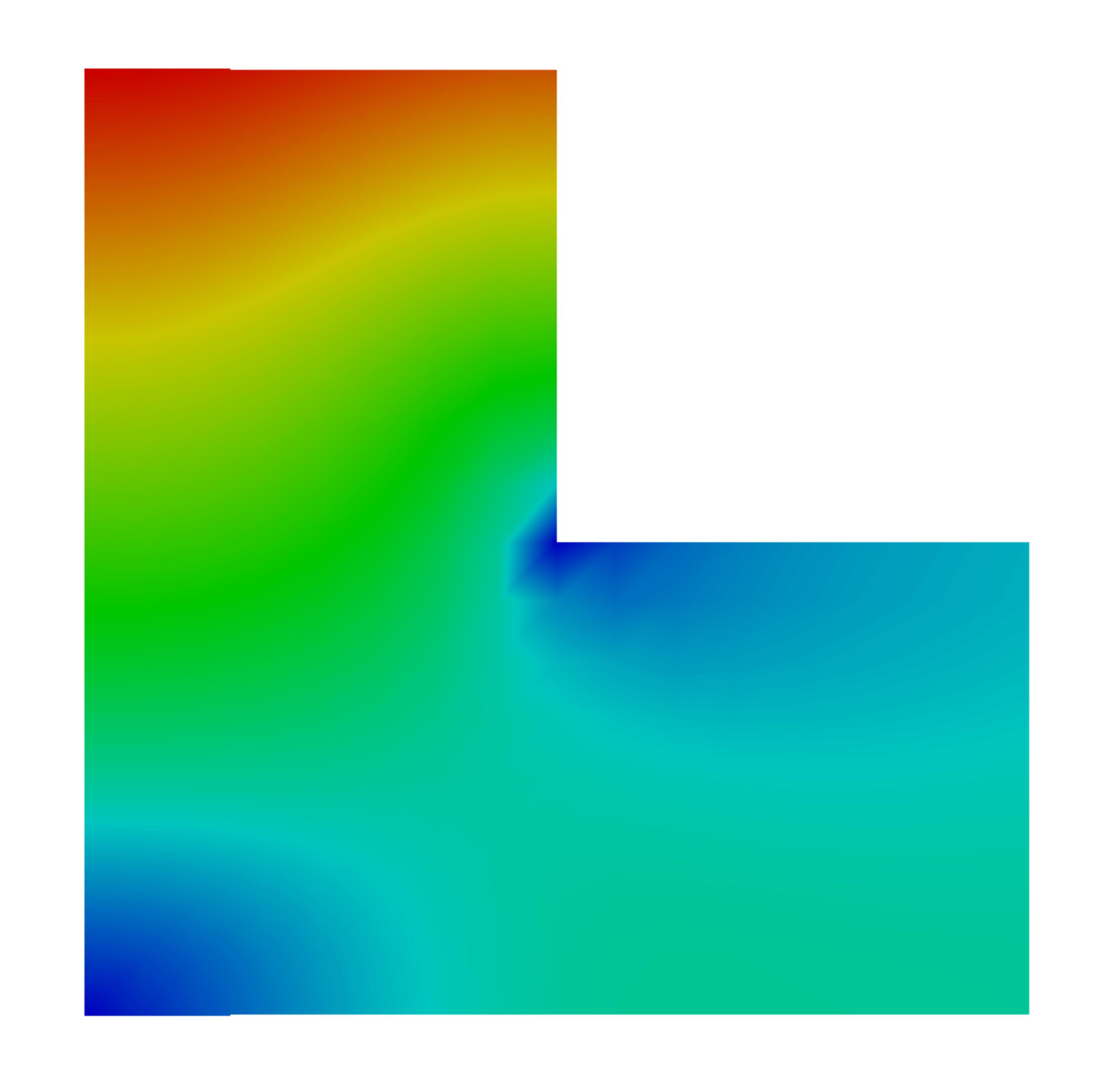}
}%
\caption{Test setting, pressure and magnitude of displacement field at $t=0.26$ for time 
step size $\tau_n=0.02$ and solver cGP(2)--MFEM(1)cG(2); cf.\ Sec.\ 
\ref{Sec:cGPMFEM}.}
\label{fig:1:SettingSolution}
\end{figure}

The problem setting of our test configuration with an L-shaped domain is sketched in 
Fig.~\ref{fig:1:SettingSolution}. We consider solving the Biot 
problem in the time interval $I=(0,0.5)$. We prescribed homogeneous initial conditions. 
The solid lines describe an undrained flow boundary (i.e.\ $\vec q \cdot \vec n = 0$ with 
outer unit normal vector $\vec n$) and the dashed line on the top describes a open flow 
boundary with a prescribed pressure value $p = 0$. At the open flow boundary at the top 
we prescribe a time-dependent traction boundary condition for mechanical deformation 
given by $\vec \sigma \vec n  = (0, h(t))^\top$, with $h(t) = -2560\, t^2\, (t-0.5)^2$. 
At the lower right boundary a homogeneous traction boundary condition is imposed. At all 
remaining boundaries we prescribe one displacement component to fulfil a homogeneous 
Dirichlet condition and the remaining component to fulfil a homogeneous traction boundary 
condition. The physical parameters are chosen as $M=100$, $b=100$, $\mu = E/(2\cdot 
(1+\nu))$ and $\lambda = E\nu /((1-2\nu)\cdot (1+\nu))$ with $E=100$ and $\nu =0.35$ 
such that $\mu = 37.037$ and $\lambda =86.42$. Further we put $\vec K = 0.1 \cdot \vec I$ 
with the identity matrix $\vec I$. Gravity is not considered, i.e.\ $\vec g\equiv \vec 
0$. The calculated profiles for fluid pressure and magnitude of displacement are 
illustrated exemplarily for $t=0.26$ in 
Fig.\ \ref{fig:1:SettingSolution}. For the pressure distribution the 
green coloured region corresponds to $p(\cdot,0.26)=0$
and rises up to $p(\cdot,0.26) \approx 0.5$ in the orange coloured 
regions. For the displacement field magnitude distribution the blue 
coloured region corresponds to $\|\boldsymbol u(\cdot,0.26)\|_2 = 0$ 
and rises up to $\|\boldsymbol u(\cdot,0.26)\|_2 \approx 0.06$ in the 
red coloured region at the top of the domain.

As a stopping criterion for the fixed-stress iteration we prescribed a tolerance of  
$\text{tol}_\text{fixed} = 1\mathrm{e}-8$, measured in the $l^2$ norm, between two 
successive solution vectors for each of the unknown variables, i.e.\ pressure, flux and 
displacement field.
For the lower order time discretizations dG(0) and cG(1) we chose
$\text{tol}_\text{flow} = 1\mathrm{e}-14$ and
$\text{tol}_\text{mechanics} = 1\mathrm{e}-12$ for the iterative solvers of the 
subproblems.  
For the higher order time discretisations dG(1) and cG(2) we put  
$\text{tol}_\text{flow} = 1\mathrm{e}-12$, and
$\text{tol}_\text{mechanics} = 1\mathrm{e}-12$ for the iterative solver tolerances. 
\input{./h-independent}

In our first numerical study the sensitivity (cf.\ Fig.\ \ref{fig:2:dG0:h}) of the 
iteration process with respect to choice of the spatial discretization step size $h$ is 
analyzed. This is done for a lowest order in time discontinuous Galerkin discretization 
dG(0) and a MFEM(0)cG(1) approximation in space; cf.\ Sec.\ \ref{Sec:dGMFEM}. In Fig.\ 
\ref{fig:2:dG0:h} the total number of iterations for all time steps in the interval $I$ 
and step size $\tau_n = 0.01$ is illustrated versus a perturbation $\omega$ 
of our optimal choice of the numerical tuning parameter. Precisely, we performed our 
iterations with $L = \omega \widehat L$ where $\widehat L = b^2/(2\lambda)$ is the choice 
that is proposed by our analysis such that $\omega=1$ represents the theoretically 
expected result for the best performance of the iteration scheme with a minimum number of 
iterations. In Fig.\ \ref{fig:2:dG0:h} we observe a convergence behaviour that is almost 
independent of the refinement level $m$ with $h=2^{-(m+1)}$. For all refinement levels 
the computations show the optimal convergence behaviour for values slightly greater 
than one for the perturbation parameter, $\omega \approx 1.05$, such that our proposed 
choice of Sec.\ \ref{Sec:Convergence} corresponding to $\omega=1$ fits quite well. We note 
that for stronger perturbations of $\omega=1$ the number of required iterations increases 
strongly which leads to additional numerical costs.   

\input{./p-independent}

In our second numerical study the sensitivity (cf.\ Fig.\ \ref{fig:3:cG1:p}) of the 
iteration process with respect to a variation of the polynomial degree $s$ of the spatial 
discretization is analyzed; cf.\ Sec.\ \ref{Sec:cGPMFEM}. We vary the parameter 
$s$ from  $s=1$ to $s=4$. For the time discretization the lowest order continuous 
Galerkin apporach cGP(1) is applied; cf.\ Sec.\ \ref{Sec:cGPMFEM}. In Fig.\ 
\ref{fig:3:cG1:p} we illustrate the total number of iterations for all time steps in 
the time interval $I$ versus a perturbation of our proposed choice of the tuning 
parameter. As before, $\omega=1$ corresponds to the proposed value of our analysis in 
Sec.\ \ref{Sec:Convergence}. Again, in our computations the iterations show strong 
robustness with respect to the choice of $s$. Therefore, the result of our 
analysis, corresponding to $\omega =1$, is close to the optimal point of a 
minimum number of iterations. For the higher order variants a value 
of $\omega$ slighlty greater than $\omega = 1$ seems to be 
advantageous. Nevertheless, the great impact of our analysis for 
the choice of the optimal numerical parameter $L$ is obvious.        

\input{./tau-independent}

Next, in our third numerical study the sensitivity (cf.\ Fig.\ \ref{fig:4:dG1:tau}) of 
the iteration process with respect to the choice of the time step size is analyzed. This 
is done for the dG(1) time discretization scheme; cf.\ Sec.\ \ref{Sec:dGMFEM}. Halfening 
the time step size and thereby doubling the number of time steps doubles the total number 
of iterations for the fixed-stress splitting solution in the interval 
$I$. Again, the results of our analyses in Sec.\ \ref{Sec:IterScheme} 
for the continuous case and in Sec.\ \ref{Sec:Convergence} for the 
discrete case are confirmed by the illustrated dependence of the 
number of iterations on the perturbation $\omega$. No significant 
difference is observed in the convergence behavior whether a 
continuous cGP(1) or discontinuous dG(1) time discretization is are 
applied. 

Finally in Fig.\ \ref{fig:5:cG2:tau} the same study is presented for 
the higher order cGP(2) approach with a continuous approximation in 
time with piecewise quadratic polynomials. For comparison the total 
number of iterations depending on the perturbation $\omega$ are 
illustrated for the cGP(1) and cGP(2) approach. No significant 
deviations are observed. 

Summarizing, we can state that the numerical results 
nicely confirm our analyses and conjectures given in Sec.\ 
\ref{Sec:IterScheme} and in Sec.\ \ref{Sec:Convergence}, respectively. An almost 
optimal choice of the numerical tuning parameter $L$ in the 
fixed-stress iteration schemes \eqref{Eq:WPF1}--\eqref{Eq:WPM} as 
well as \eqref{Eq:FDM}--\eqref{Eq:FDD} 
and \eqref{Eq:D_FDM}--\eqref{Eq:D_FDMD} is given by 
$L=b^2/(2\lambda)$. This choice only depends on modelling and not on 
discretization parameters.

\section{Summary}
\label{Sec:Summary}

In this work we presented and analyzed an iterative splitting scheme for the numerical 
approximation of the quasi-static Biot system of poroelasticity. For the discretization 
of the separated subproblems of fluid flow and mechanical deformation space-time finite 
element methods of arbitrary polynomial order are used. For the approximation of the 
time variable continuous and discontinuous Galerkin approaches are considered. The 
convergence of the iterative coupling scheme is shown for the continuous model of partial 
differential equations and the fully discrete set of algebraic equations. For both cases 
our analyses propose the same optimal choice of an inherent stabilization or tuning 
parameter of the iterative approach. In particular, the parameter is independent of the 
numerical discretization parameters. Our presented numerical results nicely confirm the 
theoretical results and the expected convergence behaviour. Moreover, they underline the 
efficiency and stability of the proposed approaches for simulating flow in deformable 
porous media modelled by the Biot system. Next, we plan to apply the optimized 
fixed-stress iterative coupling strategy to more complex physical models of flow in 
deformable porous media. In particular, variably saturated and multiphase flow 
\cite{L16,P04,R15} as well as non-linear poroelasticity are in the scope of our interest.
    
\section*{Acknowledgements}

This work was supported by the German Academic Exchange Service (DAAD) under
the grant ID 57238185, by the Research Council of Norway under the grant ID DAADppp255715 
and the Toppforsk projekt under the grant ID 250223.

\end{document}

%% file: h-independent.tex
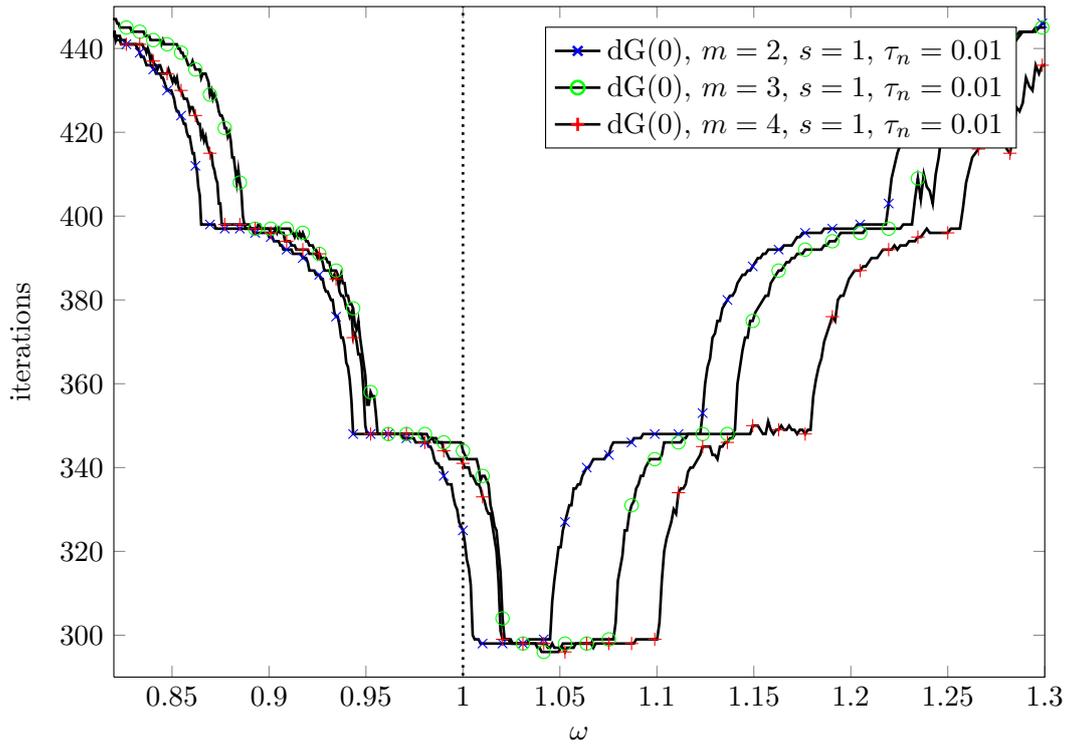
\begin{figure}[p]
\centering
%
\begin{tikzpicture}

\begin{axis}[%
width=4.82222222222222in,
height=3.5in,
scale only axis,
/pgf/number format/.cd, 1000 sep={},
xlabel={$\omega$},
ylabel={iterations},
xmin=0.82,
xmax=1.3,
ymin=290,
ymax=450,
yminorticks=true,
legend style={legend pos=north east},
legend style={draw=black,fill=white,legend cell align=left},
legend entries = {
  {dG(0), $m=2$, $s=1$, $\tau_n=0.01$},
  {dG(0), $m=3$, $s=1$, $\tau_n=0.01$},
  {dG(0), $m=4$, $s=1$, $\tau_n=0.01$}
}
]

\addplot [
color=black,
solid,
line width=1.0pt,
mark=x,
mark size = 2.5,
mark options={solid,blue}
]
table[row sep=crcr]{
0.0 0\\
0.1 0\\
};

\addplot [
color=black,
solid,
line width=1.0pt,
mark=o,
mark size = 2.5,
mark options={solid,green}
]
table[row sep=crcr]{
0.0 0\\
0.1 0\\
};

\addplot [
color=black,
solid,
line width=1.0pt,
mark=+,
mark size = 2.5,
mark options={solid,red}
]
table[row sep=crcr]{
0.0 0\\
0.1 0\\
};

\addplot [
color=black,
dotted,
line width=1.0pt
]
table[row sep=crcr]{
1.0 0\\
1.0 2000\\
};

\addplot [
color=black,
solid,
line width=1.0pt
]
table[row sep=crcr]{
0.769231 472\\
0.769823 471\\
0.770416 471\\
0.77101 471\\
0.771605 470\\
0.772201 468\\
0.772798 468\\
0.773395 468\\
0.773994 467\\
0.774593 465\\
0.775194 464\\
0.775795 464\\
0.776398 464\\
0.777001 461\\
0.777605 460\\
0.77821 460\\
0.778816 457\\
0.779423 456\\
0.780031 452\\
0.78064 447\\
0.78125 447\\
0.781861 447\\
0.782473 447\\
0.783085 447\\
0.783699 447\\
0.784314 447\\
0.784929 447\\
0.785546 447\\
0.786164 447\\
0.786782 447\\
0.787402 447\\
0.788022 447\\
0.788644 447\\
0.789266 447\\
0.789889 447\\
0.790514 447\\
0.791139 447\\
0.791766 447\\
0.792393 447\\
0.793021 447\\
0.793651 447\\
0.794281 447\\
0.794913 447\\
0.795545 447\\
0.796178 447\\
0.796813 447\\
0.797448 447\\
0.798085 447\\
0.798722 447\\
0.799361 447\\
0.8 447\\
0.800641 447\\
0.801282 447\\
0.801925 447\\
0.802568 446\\
0.803213 446\\
0.803859 446\\
0.804505 446\\
0.805153 446\\
0.805802 446\\
0.806452 445\\
0.807103 445\\
0.807754 445\\
0.808407 445\\
0.809061 445\\
0.809717 445\\
0.810373 445\\
0.81103 445\\
0.811688 445\\
0.812348 444\\
0.813008 444\\
0.81367 444\\
0.814332 444\\
0.814996 444\\
0.815661 444\\
0.816327 444\\
0.816993 444\\
0.817661 444\\
0.818331 443\\
0.819001 443\\
0.819672 443\\
0.820345 443\\
0.821018 442\\
0.821693 442\\
0.822368 442\\
0.823045 442\\
0.823723 441\\
0.824402 441\\
0.825083 441\\
0.825764 441\\
0.826446 441\\
0.82713 441\\
0.827815 441\\
0.8285 441\\
0.829187 441\\
0.829876 441\\
0.830565 441\\
0.831255 441\\
0.831947 440\\
0.832639 440\\
0.833333 439\\
0.834028 438\\
0.834725 438\\
0.835422 438\\
0.83612 437\\
0.83682 436\\
0.837521 436\\
0.838223 436\\
0.838926 436\\
0.839631 436\\
0.840336 435\\
0.841043 435\\
0.841751 434\\
0.84246 434\\
0.84317 434\\
0.843882 434\\
0.844595 433\\
0.845309 433\\
0.846024 432\\
0.84674 431\\
0.847458 430\\
0.848176 430\\
0.848896 430\\
0.849618 430\\
0.85034 429\\
0.851064 428\\
0.851789 426\\
0.852515 425\\
0.853242 425\\
0.853971 424\\
0.854701 424\\
0.855432 423\\
0.856164 422\\
0.856898 422\\
0.857633 420\\
0.858369 419\\
0.859107 419\\
0.859845 417\\
0.860585 415\\
0.861326 415\\
0.862069 412\\
0.862813 411\\
0.863558 407\\
0.864304 405\\
0.865052 398\\
0.865801 398\\
0.866551 398\\
0.867303 398\\
0.868056 398\\
0.86881 398\\
0.869565 398\\
0.870322 398\\
0.87108 398\\
0.87184 398\\
0.8726 397\\
0.873362 397\\
0.874126 397\\
0.874891 397\\
0.875657 397\\
0.876424 397\\
0.877193 397\\
0.877963 397\\
0.878735 397\\
0.879507 397\\
0.880282 397\\
0.881057 397\\
0.881834 397\\
0.882613 397\\
0.883392 397\\
0.884173 397\\
0.884956 397\\
0.88574 397\\
0.886525 397\\
0.887311 397\\
0.888099 397\\
0.888889 397\\
0.88968 397\\
0.890472 397\\
0.891266 397\\
0.892061 396\\
0.892857 396\\
0.893655 396\\
0.894454 396\\
0.895255 396\\
0.896057 396\\
0.896861 396\\
0.897666 396\\
0.898473 396\\
0.899281 396\\
0.90009 395\\
0.900901 395\\
0.901713 395\\
0.902527 394\\
0.903342 394\\
0.904159 394\\
0.904977 394\\
0.905797 393\\
0.906618 393\\
0.907441 392\\
0.908265 392\\
0.909091 392\\
0.909918 392\\
0.910747 392\\
0.911577 392\\
0.912409 391\\
0.913242 391\\
0.914077 391\\
0.914913 391\\
0.915751 391\\
0.91659 391\\
0.917431 390\\
0.918274 390\\
0.919118 390\\
0.919963 389\\
0.92081 387\\
0.921659 387\\
0.922509 387\\
0.923361 387\\
0.924214 386\\
0.925069 386\\
0.925926 386\\
0.926784 385\\
0.927644 385\\
0.928505 383\\
0.929368 382\\
0.930233 382\\
0.931099 381\\
0.931966 380\\
0.932836 380\\
0.933707 377\\
0.934579 376\\
0.935454 376\\
0.93633 374\\
0.937207 371\\
0.938086 371\\
0.938967 369\\
0.93985 365\\
0.940734 364\\
0.94162 361\\
0.942507 356\\
0.943396 348\\
0.944287 348\\
0.94518 348\\
0.946074 348\\
0.94697 348\\
0.947867 348\\
0.948767 348\\
0.949668 348\\
0.95057 348\\
0.951475 348\\
0.952381 348\\
0.953289 348\\
0.954198 348\\
0.95511 348\\
0.956023 348\\
0.956938 348\\
0.957854 348\\
0.958773 348\\
0.959693 348\\
0.960615 348\\
0.961538 348\\
0.962464 348\\
0.963391 348\\
0.96432 348\\
0.965251 348\\
0.966184 348\\
0.967118 348\\
0.968054 348\\
0.968992 347\\
0.969932 347\\
0.970874 347\\
0.971817 347\\
0.972763 347\\
0.97371 347\\
0.974659 347\\
0.97561 346\\
0.976562 346\\
0.977517 346\\
0.978474 346\\
0.979432 346\\
0.980392 346\\
0.981354 346\\
0.982318 345\\
0.983284 345\\
0.984252 344\\
0.985222 342\\
0.986193 341\\
0.987167 341\\
0.988142 340\\
0.98912 340\\
0.990099 338\\
0.99108 337\\
0.992063 336\\
0.993049 336\\
0.994036 336\\
0.995025 333\\
0.996016 331\\
0.997009 330\\
0.998004 327\\
0.999001 326\\
1 325\\
1.001 322\\
1.002 318\\
1.00301 316\\
1.00402 312\\
1.00503 300\\
1.00604 299\\
1.00705 299\\
1.00806 299\\
1.00908 298\\
1.0101 298\\
1.01112 298\\
1.01215 298\\
1.01317 298\\
1.0142 298\\
1.01523 298\\
1.01626 298\\
1.01729 298\\
1.01833 298\\
1.01937 298\\
1.02041 298\\
1.02145 298\\
1.02249 298\\
1.02354 298\\
1.02459 298\\
1.02564 298\\
1.02669 298\\
1.02775 298\\
1.02881 298\\
1.02987 298\\
1.03093 298\\
1.03199 298\\
1.03306 298\\
1.03413 299\\
1.0352 299\\
1.03627 299\\
1.03734 299\\
1.03842 299\\
1.0395 299\\
1.04058 299\\
1.04167 299\\
1.04275 299\\
1.04384 299\\
1.04493 299\\
1.04603 308\\
1.04712 313\\
1.04822 317\\
1.04932 321\\
1.05042 321\\
1.05152 325\\
1.05263 327\\
1.05374 329\\
1.05485 331\\
1.05597 331\\
1.05708 335\\
1.0582 335\\
1.05932 336\\
1.06045 336\\
1.06157 338\\
1.0627 339\\
1.06383 340\\
1.06496 340\\
1.0661 341\\
1.06724 342\\
1.06838 342\\
1.06952 342\\
1.07066 342\\
1.07181 342\\
1.07296 342\\
1.07411 343\\
1.07527 343\\
1.07643 345\\
1.07759 346\\
1.07875 346\\
1.07991 346\\
1.08108 346\\
1.08225 346\\
1.08342 346\\
1.0846 346\\
1.08578 346\\
1.08696 346\\
1.08814 347\\
1.08932 347\\
1.09051 347\\
1.0917 347\\
1.0929 348\\
1.09409 348\\
1.09529 348\\
1.09649 348\\
1.09769 348\\
1.0989 348\\
1.10011 348\\
1.10132 348\\
1.10254 348\\
1.10375 348\\
1.10497 348\\
1.10619 348\\
1.10742 348\\
1.10865 348\\
1.10988 348\\
1.11111 348\\
1.11235 348\\
1.11359 348\\
1.11483 348\\
1.11607 348\\
1.11732 348\\
1.11857 348\\
1.11982 348\\
1.12108 348\\
1.12233 348\\
1.1236 353\\
1.12486 360\\
1.12613 364\\
1.1274 366\\
1.12867 371\\
1.12994 372\\
1.13122 375\\
1.1325 375\\
1.13379 378\\
1.13507 379\\
1.13636 380\\
1.13766 381\\
1.13895 381\\
1.14025 382\\
1.14155 384\\
1.14286 385\\
1.14416 386\\
1.14548 386\\
1.14679 387\\
1.14811 387\\
1.14943 388\\
1.15075 389\\
1.15207 390\\
1.1534 390\\
1.15473 391\\
1.15607 391\\
1.15741 392\\
1.15875 392\\
1.16009 392\\
1.16144 392\\
1.16279 392\\
1.16414 392\\
1.1655 393\\
1.16686 393\\
1.16822 393\\
1.16959 393\\
1.17096 394\\
1.17233 395\\
1.17371 395\\
1.17509 396\\
1.17647 396\\
1.17786 396\\
1.17925 396\\
1.18064 396\\
1.18203 396\\
1.18343 396\\
1.18483 396\\
1.18624 397\\
1.18765 397\\
1.18906 397\\
1.19048 397\\
1.1919 397\\
1.19332 397\\
1.19474 397\\
1.19617 397\\
1.1976 397\\
1.19904 397\\
1.20048 397\\
1.20192 397\\
1.20337 398\\
1.20482 398\\
1.20627 398\\
1.20773 398\\
1.20919 398\\
1.21065 398\\
1.21212 398\\
1.21359 398\\
1.21507 398\\
1.21655 398\\
1.21803 398\\
1.21951 403\\
1.221 409\\
1.22249 412\\
1.22399 415\\
1.22549 417\\
1.22699 419\\
1.2285 420\\
1.23001 422\\
1.23153 424\\
1.23305 425\\
1.23457 425\\
1.23609 427\\
1.23762 429\\
1.23916 429\\
1.24069 430\\
1.24224 430\\
1.24378 432\\
1.24533 432\\
1.24688 433\\
1.24844 435\\
1.25 435\\
1.25156 435\\
1.25313 436\\
1.25471 437\\
1.25628 437\\
1.25786 439\\
1.25945 439\\
1.26103 439\\
1.26263 441\\
1.26422 441\\
1.26582 441\\
1.26743 441\\
1.26904 441\\
1.27065 441\\
1.27226 442\\
1.27389 442\\
1.27551 442\\
1.27714 442\\
1.27877 443\\
1.28041 443\\
1.28205 443\\
1.2837 443\\
1.28535 443\\
1.287 444\\
1.28866 444\\
1.29032 444\\
1.29199 444\\
1.29366 444\\
1.29534 445\\
1.29702 445\\
1.2987 446\\
1.30039 446\\
1.30208 446\\
1.30378 446\\
1.30548 447\\
1.30719 447\\
1.3089 447\\
1.31062 447\\
1.31234 447\\
1.31406 447\\
1.31579 447\\
1.31752 447\\
1.31926 447\\
1.321 447\\
1.32275 447\\
1.3245 447\\
1.32626 447\\
1.32802 447\\
1.32979 447\\
1.33156 447\\
1.33333 447\\
1.33511 453\\
1.3369 458\\
1.33869 460\\
1.34048 464\\
1.34228 465\\
1.34409 467\\
1.3459 468\\
1.34771 470\\
1.34953 471\\
1.35135 473\\
1.35318 475\\
1.35501 475\\
1.35685 476\\
1.3587 477\\
1.36054 479\\
1.3624 479\\
1.36426 479\\
1.36612 480\\
1.36799 482\\
1.36986 482\\
1.37174 483\\
1.37363 484\\
1.37552 485\\
1.37741 486\\
1.37931 486\\
1.38122 486\\
1.38313 486\\
1.38504 486\\
1.38696 487\\
1.38889 487\\
1.39082 489\\
1.39276 489\\
1.3947 489\\
1.39665 490\\
1.3986 490\\
1.40056 490\\
1.40252 490\\
1.40449 490\\
1.40647 490\\
1.40845 492\\
1.41044 493\\
1.41243 493\\
1.41443 493\\
1.41643 493\\
1.41844 493\\
1.42045 493\\
1.42248 494\\
1.4245 494\\
1.42653 494\\
1.42857 494\\
};

\addplot [
color=black,
solid,
line width=1.0pt
]
table[row sep=crcr]{
0.769231 484\\
0.769823 484\\
0.770416 482\\
0.77101 481\\
0.771605 481\\
0.772201 479\\
0.772798 479\\
0.773395 479\\
0.773994 479\\
0.774593 479\\
0.775194 479\\
0.775795 478\\
0.776398 478\\
0.777001 478\\
0.777605 478\\
0.77821 476\\
0.778816 476\\
0.779423 475\\
0.780031 475\\
0.78064 474\\
0.78125 474\\
0.781861 474\\
0.782473 474\\
0.783085 473\\
0.783699 473\\
0.784314 471\\
0.784929 471\\
0.785546 469\\
0.786164 470\\
0.786782 469\\
0.787402 469\\
0.788022 467\\
0.788644 467\\
0.789266 464\\
0.789889 465\\
0.790514 464\\
0.791139 464\\
0.791766 461\\
0.792393 460\\
0.793021 460\\
0.793651 457\\
0.794281 456\\
0.794913 453\\
0.795545 449\\
0.796178 447\\
0.796813 447\\
0.797448 447\\
0.798085 447\\
0.798722 447\\
0.799361 447\\
0.8 447\\
0.800641 447\\
0.801282 447\\
0.801925 447\\
0.802568 447\\
0.803213 447\\
0.803859 447\\
0.804505 447\\
0.805153 447\\
0.805802 447\\
0.806452 447\\
0.807103 447\\
0.807754 447\\
0.808407 447\\
0.809061 447\\
0.809717 447\\
0.810373 447\\
0.81103 447\\
0.811688 447\\
0.812348 447\\
0.813008 447\\
0.81367 447\\
0.814332 447\\
0.814996 447\\
0.815661 447\\
0.816327 447\\
0.816993 446\\
0.817661 447\\
0.818331 447\\
0.819001 447\\
0.819672 447\\
0.820345 447\\
0.821018 447\\
0.821693 446\\
0.822368 446\\
0.823045 445\\
0.823723 445\\
0.824402 445\\
0.825083 445\\
0.825764 445\\
0.826446 445\\
0.82713 445\\
0.827815 445\\
0.8285 445\\
0.829187 444\\
0.829876 444\\
0.830565 444\\
0.831255 444\\
0.831947 444\\
0.832639 444\\
0.833333 444\\
0.834028 444\\
0.834725 444\\
0.835422 444\\
0.83612 443\\
0.83682 443\\
0.837521 443\\
0.838223 443\\
0.838926 442\\
0.839631 442\\
0.840336 442\\
0.841043 442\\
0.841751 442\\
0.84246 442\\
0.84317 442\\
0.843882 442\\
0.844595 441\\
0.845309 442\\
0.846024 441\\
0.84674 441\\
0.847458 441\\
0.848176 441\\
0.848896 441\\
0.849618 441\\
0.85034 441\\
0.851064 441\\
0.851789 441\\
0.852515 441\\
0.853242 440\\
0.853971 440\\
0.854701 439\\
0.855432 439\\
0.856164 438\\
0.856898 438\\
0.857633 437\\
0.858369 436\\
0.859107 436\\
0.859845 436\\
0.860585 437\\
0.861326 436\\
0.862069 435\\
0.862813 435\\
0.863558 435\\
0.864304 435\\
0.865052 433\\
0.865801 434\\
0.866551 434\\
0.867303 434\\
0.868056 433\\
0.86881 432\\
0.869565 429\\
0.870322 430\\
0.87108 428\\
0.87184 429\\
0.8726 428\\
0.873362 427\\
0.874126 425\\
0.874891 424\\
0.875657 424\\
0.876424 423\\
0.877193 421\\
0.877963 421\\
0.878735 419\\
0.879507 420\\
0.880282 418\\
0.881057 414\\
0.881834 414\\
0.882613 411\\
0.883392 408\\
0.884173 411\\
0.884956 408\\
0.88574 406\\
0.886525 400\\
0.887311 398\\
0.888099 398\\
0.888889 398\\
0.88968 398\\
0.890472 397\\
0.891266 397\\
0.892061 397\\
0.892857 397\\
0.893655 397\\
0.894454 398\\
0.895255 398\\
0.896057 397\\
0.896861 397\\
0.897666 397\\
0.898473 397\\
0.899281 397\\
0.90009 397\\
0.900901 397\\
0.901713 397\\
0.902527 397\\
0.903342 397\\
0.904159 397\\
0.904977 397\\
0.905797 397\\
0.906618 397\\
0.907441 397\\
0.908265 397\\
0.909091 397\\
0.909918 397\\
0.910747 397\\
0.911577 396\\
0.912409 396\\
0.913242 395\\
0.914077 396\\
0.914913 395\\
0.915751 396\\
0.91659 395\\
0.917431 396\\
0.918274 395\\
0.919118 394\\
0.919963 394\\
0.92081 393\\
0.921659 393\\
0.922509 393\\
0.923361 392\\
0.924214 392\\
0.925069 391\\
0.925926 391\\
0.926784 390\\
0.927644 390\\
0.928505 390\\
0.929368 390\\
0.930233 389\\
0.931099 389\\
0.931966 387\\
0.932836 387\\
0.933707 387\\
0.934579 387\\
0.935454 386\\
0.93633 384\\
0.937207 386\\
0.938086 384\\
0.938967 384\\
0.93985 381\\
0.940734 381\\
0.94162 379\\
0.942507 379\\
0.943396 378\\
0.944287 374\\
0.94518 373\\
0.946074 375\\
0.94697 371\\
0.947867 368\\
0.948767 365\\
0.949668 362\\
0.95057 355\\
0.951475 355\\
0.952381 358\\
0.953289 357\\
0.954198 357\\
0.95511 353\\
0.956023 348\\
0.956938 348\\
0.957854 348\\
0.958773 348\\
0.959693 348\\
0.960615 348\\
0.961538 348\\
0.962464 348\\
0.963391 348\\
0.96432 348\\
0.965251 348\\
0.966184 348\\
0.967118 348\\
0.968054 348\\
0.968992 348\\
0.969932 348\\
0.970874 348\\
0.971817 348\\
0.972763 348\\
0.97371 348\\
0.974659 348\\
0.97561 348\\
0.976562 348\\
0.977517 348\\
0.978474 348\\
0.979432 348\\
0.980392 348\\
0.981354 348\\
0.982318 347\\
0.983284 347\\
0.984252 347\\
0.985222 347\\
0.986193 347\\
0.987167 346\\
0.988142 346\\
0.98912 346\\
0.990099 346\\
0.99108 346\\
0.992063 346\\
0.993049 346\\
0.994036 346\\
0.995025 346\\
0.996016 346\\
0.997009 346\\
0.998004 346\\
0.999001 346\\
1 344\\
1.001 343\\
1.002 342\\
1.00301 342\\
1.00402 342\\
1.00503 342\\
1.00604 342\\
1.00705 342\\
1.00806 340\\
1.00908 337\\
1.0101 338\\
1.01112 338\\
1.01215 337\\
1.01317 336\\
1.0142 331\\
1.01523 329\\
1.01626 326\\
1.01729 325\\
1.01833 320\\
1.01937 315\\
1.02041 304\\
1.02145 299\\
1.02249 299\\
1.02354 299\\
1.02459 299\\
1.02564 298\\
1.02669 298\\
1.02775 298\\
1.02881 298\\
1.02987 298\\
1.03093 298\\
1.03199 298\\
1.03306 298\\
1.03413 298\\
1.0352 298\\
1.03627 298\\
1.03734 298\\
1.03842 297\\
1.0395 297\\
1.04058 296\\
1.04167 296\\
1.04275 296\\
1.04384 296\\
1.04493 296\\
1.04603 296\\
1.04712 296\\
1.04822 296\\
1.04932 296\\
1.05042 297\\
1.05152 298\\
1.05263 298\\
1.05374 298\\
1.05485 298\\
1.05597 298\\
1.05708 298\\
1.0582 298\\
1.05932 298\\
1.06045 298\\
1.06157 298\\
1.0627 298\\
1.06383 298\\
1.06496 298\\
1.0661 298\\
1.06724 299\\
1.06838 299\\
1.06952 299\\
1.07066 299\\
1.07181 299\\
1.07296 299\\
1.07411 299\\
1.07527 299\\
1.07643 299\\
1.07759 299\\
1.07875 303\\
1.07991 313\\
1.08108 316\\
1.08225 319\\
1.08342 325\\
1.0846 325\\
1.08578 329\\
1.08696 331\\
1.08814 332\\
1.08932 335\\
1.09051 336\\
1.0917 338\\
1.0929 338\\
1.09409 340\\
1.09529 341\\
1.09649 342\\
1.09769 342\\
1.0989 342\\
1.10011 342\\
1.10132 342\\
1.10254 343\\
1.10375 346\\
1.10497 346\\
1.10619 346\\
1.10742 346\\
1.10865 346\\
1.10988 346\\
1.11111 346\\
1.11235 346\\
1.11359 347\\
1.11483 347\\
1.11607 348\\
1.11732 348\\
1.11857 348\\
1.11982 348\\
1.12108 348\\
1.12233 348\\
1.1236 348\\
1.12486 348\\
1.12613 348\\
1.1274 348\\
1.12867 348\\
1.12994 348\\
1.13122 348\\
1.1325 348\\
1.13379 348\\
1.13507 348\\
1.13636 348\\
1.13766 348\\
1.13895 348\\
1.14025 349\\
1.14155 360\\
1.14286 364\\
1.14416 368\\
1.14548 371\\
1.14679 371\\
1.14811 375\\
1.14943 375\\
1.15075 377\\
1.15207 378\\
1.1534 379\\
1.15473 380\\
1.15607 381\\
1.15741 381\\
1.15875 384\\
1.16009 385\\
1.16144 386\\
1.16279 387\\
1.16414 387\\
1.1655 388\\
1.16686 389\\
1.16822 389\\
1.16959 390\\
1.17096 390\\
1.17233 391\\
1.17371 391\\
1.17509 391\\
1.17647 392\\
1.17786 392\\
1.17925 392\\
1.18064 392\\
1.18203 392\\
1.18343 392\\
1.18483 392\\
1.18624 393\\
1.18765 393\\
1.18906 394\\
1.19048 394\\
1.1919 394\\
1.19332 395\\
1.19474 396\\
1.19617 396\\
1.1976 396\\
1.19904 396\\
1.20048 396\\
1.20192 396\\
1.20337 396\\
1.20482 396\\
1.20627 396\\
1.20773 397\\
1.20919 397\\
1.21065 397\\
1.21212 397\\
1.21359 397\\
1.21507 397\\
1.21655 397\\
1.21803 397\\
1.21951 397\\
1.221 397\\
1.22249 397\\
1.22399 397\\
1.22549 397\\
1.22699 398\\
1.2285 398\\
1.23001 398\\
1.23153 398\\
1.23305 404\\
1.23457 409\\
1.23609 404\\
1.23762 410\\
1.23916 407\\
1.24069 406\\
1.24224 403\\
1.24378 410\\
1.24533 414\\
1.24688 420\\
1.24844 424\\
1.25 425\\
1.25156 429\\
1.25313 429\\
1.25471 432\\
1.25628 431\\
1.25786 432\\
1.25945 434\\
1.26103 436\\
1.26263 438\\
1.26422 438\\
1.26582 439\\
1.26743 439\\
1.26904 439\\
1.27065 438\\
1.27226 440\\
1.27389 441\\
1.27551 441\\
1.27714 442\\
1.27877 441\\
1.28041 442\\
1.28205 443\\
1.2837 443\\
1.28535 442\\
1.287 443\\
1.28866 444\\
1.29032 444\\
1.29199 444\\
1.29366 444\\
1.29534 444\\
1.29702 445\\
1.2987 445\\
1.30039 445\\
1.30208 446\\
1.30378 446\\
1.30548 446\\
1.30719 446\\
1.3089 446\\
1.31062 446\\
1.31234 447\\
1.31406 447\\
1.31579 447\\
1.31752 447\\
1.31926 447\\
1.321 447\\
1.32275 447\\
1.3245 447\\
1.32626 447\\
1.32802 447\\
1.32979 447\\
1.33156 447\\
1.33333 447\\
1.33511 447\\
1.3369 447\\
1.33869 452\\
1.34048 456\\
1.34228 460\\
1.34409 460\\
1.3459 461\\
1.34771 463\\
1.34953 467\\
1.35135 467\\
1.35318 469\\
1.35501 471\\
1.35685 472\\
1.3587 473\\
1.36054 474\\
1.3624 476\\
1.36426 476\\
1.36612 476\\
1.36799 477\\
1.36986 479\\
1.37174 480\\
1.37363 481\\
1.37552 481\\
1.37741 480\\
1.37931 483\\
1.38122 484\\
1.38313 485\\
1.38504 485\\
1.38696 486\\
1.38889 486\\
1.39082 486\\
1.39276 486\\
1.3947 486\\
1.39665 488\\
1.3986 489\\
1.40056 488\\
1.40252 490\\
1.40449 490\\
1.40647 490\\
1.40845 490\\
1.41044 490\\
1.41243 491\\
1.41443 491\\
1.41643 491\\
1.41844 492\\
1.42045 493\\
1.42248 493\\
1.4245 493\\
1.42653 493\\
1.42857 493\\
};

\addplot [
color=black,
solid,
line width=1.0pt
]
table[row sep=crcr]{
0.769231 472\\
0.769823 471\\
0.770416 470\\
0.77101 471\\
0.771605 473\\
0.772201 473\\
0.772798 471\\
0.773395 470\\
0.773994 467\\
0.774593 467\\
0.775194 467\\
0.775795 466\\
0.776398 468\\
0.777001 464\\
0.777605 462\\
0.77821 462\\
0.778816 462\\
0.779423 460\\
0.780031 458\\
0.78064 455\\
0.78125 453\\
0.781861 449\\
0.782473 448\\
0.783085 448\\
0.783699 449\\
0.784314 449\\
0.784929 448\\
0.785546 448\\
0.786164 447\\
0.786782 447\\
0.787402 447\\
0.788022 447\\
0.788644 447\\
0.789266 447\\
0.789889 447\\
0.790514 447\\
0.791139 447\\
0.791766 447\\
0.792393 447\\
0.793021 447\\
0.793651 447\\
0.794281 447\\
0.794913 447\\
0.795545 447\\
0.796178 447\\
0.796813 447\\
0.797448 447\\
0.798085 447\\
0.798722 447\\
0.799361 447\\
0.8 447\\
0.800641 446\\
0.801282 447\\
0.801925 447\\
0.802568 446\\
0.803213 447\\
0.803859 446\\
0.804505 446\\
0.805153 446\\
0.805802 446\\
0.806452 446\\
0.807103 446\\
0.807754 445\\
0.808407 446\\
0.809061 446\\
0.809717 446\\
0.810373 445\\
0.81103 445\\
0.811688 445\\
0.812348 444\\
0.813008 444\\
0.81367 445\\
0.814332 445\\
0.814996 445\\
0.815661 444\\
0.816327 444\\
0.816993 444\\
0.817661 444\\
0.818331 443\\
0.819001 443\\
0.819672 444\\
0.820345 444\\
0.821018 442\\
0.821693 443\\
0.822368 443\\
0.823045 443\\
0.823723 442\\
0.824402 442\\
0.825083 442\\
0.825764 442\\
0.826446 441\\
0.82713 441\\
0.827815 441\\
0.8285 441\\
0.829187 441\\
0.829876 441\\
0.830565 441\\
0.831255 441\\
0.831947 441\\
0.832639 441\\
0.833333 441\\
0.834028 441\\
0.834725 441\\
0.835422 440\\
0.83612 440\\
0.83682 439\\
0.837521 440\\
0.838223 439\\
0.838926 437\\
0.839631 438\\
0.840336 437\\
0.841043 436\\
0.841751 436\\
0.84246 435\\
0.84317 436\\
0.843882 436\\
0.844595 435\\
0.845309 434\\
0.846024 434\\
0.84674 434\\
0.847458 434\\
0.848176 435\\
0.848896 434\\
0.849618 433\\
0.85034 432\\
0.851064 433\\
0.851789 432\\
0.852515 433\\
0.853242 431\\
0.853971 430\\
0.854701 430\\
0.855432 429\\
0.856164 428\\
0.856898 428\\
0.857633 428\\
0.858369 427\\
0.859107 425\\
0.859845 426\\
0.860585 426\\
0.861326 425\\
0.862069 424\\
0.862813 423\\
0.863558 422\\
0.864304 422\\
0.865052 419\\
0.865801 420\\
0.866551 418\\
0.867303 419\\
0.868056 417\\
0.86881 416\\
0.869565 415\\
0.870322 415\\
0.87108 412\\
0.87184 411\\
0.8726 409\\
0.873362 408\\
0.874126 403\\
0.874891 403\\
0.875657 398\\
0.876424 398\\
0.877193 398\\
0.877963 398\\
0.878735 398\\
0.879507 398\\
0.880282 398\\
0.881057 398\\
0.881834 398\\
0.882613 398\\
0.883392 398\\
0.884173 398\\
0.884956 398\\
0.88574 398\\
0.886525 398\\
0.887311 398\\
0.888099 397\\
0.888889 397\\
0.88968 397\\
0.890472 397\\
0.891266 397\\
0.892061 397\\
0.892857 397\\
0.893655 397\\
0.894454 397\\
0.895255 397\\
0.896057 397\\
0.896861 397\\
0.897666 397\\
0.898473 397\\
0.899281 397\\
0.90009 397\\
0.900901 396\\
0.901713 396\\
0.902527 396\\
0.903342 397\\
0.904159 396\\
0.904977 397\\
0.905797 396\\
0.906618 396\\
0.907441 396\\
0.908265 394\\
0.909091 394\\
0.909918 394\\
0.910747 394\\
0.911577 393\\
0.912409 393\\
0.913242 393\\
0.914077 392\\
0.914913 392\\
0.915751 392\\
0.91659 392\\
0.917431 392\\
0.918274 392\\
0.919118 392\\
0.919963 392\\
0.92081 391\\
0.921659 391\\
0.922509 391\\
0.923361 391\\
0.924214 391\\
0.925069 390\\
0.925926 391\\
0.926784 390\\
0.927644 389\\
0.928505 387\\
0.929368 387\\
0.930233 387\\
0.931099 387\\
0.931966 386\\
0.932836 385\\
0.933707 386\\
0.934579 385\\
0.935454 385\\
0.93633 382\\
0.937207 381\\
0.938086 381\\
0.938967 380\\
0.93985 378\\
0.940734 378\\
0.94162 376\\
0.942507 377\\
0.943396 371\\
0.944287 372\\
0.94518 370\\
0.946074 367\\
0.94697 367\\
0.947867 364\\
0.948767 358\\
0.949668 350\\
0.95057 348\\
0.951475 348\\
0.952381 348\\
0.953289 348\\
0.954198 348\\
0.95511 348\\
0.956023 348\\
0.956938 348\\
0.957854 348\\
0.958773 348\\
0.959693 348\\
0.960615 348\\
0.961538 348\\
0.962464 348\\
0.963391 348\\
0.96432 348\\
0.965251 348\\
0.966184 348\\
0.967118 348\\
0.968054 348\\
0.968992 348\\
0.969932 348\\
0.970874 348\\
0.971817 348\\
0.972763 348\\
0.97371 348\\
0.974659 348\\
0.97561 348\\
0.976562 348\\
0.977517 348\\
0.978474 348\\
0.979432 347\\
0.980392 346\\
0.981354 346\\
0.982318 346\\
0.983284 347\\
0.984252 347\\
0.985222 346\\
0.986193 346\\
0.987167 346\\
0.988142 346\\
0.98912 346\\
0.990099 344\\
0.99108 344\\
0.992063 343\\
0.993049 342\\
0.994036 342\\
0.995025 342\\
0.996016 342\\
0.997009 342\\
0.998004 342\\
0.999001 342\\
1 341\\
1.001 340\\
1.002 340\\
1.00301 340\\
1.00402 338\\
1.00503 338\\
1.00604 337\\
1.00705 336\\
1.00806 336\\
1.00908 335\\
1.0101 333\\
1.01112 332\\
1.01215 330\\
1.01317 329\\
1.0142 329\\
1.01523 326\\
1.01626 323\\
1.01729 320\\
1.01833 313\\
1.01937 301\\
1.02041 299\\
1.02145 299\\
1.02249 299\\
1.02354 299\\
1.02459 298\\
1.02564 298\\
1.02669 298\\
1.02775 298\\
1.02881 298\\
1.02987 298\\
1.03093 298\\
1.03199 298\\
1.03306 298\\
1.03413 298\\
1.0352 298\\
1.03627 298\\
1.03734 298\\
1.03842 298\\
1.0395 298\\
1.04058 298\\
1.04167 298\\
1.04275 298\\
1.04384 297\\
1.04493 297\\
1.04603 298\\
1.04712 297\\
1.04822 297\\
1.04932 297\\
1.05042 297\\
1.05152 297\\
1.05263 296\\
1.05374 297\\
1.05485 297\\
1.05597 297\\
1.05708 297\\
1.0582 298\\
1.05932 298\\
1.06045 298\\
1.06157 298\\
1.0627 298\\
1.06383 298\\
1.06496 298\\
1.0661 298\\
1.06724 297\\
1.06838 298\\
1.06952 298\\
1.07066 298\\
1.07181 298\\
1.07296 298\\
1.07411 298\\
1.07527 298\\
1.07643 298\\
1.07759 298\\
1.07875 298\\
1.07991 298\\
1.08108 298\\
1.08225 298\\
1.08342 298\\
1.0846 298\\
1.08578 298\\
1.08696 298\\
1.08814 298\\
1.08932 299\\
1.09051 299\\
1.0917 299\\
1.0929 299\\
1.09409 299\\
1.09529 299\\
1.09649 299\\
1.09769 299\\
1.0989 299\\
1.10011 299\\
1.10132 303\\
1.10254 313\\
1.10375 320\\
1.10497 323\\
1.10619 325\\
1.10742 327\\
1.10865 330\\
1.10988 329\\
1.11111 334\\
1.11235 335\\
1.11359 335\\
1.11483 336\\
1.11607 339\\
1.11732 340\\
1.11857 340\\
1.11982 341\\
1.12108 342\\
1.12233 343\\
1.1236 345\\
1.12486 345\\
1.12613 345\\
1.1274 345\\
1.12867 343\\
1.12994 343\\
1.13122 342\\
1.1325 345\\
1.13379 345\\
1.13507 346\\
1.13636 346\\
1.13766 347\\
1.13895 347\\
1.14025 348\\
1.14155 348\\
1.14286 348\\
1.14416 348\\
1.14548 348\\
1.14679 348\\
1.14811 348\\
1.14943 350\\
1.15075 350\\
1.15207 350\\
1.1534 349\\
1.15473 348\\
1.15607 348\\
1.15741 351\\
1.15875 349\\
1.16009 348\\
1.16144 349\\
1.16279 349\\
1.16414 350\\
1.1655 349\\
1.16686 349\\
1.16822 350\\
1.16959 349\\
1.17096 349\\
1.17233 348\\
1.17371 349\\
1.17509 349\\
1.17647 348\\
1.17786 349\\
1.17925 349\\
1.18064 357\\
1.18203 361\\
1.18343 364\\
1.18483 366\\
1.18624 370\\
1.18765 372\\
1.18906 374\\
1.19048 376\\
1.1919 375\\
1.19332 380\\
1.19474 381\\
1.19617 381\\
1.1976 383\\
1.19904 385\\
1.20048 386\\
1.20192 387\\
1.20337 387\\
1.20482 387\\
1.20627 388\\
1.20773 388\\
1.20919 389\\
1.21065 390\\
1.21212 390\\
1.21359 390\\
1.21507 391\\
1.21655 392\\
1.21803 392\\
1.21951 392\\
1.221 392\\
1.22249 393\\
1.22399 392\\
1.22549 393\\
1.22699 393\\
1.2285 393\\
1.23001 394\\
1.23153 394\\
1.23305 394\\
1.23457 395\\
1.23609 395\\
1.23762 396\\
1.23916 397\\
1.24069 396\\
1.24224 396\\
1.24378 396\\
1.24533 396\\
1.24688 396\\
1.24844 396\\
1.25 396\\
1.25156 397\\
1.25313 397\\
1.25471 397\\
1.25628 397\\
1.25786 402\\
1.25945 408\\
1.26103 411\\
1.26263 414\\
1.26422 415\\
1.26582 416\\
1.26743 418\\
1.26904 420\\
1.27065 422\\
1.27226 424\\
1.27389 423\\
1.27551 425\\
1.27714 422\\
1.27877 419\\
1.28041 416\\
1.28205 415\\
1.2837 419\\
1.28535 426\\
1.287 425\\
1.28866 427\\
1.29032 433\\
1.29199 431\\
1.29366 430\\
1.29534 436\\
1.29702 434\\
1.2987 436\\
1.30039 436\\
1.30208 437\\
1.30378 438\\
1.30548 438\\
1.30719 437\\
1.3089 436\\
1.31062 441\\
1.31234 441\\
1.31406 440\\
1.31579 440\\
1.31752 441\\
1.31926 441\\
1.321 441\\
1.32275 442\\
1.3245 442\\
1.32626 442\\
1.32802 443\\
1.32979 443\\
1.33156 443\\
1.33333 443\\
1.33511 443\\
1.3369 443\\
1.33869 443\\
1.34048 444\\
1.34228 444\\
1.34409 444\\
1.3459 445\\
1.34771 444\\
1.34953 444\\
1.35135 445\\
1.35318 445\\
1.35501 445\\
1.35685 447\\
1.3587 447\\
1.36054 446\\
1.3624 447\\
1.36426 447\\
1.36612 447\\
1.36799 447\\
1.36986 447\\
1.37174 447\\
1.37363 447\\
1.37552 447\\
1.37741 447\\
1.37931 447\\
1.38122 447\\
1.38313 447\\
1.38504 447\\
1.38696 447\\
1.38889 447\\
1.39082 447\\
1.39276 447\\
1.3947 456\\
1.39665 456\\
1.3986 464\\
1.40056 462\\
1.40252 469\\
1.40449 469\\
1.40647 472\\
1.40845 472\\
1.41044 473\\
1.41243 474\\
1.41443 475\\
1.41643 476\\
1.41844 479\\
1.42045 478\\
1.42248 478\\
1.4245 479\\
1.42653 480\\
1.42857 482\\
};

\addplot [
color=black,
only marks,
mark=x,
mark size = 2.5,
mark options={solid,blue}
]
table[row sep=crcr]{
0.769231 472\\
0.775194 464\\
0.78125 447\\
0.787402 447\\
0.793651 447\\
0.8 447\\
0.806452 445\\
0.813008 444\\
0.819672 443\\
0.826446 441\\
0.833333 439\\
0.840336 435\\
0.847458 430\\
0.854701 424\\
0.862069 412\\
0.869565 398\\
0.877193 397\\
0.884956 397\\
0.892857 396\\
0.900901 395\\
0.909091 392\\
0.917431 390\\
0.925926 386\\
0.934579 376\\
0.943396 348\\
0.952381 348\\
0.961538 348\\
0.970874 347\\
0.980392 346\\
0.990099 338\\
1 325\\
1.0101 298\\
1.02041 298\\
1.03093 298\\
1.04167 299\\
1.05263 327\\
1.06383 340\\
1.07527 343\\
1.08696 346\\
1.0989 348\\
1.11111 348\\
1.1236 353\\
1.13636 380\\
1.14943 388\\
1.16279 392\\
1.17647 396\\
1.19048 397\\
1.20482 398\\
1.21951 403\\
1.23457 425\\
1.25 435\\
1.26582 441\\
1.28205 443\\
1.2987 446\\
1.31579 447\\
1.33333 447\\
1.35135 473\\
1.36986 482\\
1.38889 487\\
1.40845 492\\
1.42857 494\\
};

\addplot [
color=black,
only marks,
mark=o,
mark size = 2.5,
mark options={solid,green}
]
table[row sep=crcr]{
0.769231 484\\
0.775194 479\\
0.78125 474\\
0.787402 469\\
0.793651 457\\
0.8 447\\
0.806452 447\\
0.813008 447\\
0.819672 447\\
0.826446 445\\
0.833333 444\\
0.840336 442\\
0.847458 441\\
0.854701 439\\
0.862069 435\\
0.869565 429\\
0.877193 421\\
0.884956 408\\
0.892857 397\\
0.900901 397\\
0.909091 397\\
0.917431 396\\
0.925926 391\\
0.934579 387\\
0.943396 378\\
0.952381 358\\
0.961538 348\\
0.970874 348\\
0.980392 348\\
0.990099 346\\
1 344\\
1.0101 338\\
1.02041 304\\
1.03093 298\\
1.04167 296\\
1.05263 298\\
1.06383 298\\
1.07527 299\\
1.08696 331\\
1.0989 342\\
1.11111 346\\
1.1236 348\\
1.13636 348\\
1.14943 375\\
1.16279 387\\
1.17647 392\\
1.19048 394\\
1.20482 396\\
1.21951 397\\
1.23457 409\\
1.25 425\\
1.26582 439\\
1.28205 443\\
1.2987 445\\
1.31579 447\\
1.33333 447\\
1.35135 467\\
1.36986 479\\
1.38889 486\\
1.40845 490\\
1.42857 493\\
};

\addplot [
color=black,
only marks,
mark=+,
mark size = 2.5,
mark options={solid,red}
]
table[row sep=crcr]{
0.769231 472\\
0.775194 467\\
0.78125 453\\
0.787402 447\\
0.793651 447\\
0.8 447\\
0.806452 446\\
0.813008 444\\
0.819672 444\\
0.826446 441\\
0.833333 441\\
0.840336 437\\
0.847458 434\\
0.854701 430\\
0.862069 424\\
0.869565 415\\
0.877193 398\\
0.884956 398\\
0.892857 397\\
0.900901 396\\
0.909091 394\\
0.917431 392\\
0.925926 391\\
0.934579 385\\
0.943396 371\\
0.952381 348\\
0.961538 348\\
0.970874 348\\
0.980392 346\\
0.990099 344\\
1 341\\
1.0101 333\\
1.02041 299\\
1.03093 298\\
1.04167 298\\
1.05263 296\\
1.06383 298\\
1.07527 298\\
1.08696 298\\
1.0989 299\\
1.11111 334\\
1.1236 345\\
1.13636 346\\
1.14943 350\\
1.16279 349\\
1.17647 348\\
1.19048 376\\
1.20482 387\\
1.21951 392\\
1.23457 395\\
1.25 396\\
1.26582 416\\
1.28205 415\\
1.2987 436\\
1.31579 440\\
1.33333 443\\
1.35135 445\\
1.36986 447\\
1.38889 447\\
1.40845 472\\
1.42857 482\\
};

\end{axis}
\end{tikzpicture}
\caption{Total fixed-stress iterations for varying mesh size $h = 2^{-(m+1)}$
for dG(0) in time.}
\label{fig:2:dG0:h}
\end{figure}

%% file: p-independent.tex
\begin{figure}[p]
\centering
%
\begin{tikzpicture}

\begin{axis}[%
width=4.82222222222222in,
height=3.5in,
scale only axis,
/pgf/number format/.cd, 1000 sep={},
xlabel={$\omega$},
ylabel={iterations},
xmin=0.82,
xmax=1.3,
ymin=290,
ymax=450,
yminorticks=true,
legend style={legend pos=north east},
legend style={draw=black,fill=white,legend cell align=left},
legend entries = {
  {cGP(1), $m=2$, $s=1$, $\tau_n=0.01$},
  {cGP(1), $m=2$, $s=2$, $\tau_n=0.01$},
  {cGP(1), $m=2$, $s=3$, $\tau_n=0.01$},
  {cGP(1), $m=2$, $s=4$, $\tau_n=0.01$}
}
]

\addplot [
color=black,
solid,
line width=1.0pt,
mark=x,
mark size = 2.5,
mark options={solid,blue}
]
table[row sep=crcr]{
0.0 0\\
0.1 0\\
};

\addplot [
color=black,
solid,
line width=1.0pt,
mark=o,
mark size = 2.5,
mark options={solid,green}
]
table[row sep=crcr]{
0.0 0\\
0.1 0\\
};

\addplot [
color=black,
solid,
line width=1.0pt,
mark=+,
mark size = 2.5,
mark options={solid,red}
]
table[row sep=crcr]{
0.0 0\\
0.1 0\\
};

\addplot [
color=black,
solid,
line width=1.0pt,
mark=star,
mark size = 2.5,
mark options={solid,orange}
]
table[row sep=crcr]{
0.0 0\\
0.1 0\\
};

\addplot [
color=black,
dotted,
line width=1.0pt
]
table[row sep=crcr]{
1.0 0\\
1.0 2000\\
};

\addplot [
color=black,
solid,
line width=1.0pt
]
table[row sep=crcr]{
0.769231 447\\
0.769823 447\\
0.770416 447\\
0.77101 447\\
0.771605 446\\
0.772201 446\\
0.772798 446\\
0.773395 446\\
0.773994 446\\
0.774593 446\\
0.775194 446\\
0.775795 446\\
0.776398 446\\
0.777001 446\\
0.777605 446\\
0.77821 446\\
0.778816 446\\
0.779423 446\\
0.780031 446\\
0.78064 446\\
0.78125 446\\
0.781861 446\\
0.782473 446\\
0.783085 446\\
0.783699 446\\
0.784314 446\\
0.784929 446\\
0.785546 446\\
0.786164 446\\
0.786782 446\\
0.787402 446\\
0.788022 446\\
0.788644 445\\
0.789266 444\\
0.789889 444\\
0.790514 444\\
0.791139 444\\
0.791766 444\\
0.792393 444\\
0.793021 444\\
0.793651 444\\
0.794281 444\\
0.794913 444\\
0.795545 444\\
0.796178 444\\
0.796813 444\\
0.797448 444\\
0.798085 444\\
0.798722 444\\
0.799361 444\\
0.8 444\\
0.800641 444\\
0.801282 444\\
0.801925 444\\
0.802568 444\\
0.803213 444\\
0.803859 444\\
0.804505 444\\
0.805153 444\\
0.805802 444\\
0.806452 443\\
0.807103 442\\
0.807754 442\\
0.808407 442\\
0.809061 442\\
0.809717 442\\
0.810373 442\\
0.81103 442\\
0.811688 442\\
0.812348 442\\
0.813008 442\\
0.81367 442\\
0.814332 442\\
0.814996 442\\
0.815661 442\\
0.816327 442\\
0.816993 441\\
0.817661 441\\
0.818331 441\\
0.819001 439\\
0.819672 439\\
0.820345 439\\
0.821018 438\\
0.821693 438\\
0.822368 437\\
0.823045 437\\
0.823723 437\\
0.824402 437\\
0.825083 436\\
0.825764 436\\
0.826446 435\\
0.82713 434\\
0.827815 434\\
0.8285 434\\
0.829187 434\\
0.829876 434\\
0.830565 434\\
0.831255 434\\
0.831947 434\\
0.832639 434\\
0.833333 433\\
0.834028 432\\
0.834725 432\\
0.835422 432\\
0.83612 432\\
0.83682 431\\
0.837521 430\\
0.838223 428\\
0.838926 428\\
0.839631 426\\
0.840336 425\\
0.841043 423\\
0.841751 422\\
0.84246 421\\
0.84317 421\\
0.843882 420\\
0.844595 419\\
0.845309 416\\
0.846024 415\\
0.84674 412\\
0.847458 412\\
0.848176 410\\
0.848896 407\\
0.849618 405\\
0.85034 403\\
0.851064 402\\
0.851789 401\\
0.852515 398\\
0.853242 397\\
0.853971 397\\
0.854701 397\\
0.855432 397\\
0.856164 397\\
0.856898 397\\
0.857633 397\\
0.858369 397\\
0.859107 397\\
0.859845 397\\
0.860585 397\\
0.861326 397\\
0.862069 397\\
0.862813 397\\
0.863558 396\\
0.864304 396\\
0.865052 396\\
0.865801 396\\
0.866551 396\\
0.867303 396\\
0.868056 396\\
0.86881 396\\
0.869565 396\\
0.870322 396\\
0.87108 396\\
0.87184 396\\
0.8726 396\\
0.873362 396\\
0.874126 396\\
0.874891 396\\
0.875657 396\\
0.876424 396\\
0.877193 396\\
0.877963 396\\
0.878735 396\\
0.879507 396\\
0.880282 396\\
0.881057 396\\
0.881834 396\\
0.882613 396\\
0.883392 396\\
0.884173 396\\
0.884956 396\\
0.88574 396\\
0.886525 394\\
0.887311 394\\
0.888099 394\\
0.888889 394\\
0.88968 394\\
0.890472 394\\
0.891266 394\\
0.892061 394\\
0.892857 394\\
0.893655 394\\
0.894454 394\\
0.895255 394\\
0.896057 393\\
0.896861 392\\
0.897666 392\\
0.898473 392\\
0.899281 392\\
0.90009 392\\
0.900901 392\\
0.901713 392\\
0.902527 392\\
0.903342 392\\
0.904159 392\\
0.904977 392\\
0.905797 392\\
0.906618 391\\
0.907441 391\\
0.908265 390\\
0.909091 388\\
0.909918 387\\
0.910747 387\\
0.911577 386\\
0.912409 385\\
0.913242 384\\
0.914077 384\\
0.914913 384\\
0.915751 384\\
0.91659 384\\
0.917431 384\\
0.918274 383\\
0.919118 382\\
0.919963 382\\
0.92081 382\\
0.921659 381\\
0.922509 377\\
0.923361 376\\
0.924214 375\\
0.925069 374\\
0.925926 373\\
0.926784 373\\
0.927644 368\\
0.928505 364\\
0.929368 363\\
0.930233 359\\
0.931099 356\\
0.931966 353\\
0.932836 352\\
0.933707 348\\
0.934579 348\\
0.935454 348\\
0.93633 348\\
0.937207 348\\
0.938086 348\\
0.938967 348\\
0.93985 348\\
0.940734 348\\
0.94162 348\\
0.942507 348\\
0.943396 347\\
0.944287 347\\
0.94518 347\\
0.946074 347\\
0.94697 347\\
0.947867 347\\
0.948767 347\\
0.949668 347\\
0.95057 346\\
0.951475 346\\
0.952381 346\\
0.953289 346\\
0.954198 346\\
0.95511 346\\
0.956023 346\\
0.956938 346\\
0.957854 346\\
0.958773 346\\
0.959693 346\\
0.960615 346\\
0.961538 346\\
0.962464 346\\
0.963391 346\\
0.96432 346\\
0.965251 346\\
0.966184 346\\
0.967118 346\\
0.968054 346\\
0.968992 346\\
0.969932 346\\
0.970874 346\\
0.971817 346\\
0.972763 346\\
0.97371 344\\
0.974659 344\\
0.97561 342\\
0.976562 343\\
0.977517 342\\
0.978474 342\\
0.979432 342\\
0.980392 342\\
0.981354 342\\
0.982318 342\\
0.983284 342\\
0.984252 341\\
0.985222 340\\
0.986193 340\\
0.987167 339\\
0.988142 335\\
0.98912 334\\
0.990099 334\\
0.99108 334\\
0.992063 334\\
0.993049 334\\
0.994036 334\\
0.995025 333\\
0.996016 331\\
0.997009 328\\
0.998004 326\\
0.999001 326\\
1 325\\
1.001 323\\
1.002 322\\
1.00301 319\\
1.00402 315\\
1.00503 309\\
1.00604 303\\
1.00705 298\\
1.00806 298\\
1.00908 298\\
1.0101 298\\
1.01112 298\\
1.01215 298\\
1.01317 298\\
1.0142 298\\
1.01523 298\\
1.01626 298\\
1.01729 298\\
1.01833 298\\
1.01937 298\\
1.02041 298\\
1.02145 298\\
1.02249 298\\
1.02354 298\\
1.02459 298\\
1.02564 298\\
1.02669 298\\
1.02775 298\\
1.02881 298\\
1.02987 298\\
1.03093 298\\
1.03199 298\\
1.03306 298\\
1.03413 298\\
1.0352 298\\
1.03627 303\\
1.03734 308\\
1.03842 312\\
1.0395 314\\
1.04058 318\\
1.04167 320\\
1.04275 322\\
1.04384 324\\
1.04493 324\\
1.04603 328\\
1.04712 330\\
1.04822 332\\
1.04932 332\\
1.05042 332\\
1.05152 332\\
1.05263 333\\
1.05374 333\\
1.05485 338\\
1.05597 340\\
1.05708 340\\
1.0582 341\\
1.05932 342\\
1.06045 342\\
1.06157 342\\
1.0627 342\\
1.06383 342\\
1.06496 342\\
1.0661 342\\
1.06724 342\\
1.06838 343\\
1.06952 344\\
1.07066 344\\
1.07181 346\\
1.07296 346\\
1.07411 346\\
1.07527 346\\
1.07643 346\\
1.07759 346\\
1.07875 346\\
1.07991 346\\
1.08108 346\\
1.08225 346\\
1.08342 346\\
1.0846 346\\
1.08578 346\\
1.08696 346\\
1.08814 346\\
1.08932 346\\
1.09051 346\\
1.0917 346\\
1.0929 346\\
1.09409 346\\
1.09529 346\\
1.09649 346\\
1.09769 346\\
1.0989 346\\
1.10011 346\\
1.10132 346\\
1.10254 346\\
1.10375 346\\
1.10497 346\\
1.10619 346\\
1.10742 347\\
1.10865 347\\
1.10988 347\\
1.11111 347\\
1.11235 347\\
1.11359 347\\
1.11483 347\\
1.11607 347\\
1.11732 347\\
1.11857 351\\
1.11982 353\\
1.12108 358\\
1.12233 363\\
1.1236 364\\
1.12486 369\\
1.12613 370\\
1.1274 372\\
1.12867 374\\
1.12994 374\\
1.13122 374\\
1.1325 377\\
1.13379 379\\
1.13507 381\\
1.13636 382\\
1.13766 382\\
1.13895 382\\
1.14025 382\\
1.14155 382\\
1.14286 383\\
1.14416 383\\
1.14548 384\\
1.14679 387\\
1.14811 389\\
1.14943 390\\
1.15075 390\\
1.15207 391\\
1.1534 391\\
1.15473 392\\
1.15607 392\\
1.15741 392\\
1.15875 392\\
1.16009 392\\
1.16144 392\\
1.16279 392\\
1.16414 392\\
1.1655 392\\
1.16686 394\\
1.16822 394\\
1.16959 394\\
1.17096 394\\
1.17233 394\\
1.17371 395\\
1.17509 396\\
1.17647 396\\
1.17786 396\\
1.17925 396\\
1.18064 396\\
1.18203 396\\
1.18343 396\\
1.18483 396\\
1.18624 396\\
1.18765 396\\
1.18906 396\\
1.19048 396\\
1.1919 396\\
1.19332 396\\
1.19474 396\\
1.19617 396\\
1.1976 396\\
1.19904 396\\
1.20048 396\\
1.20192 396\\
1.20337 396\\
1.20482 396\\
1.20627 396\\
1.20773 396\\
1.20919 396\\
1.21065 396\\
1.21212 396\\
1.21359 396\\
1.21507 396\\
1.21655 396\\
1.21803 396\\
1.21951 400\\
1.221 401\\
1.22249 404\\
1.22399 407\\
1.22549 411\\
1.22699 413\\
1.2285 417\\
1.23001 419\\
1.23153 420\\
1.23305 422\\
1.23457 422\\
1.23609 423\\
1.23762 423\\
1.23916 425\\
1.24069 427\\
1.24224 429\\
1.24378 429\\
1.24533 430\\
1.24688 431\\
1.24844 431\\
1.25 431\\
1.25156 431\\
1.25313 432\\
1.25471 432\\
1.25628 433\\
1.25786 435\\
1.25945 435\\
1.26103 436\\
1.26263 437\\
1.26422 438\\
1.26582 439\\
1.26743 440\\
1.26904 440\\
1.27065 440\\
1.27226 440\\
1.27389 441\\
1.27551 441\\
1.27714 441\\
1.27877 441\\
1.28041 441\\
1.28205 441\\
1.2837 441\\
1.28535 442\\
1.287 443\\
1.28866 443\\
1.29032 443\\
1.29199 443\\
1.29366 443\\
1.29534 443\\
1.29702 443\\
1.2987 443\\
1.30039 443\\
1.30208 443\\
1.30378 444\\
1.30548 445\\
1.30719 445\\
1.3089 445\\
1.31062 445\\
1.31234 445\\
1.31406 445\\
1.31579 445\\
1.31752 445\\
1.31926 445\\
1.321 445\\
1.32275 445\\
1.3245 445\\
1.32626 445\\
1.32802 445\\
1.32979 445\\
1.33156 445\\
1.33333 445\\
1.33511 445\\
1.3369 445\\
1.33869 445\\
1.34048 449\\
1.34228 451\\
1.34409 453\\
1.3459 456\\
1.34771 460\\
1.34953 463\\
1.35135 463\\
1.35318 466\\
1.35501 468\\
1.35685 468\\
1.3587 469\\
1.36054 470\\
1.3624 471\\
1.36426 471\\
1.36612 473\\
1.36799 474\\
1.36986 476\\
1.37174 477\\
1.37363 477\\
1.37552 478\\
1.37741 480\\
1.37931 480\\
1.38122 481\\
1.38313 481\\
1.38504 481\\
1.38696 482\\
1.38889 484\\
1.39082 484\\
1.39276 484\\
1.3947 485\\
1.39665 485\\
1.3986 486\\
1.40056 487\\
1.40252 488\\
1.40449 490\\
1.40647 490\\
1.40845 490\\
1.41044 490\\
1.41243 490\\
1.41443 490\\
1.41643 491\\
1.41844 491\\
1.42045 491\\
1.42248 491\\
1.4245 492\\
1.42653 493\\
1.42857 493\\
};

\addplot [
color=black,
solid,
line width=1.0pt
]
table[row sep=crcr]{
0.769231 488\\
0.769823 487\\
0.770416 487\\
0.77101 487\\
0.771605 487\\
0.772201 487\\
0.772798 487\\
0.773395 487\\
0.773994 486\\
0.774593 485\\
0.775194 484\\
0.775795 484\\
0.776398 484\\
0.777001 484\\
0.777605 484\\
0.77821 484\\
0.778816 484\\
0.779423 484\\
0.780031 484\\
0.78064 484\\
0.78125 484\\
0.781861 483\\
0.782473 481\\
0.783085 481\\
0.783699 482\\
0.784314 482\\
0.784929 480\\
0.785546 479\\
0.786164 478\\
0.786782 476\\
0.787402 476\\
0.788022 475\\
0.788644 474\\
0.789266 474\\
0.789889 474\\
0.790514 473\\
0.791139 473\\
0.791766 472\\
0.792393 472\\
0.793021 471\\
0.793651 471\\
0.794281 470\\
0.794913 470\\
0.795545 469\\
0.796178 468\\
0.796813 467\\
0.797448 465\\
0.798085 465\\
0.798722 462\\
0.799361 462\\
0.8 461\\
0.800641 458\\
0.801282 456\\
0.801925 455\\
0.802568 453\\
0.803213 452\\
0.803859 451\\
0.804505 451\\
0.805153 447\\
0.805802 447\\
0.806452 447\\
0.807103 447\\
0.807754 447\\
0.808407 447\\
0.809061 446\\
0.809717 446\\
0.810373 446\\
0.81103 446\\
0.811688 446\\
0.812348 446\\
0.813008 446\\
0.81367 446\\
0.814332 446\\
0.814996 446\\
0.815661 446\\
0.816327 446\\
0.816993 446\\
0.817661 446\\
0.818331 446\\
0.819001 446\\
0.819672 446\\
0.820345 446\\
0.821018 446\\
0.821693 446\\
0.822368 446\\
0.823045 445\\
0.823723 445\\
0.824402 444\\
0.825083 444\\
0.825764 444\\
0.826446 444\\
0.82713 444\\
0.827815 444\\
0.8285 444\\
0.829187 444\\
0.829876 444\\
0.830565 445\\
0.831255 446\\
0.831947 445\\
0.832639 445\\
0.833333 444\\
0.834028 444\\
0.834725 444\\
0.835422 444\\
0.83612 442\\
0.83682 442\\
0.837521 443\\
0.838223 443\\
0.838926 444\\
0.839631 444\\
0.840336 443\\
0.841043 442\\
0.841751 444\\
0.84246 444\\
0.84317 444\\
0.843882 442\\
0.844595 442\\
0.845309 442\\
0.846024 442\\
0.84674 442\\
0.847458 441\\
0.848176 441\\
0.848896 441\\
0.849618 441\\
0.85034 439\\
0.851064 438\\
0.851789 437\\
0.852515 437\\
0.853242 437\\
0.853971 437\\
0.854701 436\\
0.855432 435\\
0.856164 434\\
0.856898 434\\
0.857633 434\\
0.858369 434\\
0.859107 434\\
0.859845 434\\
0.860585 434\\
0.861326 434\\
0.862069 434\\
0.862813 434\\
0.863558 433\\
0.864304 432\\
0.865052 432\\
0.865801 431\\
0.866551 431\\
0.867303 429\\
0.868056 427\\
0.86881 424\\
0.869565 424\\
0.870322 423\\
0.87108 423\\
0.87184 423\\
0.8726 422\\
0.873362 421\\
0.874126 421\\
0.874891 418\\
0.875657 416\\
0.876424 414\\
0.877193 412\\
0.877963 409\\
0.878735 407\\
0.879507 405\\
0.880282 399\\
0.881057 401\\
0.881834 398\\
0.882613 397\\
0.883392 397\\
0.884173 397\\
0.884956 397\\
0.88574 397\\
0.886525 397\\
0.887311 397\\
0.888099 397\\
0.888889 397\\
0.88968 397\\
0.890472 397\\
0.891266 397\\
0.892061 396\\
0.892857 396\\
0.893655 396\\
0.894454 396\\
0.895255 396\\
0.896057 396\\
0.896861 396\\
0.897666 396\\
0.898473 396\\
0.899281 396\\
0.90009 396\\
0.900901 396\\
0.901713 396\\
0.902527 396\\
0.903342 396\\
0.904159 396\\
0.904977 396\\
0.905797 396\\
0.906618 396\\
0.907441 396\\
0.908265 396\\
0.909091 396\\
0.909918 396\\
0.910747 396\\
0.911577 396\\
0.912409 396\\
0.913242 396\\
0.914077 396\\
0.914913 396\\
0.915751 396\\
0.91659 396\\
0.917431 394\\
0.918274 394\\
0.919118 394\\
0.919963 394\\
0.92081 394\\
0.921659 394\\
0.922509 394\\
0.923361 392\\
0.924214 392\\
0.925069 392\\
0.925926 392\\
0.926784 392\\
0.927644 392\\
0.928505 392\\
0.929368 392\\
0.930233 392\\
0.931099 392\\
0.931966 392\\
0.932836 392\\
0.933707 392\\
0.934579 391\\
0.935454 390\\
0.93633 390\\
0.937207 390\\
0.938086 387\\
0.938967 387\\
0.93985 385\\
0.940734 384\\
0.94162 384\\
0.942507 384\\
0.943396 384\\
0.944287 384\\
0.94518 384\\
0.946074 383\\
0.94697 382\\
0.947867 381\\
0.948767 378\\
0.949668 375\\
0.95057 374\\
0.951475 374\\
0.952381 371\\
0.953289 370\\
0.954198 369\\
0.95511 364\\
0.956023 362\\
0.956938 358\\
0.957854 356\\
0.958773 351\\
0.959693 348\\
0.960615 348\\
0.961538 348\\
0.962464 348\\
0.963391 348\\
0.96432 348\\
0.965251 348\\
0.966184 348\\
0.967118 348\\
0.968054 348\\
0.968992 348\\
0.969932 348\\
0.970874 347\\
0.971817 347\\
0.972763 347\\
0.97371 347\\
0.974659 347\\
0.97561 346\\
0.976562 346\\
0.977517 346\\
0.978474 346\\
0.979432 346\\
0.980392 346\\
0.981354 346\\
0.982318 346\\
0.983284 346\\
0.984252 346\\
0.985222 346\\
0.986193 346\\
0.987167 346\\
0.988142 346\\
0.98912 346\\
0.990099 346\\
0.99108 346\\
0.992063 346\\
0.993049 346\\
0.994036 346\\
0.995025 346\\
0.996016 346\\
0.997009 346\\
0.998004 346\\
0.999001 344\\
1 344\\
1.001 344\\
1.002 342\\
1.00301 342\\
1.00402 342\\
1.00503 342\\
1.00604 342\\
1.00705 342\\
1.00806 342\\
1.00908 340\\
1.0101 340\\
1.01112 337\\
1.01215 335\\
1.01317 334\\
1.0142 334\\
1.01523 334\\
1.01626 334\\
1.01729 334\\
1.01833 334\\
1.01937 332\\
1.02041 331\\
1.02145 330\\
1.02249 329\\
1.02354 326\\
1.02459 326\\
1.02564 326\\
1.02669 325\\
1.02775 323\\
1.02881 322\\
1.02987 319\\
1.03093 317\\
1.03199 314\\
1.03306 313\\
1.03413 313\\
1.0352 310\\
1.03627 307\\
1.03734 306\\
1.03842 306\\
1.0395 304\\
1.04058 303\\
1.04167 302\\
1.04275 302\\
1.04384 302\\
1.04493 300\\
1.04603 298\\
1.04712 298\\
1.04822 298\\
1.04932 298\\
1.05042 298\\
1.05152 298\\
1.05263 298\\
1.05374 298\\
1.05485 298\\
1.05597 298\\
1.05708 298\\
1.0582 298\\
1.05932 298\\
1.06045 298\\
1.06157 298\\
1.0627 298\\
1.06383 298\\
1.06496 298\\
1.0661 298\\
1.06724 298\\
1.06838 298\\
1.06952 298\\
1.07066 298\\
1.07181 298\\
1.07296 298\\
1.07411 298\\
1.07527 298\\
1.07643 298\\
1.07759 298\\
1.07875 298\\
1.07991 298\\
1.08108 298\\
1.08225 298\\
1.08342 302\\
1.0846 308\\
1.08578 314\\
1.08696 320\\
1.08814 323\\
1.08932 324\\
1.09051 328\\
1.0917 332\\
1.0929 332\\
1.09409 332\\
1.09529 333\\
1.09649 336\\
1.09769 340\\
1.0989 341\\
1.10011 342\\
1.10132 342\\
1.10254 342\\
1.10375 342\\
1.10497 342\\
1.10619 344\\
1.10742 346\\
1.10865 346\\
1.10988 346\\
1.11111 346\\
1.11235 346\\
1.11359 346\\
1.11483 346\\
1.11607 346\\
1.11732 346\\
1.11857 346\\
1.11982 346\\
1.12108 346\\
1.12233 346\\
1.1236 346\\
1.12486 346\\
1.12613 346\\
1.1274 346\\
1.12867 346\\
1.12994 346\\
1.13122 346\\
1.1325 346\\
1.13379 347\\
1.13507 347\\
1.13636 347\\
1.13766 347\\
1.13895 347\\
1.14025 347\\
1.14155 348\\
1.14286 348\\
1.14416 348\\
1.14548 348\\
1.14679 352\\
1.14811 356\\
1.14943 358\\
1.15075 364\\
1.15207 368\\
1.1534 370\\
1.15473 372\\
1.15607 374\\
1.15741 374\\
1.15875 378\\
1.16009 379\\
1.16144 382\\
1.16279 382\\
1.16414 382\\
1.1655 382\\
1.16686 383\\
1.16822 383\\
1.16959 386\\
1.17096 388\\
1.17233 390\\
1.17371 390\\
1.17509 391\\
1.17647 392\\
1.17786 392\\
1.17925 392\\
1.18064 392\\
1.18203 392\\
1.18343 392\\
1.18483 392\\
1.18624 393\\
1.18765 394\\
1.18906 394\\
1.19048 394\\
1.1919 396\\
1.19332 396\\
1.19474 396\\
1.19617 396\\
1.1976 396\\
1.19904 396\\
1.20048 396\\
1.20192 396\\
1.20337 396\\
1.20482 396\\
1.20627 396\\
1.20773 396\\
1.20919 396\\
1.21065 396\\
1.21212 396\\
1.21359 396\\
1.21507 396\\
1.21655 396\\
1.21803 396\\
1.21951 396\\
1.221 396\\
1.22249 396\\
1.22399 396\\
1.22549 396\\
1.22699 396\\
1.2285 396\\
1.23001 396\\
1.23153 396\\
1.23305 396\\
1.23457 396\\
1.23609 397\\
1.23762 402\\
1.23916 407\\
1.24069 412\\
1.24224 418\\
1.24378 419\\
1.24533 424\\
1.24688 421\\
1.24844 423\\
1.25 422\\
1.25156 424\\
1.25313 425\\
1.25471 425\\
1.25628 430\\
1.25786 433\\
1.25945 435\\
1.26103 436\\
1.26263 435\\
1.26422 437\\
1.26582 437\\
1.26743 437\\
1.26904 438\\
1.27065 436\\
1.27226 437\\
1.27389 436\\
1.27551 440\\
1.27714 440\\
1.27877 441\\
1.28041 443\\
1.28205 443\\
1.2837 443\\
1.28535 442\\
1.287 441\\
1.28866 441\\
1.29032 443\\
1.29199 445\\
1.29366 446\\
1.29534 446\\
1.29702 446\\
1.2987 446\\
1.30039 446\\
1.30208 446\\
1.30378 445\\
1.30548 445\\
1.30719 446\\
1.3089 446\\
1.31062 446\\
1.31234 446\\
1.31406 446\\
1.31579 446\\
1.31752 446\\
1.31926 446\\
1.321 446\\
1.32275 446\\
1.3245 446\\
1.32626 446\\
1.32802 446\\
1.32979 446\\
1.33156 446\\
1.33333 449\\
1.33511 451\\
1.3369 453\\
1.33869 455\\
1.34048 461\\
1.34228 460\\
1.34409 466\\
1.3459 468\\
1.34771 469\\
1.34953 471\\
1.35135 471\\
1.35318 472\\
1.35501 473\\
1.35685 475\\
1.3587 476\\
1.36054 479\\
1.3624 479\\
1.36426 480\\
1.36612 481\\
1.36799 481\\
1.36986 481\\
1.37174 481\\
1.37363 482\\
1.37552 483\\
1.37741 484\\
1.37931 484\\
1.38122 486\\
1.38313 486\\
1.38504 486\\
1.38696 487\\
1.38889 489\\
1.39082 490\\
1.39276 490\\
1.3947 490\\
1.39665 490\\
1.3986 490\\
1.40056 491\\
1.40252 491\\
1.40449 491\\
1.40647 491\\
1.40845 491\\
1.41044 492\\
1.41243 493\\
1.41443 493\\
1.41643 493\\
1.41844 493\\
1.42045 493\\
1.42248 493\\
1.4245 493\\
1.42653 493\\
1.42857 493\\
};

\addplot [
color=black,
solid,
line width=1.0pt
]
table[row sep=crcr]{
0.769231 492\\
0.769823 492\\
0.770416 491\\
0.77101 491\\
0.771605 490\\
0.772201 489\\
0.772798 488\\
0.773395 487\\
0.773994 485\\
0.774593 486\\
0.775194 487\\
0.775795 488\\
0.776398 488\\
0.777001 488\\
0.777605 488\\
0.77821 487\\
0.778816 486\\
0.779423 483\\
0.780031 484\\
0.78064 487\\
0.78125 487\\
0.781861 486\\
0.782473 485\\
0.783085 482\\
0.783699 481\\
0.784314 481\\
0.784929 481\\
0.785546 480\\
0.786164 483\\
0.786782 481\\
0.787402 479\\
0.788022 476\\
0.788644 476\\
0.789266 477\\
0.789889 477\\
0.790514 477\\
0.791139 476\\
0.791766 475\\
0.792393 475\\
0.793021 476\\
0.793651 473\\
0.794281 474\\
0.794913 474\\
0.795545 471\\
0.796178 471\\
0.796813 470\\
0.797448 472\\
0.798085 470\\
0.798722 469\\
0.799361 472\\
0.8 467\\
0.800641 465\\
0.801282 466\\
0.801925 465\\
0.802568 462\\
0.803213 460\\
0.803859 462\\
0.804505 457\\
0.805153 458\\
0.805802 453\\
0.806452 452\\
0.807103 454\\
0.807754 452\\
0.808407 448\\
0.809061 449\\
0.809717 449\\
0.810373 448\\
0.81103 448\\
0.811688 448\\
0.812348 447\\
0.813008 447\\
0.81367 447\\
0.814332 446\\
0.814996 446\\
0.815661 446\\
0.816327 446\\
0.816993 446\\
0.817661 446\\
0.818331 446\\
0.819001 446\\
0.819672 446\\
0.820345 446\\
0.821018 446\\
0.821693 446\\
0.822368 446\\
0.823045 446\\
0.823723 446\\
0.824402 446\\
0.825083 446\\
0.825764 446\\
0.826446 446\\
0.82713 446\\
0.827815 446\\
0.8285 446\\
0.829187 446\\
0.829876 446\\
0.830565 446\\
0.831255 446\\
0.831947 446\\
0.832639 446\\
0.833333 446\\
0.834028 446\\
0.834725 446\\
0.835422 446\\
0.83612 445\\
0.83682 445\\
0.837521 445\\
0.838223 446\\
0.838926 446\\
0.839631 446\\
0.840336 446\\
0.841043 446\\
0.841751 446\\
0.84246 444\\
0.84317 444\\
0.843882 444\\
0.844595 444\\
0.845309 444\\
0.846024 444\\
0.84674 444\\
0.847458 444\\
0.848176 444\\
0.848896 443\\
0.849618 442\\
0.85034 442\\
0.851064 442\\
0.851789 441\\
0.852515 442\\
0.853242 439\\
0.853971 442\\
0.854701 441\\
0.855432 441\\
0.856164 442\\
0.856898 440\\
0.857633 442\\
0.858369 442\\
0.859107 442\\
0.859845 440\\
0.860585 440\\
0.861326 439\\
0.862069 440\\
0.862813 442\\
0.863558 440\\
0.864304 441\\
0.865052 439\\
0.865801 439\\
0.866551 437\\
0.867303 436\\
0.868056 437\\
0.86881 433\\
0.869565 431\\
0.870322 432\\
0.87108 433\\
0.87184 432\\
0.8726 432\\
0.873362 430\\
0.874126 429\\
0.874891 430\\
0.875657 428\\
0.876424 428\\
0.877193 428\\
0.877963 428\\
0.878735 426\\
0.879507 424\\
0.880282 420\\
0.881057 417\\
0.881834 418\\
0.882613 412\\
0.883392 414\\
0.884173 413\\
0.884956 410\\
0.88574 404\\
0.886525 403\\
0.887311 402\\
0.888099 402\\
0.888889 403\\
0.88968 404\\
0.890472 403\\
0.891266 404\\
0.892061 405\\
0.892857 406\\
0.893655 406\\
0.894454 402\\
0.895255 404\\
0.896057 402\\
0.896861 401\\
0.897666 396\\
0.898473 396\\
0.899281 396\\
0.90009 396\\
0.900901 396\\
0.901713 396\\
0.902527 396\\
0.903342 396\\
0.904159 396\\
0.904977 396\\
0.905797 396\\
0.906618 396\\
0.907441 396\\
0.908265 397\\
0.909091 397\\
0.909918 396\\
0.910747 396\\
0.911577 396\\
0.912409 396\\
0.913242 396\\
0.914077 396\\
0.914913 396\\
0.915751 396\\
0.91659 396\\
0.917431 396\\
0.918274 396\\
0.919118 396\\
0.919963 396\\
0.92081 396\\
0.921659 396\\
0.922509 396\\
0.923361 396\\
0.924214 394\\
0.925069 396\\
0.925926 396\\
0.926784 396\\
0.927644 394\\
0.928505 394\\
0.929368 394\\
0.930233 394\\
0.931099 394\\
0.931966 394\\
0.932836 394\\
0.933707 392\\
0.934579 392\\
0.935454 392\\
0.93633 392\\
0.937207 392\\
0.938086 391\\
0.938967 392\\
0.93985 392\\
0.940734 392\\
0.94162 392\\
0.942507 389\\
0.943396 388\\
0.944287 389\\
0.94518 388\\
0.946074 390\\
0.94697 389\\
0.947867 389\\
0.948767 387\\
0.949668 385\\
0.95057 382\\
0.951475 382\\
0.952381 381\\
0.953289 379\\
0.954198 379\\
0.95511 377\\
0.956023 377\\
0.956938 376\\
0.957854 376\\
0.958773 375\\
0.959693 372\\
0.960615 367\\
0.961538 363\\
0.962464 365\\
0.963391 356\\
0.96432 355\\
0.965251 352\\
0.966184 349\\
0.967118 350\\
0.968054 348\\
0.968992 348\\
0.969932 348\\
0.970874 348\\
0.971817 348\\
0.972763 348\\
0.97371 348\\
0.974659 348\\
0.97561 348\\
0.976562 348\\
0.977517 348\\
0.978474 347\\
0.979432 348\\
0.980392 347\\
0.981354 347\\
0.982318 347\\
0.983284 347\\
0.984252 346\\
0.985222 346\\
0.986193 346\\
0.987167 346\\
0.988142 346\\
0.98912 346\\
0.990099 346\\
0.99108 346\\
0.992063 346\\
0.993049 346\\
0.994036 346\\
0.995025 346\\
0.996016 346\\
0.997009 346\\
0.998004 346\\
0.999001 346\\
1 346\\
1.001 346\\
1.002 346\\
1.00301 346\\
1.00402 346\\
1.00503 346\\
1.00604 346\\
1.00705 346\\
1.00806 346\\
1.00908 346\\
1.0101 344\\
1.01112 344\\
1.01215 342\\
1.01317 344\\
1.0142 342\\
1.01523 342\\
1.01626 342\\
1.01729 342\\
1.01833 342\\
1.01937 342\\
1.02041 342\\
1.02145 342\\
1.02249 342\\
1.02354 342\\
1.02459 342\\
1.02564 341\\
1.02669 341\\
1.02775 340\\
1.02881 339\\
1.02987 339\\
1.03093 337\\
1.03199 336\\
1.03306 336\\
1.03413 335\\
1.0352 334\\
1.03627 334\\
1.03734 334\\
1.03842 334\\
1.0395 334\\
1.04058 334\\
1.04167 334\\
1.04275 334\\
1.04384 334\\
1.04493 334\\
1.04603 334\\
1.04712 334\\
1.04822 334\\
1.04932 334\\
1.05042 334\\
1.05152 334\\
1.05263 334\\
1.05374 334\\
1.05485 334\\
1.05597 333\\
1.05708 332\\
1.0582 331\\
1.05932 329\\
1.06045 326\\
1.06157 325\\
1.0627 324\\
1.06383 324\\
1.06496 324\\
1.0661 324\\
1.06724 324\\
1.06838 324\\
1.06952 324\\
1.07066 324\\
1.07181 324\\
1.07296 324\\
1.07411 324\\
1.07527 323\\
1.07643 321\\
1.07759 320\\
1.07875 314\\
1.07991 312\\
1.08108 308\\
1.08225 306\\
1.08342 305\\
1.0846 306\\
1.08578 303\\
1.08696 302\\
1.08814 301\\
1.08932 298\\
1.09051 298\\
1.0917 298\\
1.0929 298\\
1.09409 298\\
1.09529 298\\
1.09649 302\\
1.09769 306\\
1.0989 312\\
1.10011 318\\
1.10132 320\\
1.10254 324\\
1.10375 328\\
1.10497 331\\
1.10619 332\\
1.10742 333\\
1.10865 334\\
1.10988 340\\
1.11111 341\\
1.11235 342\\
1.11359 342\\
1.11483 342\\
1.11607 342\\
1.11732 342\\
1.11857 344\\
1.11982 346\\
1.12108 346\\
1.12233 346\\
1.1236 346\\
1.12486 346\\
1.12613 346\\
1.1274 346\\
1.12867 346\\
1.12994 346\\
1.13122 346\\
1.1325 346\\
1.13379 346\\
1.13507 346\\
1.13636 346\\
1.13766 346\\
1.13895 346\\
1.14025 346\\
1.14155 346\\
1.14286 346\\
1.14416 347\\
1.14548 347\\
1.14679 347\\
1.14811 347\\
1.14943 347\\
1.15075 348\\
1.15207 348\\
1.1534 348\\
1.15473 348\\
1.15607 348\\
1.15741 348\\
1.15875 352\\
1.16009 356\\
1.16144 363\\
1.16279 366\\
1.16414 370\\
1.1655 373\\
1.16686 374\\
1.16822 374\\
1.16959 378\\
1.17096 379\\
1.17233 381\\
1.17371 382\\
1.17509 382\\
1.17647 382\\
1.17786 383\\
1.17925 383\\
1.18064 384\\
1.18203 385\\
1.18343 387\\
1.18483 389\\
1.18624 390\\
1.18765 390\\
1.18906 391\\
1.19048 392\\
1.1919 392\\
1.19332 392\\
1.19474 392\\
1.19617 392\\
1.1976 393\\
1.19904 393\\
1.20048 394\\
1.20192 393\\
1.20337 396\\
1.20482 394\\
1.20627 396\\
1.20773 396\\
1.20919 396\\
1.21065 396\\
1.21212 396\\
1.21359 396\\
1.21507 396\\
1.21655 396\\
1.21803 396\\
1.21951 396\\
1.221 396\\
1.22249 396\\
1.22399 396\\
1.22549 396\\
1.22699 396\\
1.2285 396\\
1.23001 396\\
1.23153 396\\
1.23305 396\\
1.23457 397\\
1.23609 398\\
1.23762 405\\
1.23916 405\\
1.24069 410\\
1.24224 407\\
1.24378 408\\
1.24533 405\\
1.24688 411\\
1.24844 409\\
1.25 405\\
1.25156 407\\
1.25313 411\\
1.25471 419\\
1.25628 411\\
1.25786 422\\
1.25945 422\\
1.26103 420\\
1.26263 424\\
1.26422 427\\
1.26582 429\\
1.26743 429\\
1.26904 428\\
1.27065 428\\
1.27226 430\\
1.27389 431\\
1.27551 431\\
1.27714 432\\
1.27877 434\\
1.28041 435\\
1.28205 435\\
1.2837 436\\
1.28535 436\\
1.287 438\\
1.28866 441\\
1.29032 440\\
1.29199 441\\
1.29366 441\\
1.29534 441\\
1.29702 441\\
1.2987 443\\
1.30039 443\\
1.30208 442\\
1.30378 441\\
1.30548 442\\
1.30719 442\\
1.3089 442\\
1.31062 442\\
1.31234 443\\
1.31406 444\\
1.31579 444\\
1.31752 444\\
1.31926 444\\
1.321 444\\
1.32275 444\\
1.3245 445\\
1.32626 444\\
1.32802 444\\
1.32979 446\\
1.33156 446\\
1.33333 446\\
1.33511 446\\
1.3369 446\\
1.33869 446\\
1.34048 446\\
1.34228 446\\
1.34409 446\\
1.3459 446\\
1.34771 446\\
1.34953 446\\
1.35135 450\\
1.35318 451\\
1.35501 453\\
1.35685 454\\
1.3587 458\\
1.36054 457\\
1.3624 462\\
1.36426 465\\
1.36612 467\\
1.36799 467\\
1.36986 472\\
1.37174 472\\
1.37363 472\\
1.37552 473\\
1.37741 476\\
1.37931 476\\
1.38122 479\\
1.38313 480\\
1.38504 478\\
1.38696 479\\
1.38889 482\\
1.39082 482\\
1.39276 483\\
1.3947 483\\
1.39665 483\\
1.3986 482\\
1.40056 484\\
1.40252 484\\
1.40449 484\\
1.40647 484\\
1.40845 486\\
1.41044 486\\
1.41243 486\\
1.41443 487\\
1.41643 488\\
1.41844 490\\
1.42045 490\\
1.42248 491\\
1.4245 490\\
1.42653 490\\
1.42857 491\\
};

\addplot [
color=black,
solid,
line width=1.0pt
]
table[row sep=crcr]{
0.769231 485\\
0.769823 484\\
0.770416 483\\
0.77101 483\\
0.771605 481\\
0.772201 482\\
0.772798 480\\
0.773395 480\\
0.773994 480\\
0.774593 479\\
0.775194 478\\
0.775795 477\\
0.776398 478\\
0.777001 476\\
0.777605 476\\
0.77821 476\\
0.778816 476\\
0.779423 474\\
0.780031 475\\
0.78064 474\\
0.78125 471\\
0.781861 471\\
0.782473 470\\
0.783085 471\\
0.783699 472\\
0.784314 471\\
0.784929 473\\
0.785546 472\\
0.786164 466\\
0.786782 465\\
0.787402 462\\
0.788022 466\\
0.788644 463\\
0.789266 464\\
0.789889 463\\
0.790514 465\\
0.791139 462\\
0.791766 459\\
0.792393 456\\
0.793021 463\\
0.793651 461\\
0.794281 461\\
0.794913 462\\
0.795545 457\\
0.796178 453\\
0.796813 457\\
0.797448 457\\
0.798085 456\\
0.798722 453\\
0.799361 455\\
0.8 455\\
0.800641 452\\
0.801282 450\\
0.801925 452\\
0.802568 452\\
0.803213 451\\
0.803859 451\\
0.804505 450\\
0.805153 449\\
0.805802 447\\
0.806452 446\\
0.807103 447\\
0.807754 447\\
0.808407 446\\
0.809061 446\\
0.809717 446\\
0.810373 446\\
0.81103 446\\
0.811688 446\\
0.812348 446\\
0.813008 446\\
0.81367 446\\
0.814332 445\\
0.814996 444\\
0.815661 445\\
0.816327 445\\
0.816993 446\\
0.817661 445\\
0.818331 445\\
0.819001 445\\
0.819672 445\\
0.820345 445\\
0.821018 444\\
0.821693 444\\
0.822368 444\\
0.823045 444\\
0.823723 444\\
0.824402 444\\
0.825083 444\\
0.825764 444\\
0.826446 444\\
0.82713 444\\
0.827815 444\\
0.8285 444\\
0.829187 444\\
0.829876 444\\
0.830565 444\\
0.831255 444\\
0.831947 444\\
0.832639 442\\
0.833333 442\\
0.834028 443\\
0.834725 442\\
0.835422 442\\
0.83612 442\\
0.83682 442\\
0.837521 442\\
0.838223 442\\
0.838926 442\\
0.839631 442\\
0.840336 442\\
0.841043 441\\
0.841751 441\\
0.84246 441\\
0.84317 441\\
0.843882 441\\
0.844595 441\\
0.845309 440\\
0.846024 439\\
0.84674 438\\
0.847458 439\\
0.848176 438\\
0.848896 439\\
0.849618 437\\
0.85034 437\\
0.851064 434\\
0.851789 437\\
0.852515 437\\
0.853242 437\\
0.853971 437\\
0.854701 436\\
0.855432 439\\
0.856164 437\\
0.856898 440\\
0.857633 437\\
0.858369 436\\
0.859107 435\\
0.859845 434\\
0.860585 434\\
0.861326 434\\
0.862069 433\\
0.862813 429\\
0.863558 428\\
0.864304 427\\
0.865052 426\\
0.865801 425\\
0.866551 424\\
0.867303 424\\
0.868056 425\\
0.86881 424\\
0.869565 425\\
0.870322 423\\
0.87108 422\\
0.87184 420\\
0.8726 418\\
0.873362 416\\
0.874126 414\\
0.874891 409\\
0.875657 408\\
0.876424 405\\
0.877193 402\\
0.877963 402\\
0.878735 400\\
0.879507 398\\
0.880282 400\\
0.881057 397\\
0.881834 397\\
0.882613 397\\
0.883392 397\\
0.884173 397\\
0.884956 397\\
0.88574 397\\
0.886525 397\\
0.887311 397\\
0.888099 397\\
0.888889 397\\
0.88968 396\\
0.890472 397\\
0.891266 396\\
0.892061 396\\
0.892857 396\\
0.893655 396\\
0.894454 396\\
0.895255 396\\
0.896057 396\\
0.896861 396\\
0.897666 396\\
0.898473 396\\
0.899281 396\\
0.90009 396\\
0.900901 396\\
0.901713 396\\
0.902527 396\\
0.903342 396\\
0.904159 396\\
0.904977 396\\
0.905797 396\\
0.906618 396\\
0.907441 396\\
0.908265 396\\
0.909091 396\\
0.909918 396\\
0.910747 396\\
0.911577 396\\
0.912409 394\\
0.913242 394\\
0.914077 394\\
0.914913 394\\
0.915751 394\\
0.91659 394\\
0.917431 394\\
0.918274 394\\
0.919118 392\\
0.919963 392\\
0.92081 392\\
0.921659 392\\
0.922509 392\\
0.923361 392\\
0.924214 392\\
0.925069 392\\
0.925926 392\\
0.926784 392\\
0.927644 392\\
0.928505 392\\
0.929368 392\\
0.930233 392\\
0.931099 392\\
0.931966 391\\
0.932836 391\\
0.933707 390\\
0.934579 390\\
0.935454 390\\
0.93633 387\\
0.937207 386\\
0.938086 384\\
0.938967 384\\
0.93985 384\\
0.940734 384\\
0.94162 384\\
0.942507 384\\
0.943396 384\\
0.944287 384\\
0.94518 383\\
0.946074 383\\
0.94697 382\\
0.947867 381\\
0.948767 379\\
0.949668 378\\
0.95057 378\\
0.951475 374\\
0.952381 374\\
0.953289 373\\
0.954198 372\\
0.95511 371\\
0.956023 370\\
0.956938 370\\
0.957854 368\\
0.958773 366\\
0.959693 362\\
0.960615 363\\
0.961538 360\\
0.962464 355\\
0.963391 354\\
0.96432 353\\
0.965251 349\\
0.966184 348\\
0.967118 348\\
0.968054 348\\
0.968992 348\\
0.969932 347\\
0.970874 347\\
0.971817 347\\
0.972763 348\\
0.97371 348\\
0.974659 348\\
0.97561 348\\
0.976562 348\\
0.977517 347\\
0.978474 346\\
0.979432 346\\
0.980392 346\\
0.981354 346\\
0.982318 346\\
0.983284 346\\
0.984252 347\\
0.985222 347\\
0.986193 346\\
0.987167 346\\
0.988142 346\\
0.98912 346\\
0.990099 346\\
0.99108 346\\
0.992063 346\\
0.993049 346\\
0.994036 346\\
0.995025 346\\
0.996016 346\\
0.997009 346\\
0.998004 346\\
0.999001 346\\
1 346\\
1.001 346\\
1.002 346\\
1.00301 346\\
1.00402 346\\
1.00503 344\\
1.00604 344\\
1.00705 344\\
1.00806 343\\
1.00908 342\\
1.0101 342\\
1.01112 342\\
1.01215 342\\
1.01317 342\\
1.0142 342\\
1.01523 342\\
1.01626 342\\
1.01729 342\\
1.01833 342\\
1.01937 342\\
1.02041 341\\
1.02145 340\\
1.02249 340\\
1.02354 340\\
1.02459 338\\
1.02564 337\\
1.02669 336\\
1.02775 334\\
1.02881 334\\
1.02987 334\\
1.03093 334\\
1.03199 334\\
1.03306 334\\
1.03413 334\\
1.0352 335\\
1.03627 334\\
1.03734 334\\
1.03842 335\\
1.0395 337\\
1.04058 336\\
1.04167 339\\
1.04275 340\\
1.04384 338\\
1.04493 341\\
1.04603 339\\
1.04712 338\\
1.04822 338\\
1.04932 335\\
1.05042 334\\
1.05152 333\\
1.05263 333\\
1.05374 334\\
1.05485 330\\
1.05597 331\\
1.05708 331\\
1.0582 329\\
1.05932 330\\
1.06045 332\\
1.06157 329\\
1.0627 325\\
1.06383 325\\
1.06496 326\\
1.0661 326\\
1.06724 328\\
1.06838 324\\
1.06952 320\\
1.07066 320\\
1.07181 319\\
1.07296 315\\
1.07411 314\\
1.07527 316\\
1.07643 314\\
1.07759 317\\
1.07875 313\\
1.07991 322\\
1.08108 319\\
1.08225 312\\
1.08342 306\\
1.0846 302\\
1.08578 298\\
1.08696 298\\
1.08814 298\\
1.08932 298\\
1.09051 298\\
1.0917 298\\
1.0929 298\\
1.09409 298\\
1.09529 302\\
1.09649 304\\
1.09769 306\\
1.0989 306\\
1.10011 310\\
1.10132 317\\
1.10254 315\\
1.10375 324\\
1.10497 324\\
1.10619 324\\
1.10742 328\\
1.10865 330\\
1.10988 332\\
1.11111 332\\
1.11235 332\\
1.11359 332\\
1.11483 333\\
1.11607 334\\
1.11732 333\\
1.11857 335\\
1.11982 340\\
1.12108 340\\
1.12233 340\\
1.1236 341\\
1.12486 342\\
1.12613 342\\
1.1274 344\\
1.12867 343\\
1.12994 344\\
1.13122 342\\
1.1325 344\\
1.13379 345\\
1.13507 346\\
1.13636 346\\
1.13766 346\\
1.13895 346\\
1.14025 346\\
1.14155 346\\
1.14286 346\\
1.14416 346\\
1.14548 346\\
1.14679 346\\
1.14811 346\\
1.14943 346\\
1.15075 346\\
1.15207 346\\
1.1534 346\\
1.15473 346\\
1.15607 346\\
1.15741 346\\
1.15875 347\\
1.16009 347\\
1.16144 347\\
1.16279 347\\
1.16414 347\\
1.1655 347\\
1.16686 347\\
1.16822 347\\
1.16959 347\\
1.17096 347\\
1.17233 348\\
1.17371 348\\
1.17509 348\\
1.17647 349\\
1.17786 353\\
1.17925 356\\
1.18064 359\\
1.18203 365\\
1.18343 370\\
1.18483 371\\
1.18624 372\\
1.18765 374\\
1.18906 374\\
1.19048 374\\
1.1919 380\\
1.19332 380\\
1.19474 383\\
1.19617 382\\
1.1976 383\\
1.19904 384\\
1.20048 384\\
1.20192 385\\
1.20337 384\\
1.20482 386\\
1.20627 388\\
1.20773 391\\
1.20919 390\\
1.21065 391\\
1.21212 391\\
1.21359 392\\
1.21507 392\\
1.21655 392\\
1.21803 392\\
1.21951 394\\
1.221 394\\
1.22249 394\\
1.22399 393\\
1.22549 396\\
1.22699 396\\
1.2285 396\\
1.23001 395\\
1.23153 395\\
1.23305 396\\
1.23457 396\\
1.23609 396\\
1.23762 396\\
1.23916 396\\
1.24069 396\\
1.24224 396\\
1.24378 396\\
1.24533 396\\
1.24688 396\\
1.24844 396\\
1.25 396\\
1.25156 396\\
1.25313 396\\
1.25471 396\\
1.25628 396\\
1.25786 396\\
1.25945 396\\
1.26103 396\\
1.26263 396\\
1.26422 396\\
1.26582 396\\
1.26743 396\\
1.26904 396\\
1.27065 396\\
1.27226 396\\
1.27389 400\\
1.27551 401\\
1.27714 404\\
1.27877 407\\
1.28041 412\\
1.28205 415\\
1.2837 417\\
1.28535 419\\
1.287 421\\
1.28866 422\\
1.29032 422\\
1.29199 424\\
1.29366 426\\
1.29534 428\\
1.29702 429\\
1.2987 431\\
1.30039 432\\
1.30208 432\\
1.30378 433\\
1.30548 434\\
1.30719 434\\
1.3089 437\\
1.31062 436\\
1.31234 436\\
1.31406 436\\
1.31579 437\\
1.31752 439\\
1.31926 440\\
1.321 440\\
1.32275 440\\
1.3245 440\\
1.32626 442\\
1.32802 442\\
1.32979 442\\
1.33156 442\\
1.33333 442\\
1.33511 441\\
1.3369 442\\
1.33869 443\\
1.34048 443\\
1.34228 445\\
1.34409 445\\
1.3459 445\\
1.34771 445\\
1.34953 445\\
1.35135 446\\
1.35318 445\\
1.35501 445\\
1.35685 445\\
1.3587 446\\
1.36054 446\\
1.3624 446\\
1.36426 445\\
1.36612 445\\
1.36799 445\\
1.36986 445\\
1.37174 445\\
1.37363 445\\
1.37552 446\\
1.37741 445\\
1.37931 448\\
1.38122 449\\
1.38313 451\\
1.38504 450\\
1.38696 456\\
1.38889 454\\
1.39082 454\\
1.39276 459\\
1.3947 456\\
1.39665 463\\
1.3986 465\\
1.40056 469\\
1.40252 466\\
1.40449 470\\
1.40647 474\\
1.40845 478\\
1.41044 481\\
1.41243 480\\
1.41443 479\\
1.41643 481\\
1.41844 481\\
1.42045 481\\
1.42248 481\\
1.4245 483\\
1.42653 483\\
1.42857 483\\
};

\addplot [
color=black,
only marks,
mark=x,
mark size = 2.5,
mark options={solid,blue}
]
table[row sep=crcr]{
0.769231 447\\
0.775194 446\\
0.78125 446\\
0.787402 446\\
0.793651 444\\
0.8 444\\
0.806452 443\\
0.813008 442\\
0.819672 439\\
0.826446 435\\
0.833333 433\\
0.840336 425\\
0.847458 412\\
0.854701 397\\
0.862069 397\\
0.869565 396\\
0.877193 396\\
0.884956 396\\
0.892857 394\\
0.900901 392\\
0.909091 388\\
0.917431 384\\
0.925926 373\\
0.934579 348\\
0.943396 347\\
0.952381 346\\
0.961538 346\\
0.970874 346\\
0.980392 342\\
0.990099 334\\
1 325\\
1.0101 298\\
1.02041 298\\
1.03093 298\\
1.04167 320\\
1.05263 333\\
1.06383 342\\
1.07527 346\\
1.08696 346\\
1.0989 346\\
1.11111 347\\
1.1236 364\\
1.13636 382\\
1.14943 390\\
1.16279 392\\
1.17647 396\\
1.19048 396\\
1.20482 396\\
1.21951 400\\
1.23457 422\\
1.25 431\\
1.26582 439\\
1.28205 441\\
1.2987 443\\
1.31579 445\\
1.33333 445\\
1.35135 463\\
1.36986 476\\
1.38889 484\\
1.40845 490\\
1.42857 493\\
};

\addplot [
color=black,
only marks,
mark=o,
mark size = 2.5,
mark options={solid,green}
]
table[row sep=crcr]{
0.769231 488\\
0.775194 484\\
0.78125 484\\
0.787402 476\\
0.793651 471\\
0.8 461\\
0.806452 447\\
0.813008 446\\
0.819672 446\\
0.826446 444\\
0.833333 444\\
0.840336 443\\
0.847458 441\\
0.854701 436\\
0.862069 434\\
0.869565 424\\
0.877193 412\\
0.884956 397\\
0.892857 396\\
0.900901 396\\
0.909091 396\\
0.917431 394\\
0.925926 392\\
0.934579 391\\
0.943396 384\\
0.952381 371\\
0.961538 348\\
0.970874 347\\
0.980392 346\\
0.990099 346\\
1 344\\
1.0101 340\\
1.02041 331\\
1.03093 317\\
1.04167 302\\
1.05263 298\\
1.06383 298\\
1.07527 298\\
1.08696 320\\
1.0989 341\\
1.11111 346\\
1.1236 346\\
1.13636 347\\
1.14943 358\\
1.16279 382\\
1.17647 392\\
1.19048 394\\
1.20482 396\\
1.21951 396\\
1.23457 396\\
1.25 422\\
1.26582 437\\
1.28205 443\\
1.2987 446\\
1.31579 446\\
1.33333 449\\
1.35135 471\\
1.36986 481\\
1.38889 489\\
1.40845 491\\
1.42857 493\\
};

\addplot [
color=black,
only marks,
mark=+,
mark size = 2.5,
mark options={solid,red}
]
table[row sep=crcr]{
0.769231 492\\
0.775194 487\\
0.78125 487\\
0.787402 479\\
0.793651 473\\
0.8 467\\
0.806452 452\\
0.813008 447\\
0.819672 446\\
0.826446 446\\
0.833333 446\\
0.840336 446\\
0.847458 444\\
0.854701 441\\
0.862069 440\\
0.869565 431\\
0.877193 428\\
0.884956 410\\
0.892857 406\\
0.900901 396\\
0.909091 397\\
0.917431 396\\
0.925926 396\\
0.934579 392\\
0.943396 388\\
0.952381 381\\
0.961538 363\\
0.970874 348\\
0.980392 347\\
0.990099 346\\
1 346\\
1.0101 344\\
1.02041 342\\
1.03093 337\\
1.04167 334\\
1.05263 334\\
1.06383 324\\
1.07527 323\\
1.08696 302\\
1.0989 312\\
1.11111 341\\
1.1236 346\\
1.13636 346\\
1.14943 347\\
1.16279 366\\
1.17647 382\\
1.19048 392\\
1.20482 394\\
1.21951 396\\
1.23457 397\\
1.25 405\\
1.26582 429\\
1.28205 435\\
1.2987 443\\
1.31579 444\\
1.33333 446\\
1.35135 450\\
1.36986 472\\
1.38889 482\\
1.40845 486\\
1.42857 491\\
};

\addplot [
color=black,
only marks,
mark=star,
mark size = 2.5,
mark options={solid,orange}
]
table[row sep=crcr]{
0.769231 485\\
0.775194 478\\
0.78125 471\\
0.787402 462\\
0.793651 461\\
0.8 455\\
0.806452 446\\
0.813008 446\\
0.819672 445\\
0.826446 444\\
0.833333 442\\
0.840336 442\\
0.847458 439\\
0.854701 436\\
0.862069 433\\
0.869565 425\\
0.877193 402\\
0.884956 397\\
0.892857 396\\
0.900901 396\\
0.909091 396\\
0.917431 394\\
0.925926 392\\
0.934579 390\\
0.943396 384\\
0.952381 374\\
0.961538 360\\
0.970874 347\\
0.980392 346\\
0.990099 346\\
1 346\\
1.0101 342\\
1.02041 341\\
1.03093 334\\
1.04167 339\\
1.05263 333\\
1.06383 325\\
1.07527 316\\
1.08696 298\\
1.0989 306\\
1.11111 332\\
1.1236 341\\
1.13636 346\\
1.14943 346\\
1.16279 347\\
1.17647 349\\
1.19048 374\\
1.20482 386\\
1.21951 394\\
1.23457 396\\
1.25 396\\
1.26582 396\\
1.28205 415\\
1.2987 431\\
1.31579 437\\
1.33333 442\\
1.35135 446\\
1.36986 445\\
1.38889 454\\
1.40845 478\\
1.42857 483\\
};

\end{axis}
\end{tikzpicture}
\caption{Total fixed-stress iterations for varying polynomial degree $s$ for cGP(1) in 
time.}
\label{fig:3:cG1:p}
\end{figure}
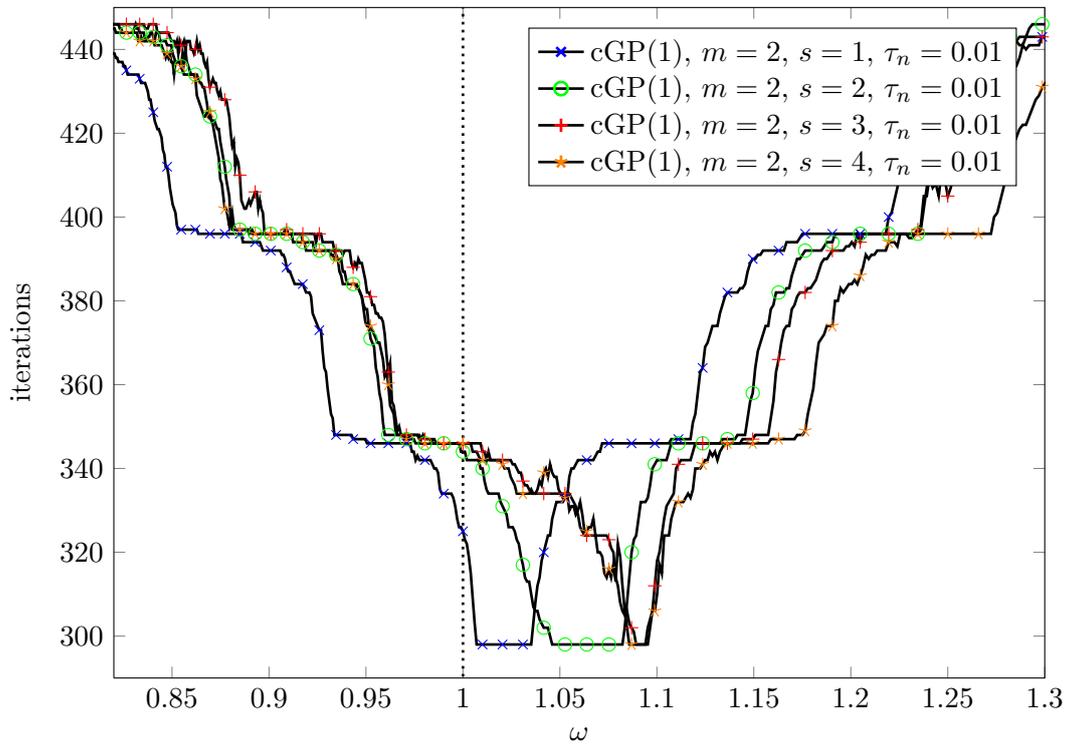

%% file: tau-independent.tex
\begin{figure}[p]
\centering
%
\begin{tikzpicture}

\begin{axis}[%
width=4.82222222222222in,
height=3.5in,
scale only axis,
/pgf/number format/.cd, 1000 sep={},
xlabel={$\omega$},
ylabel={iterations},
xmin=0.84,
xmax=1.25,
ymin=100,
ymax=1100,
yminorticks=true,
legend style={legend pos=north east},
legend style={draw=black,fill=white,legend cell align=left},
legend entries = {
  {dG(1), $m=2$, $s=2$, $\tau_n=0.005$},
  {dG(1), $m=2$, $s=2$, $\tau_n=0.010$},
  {\textcolor{black!70}{cG(1), $m=2$, $s=2$, $\tau_n=0.010$}},
  {dG(1), $m=2$, $s=2$, $\tau_n=0.020$}
}
]

\addplot [
color=black,
solid,
line width=1.0pt,
mark=+,
mark size = 2.5,
mark options={solid,red}
]
table[row sep=crcr]{
0.0 0\\
0.1 0\\
};

\addplot [
color=black,
solid,
line width=1.0pt,
mark=o,
mark size = 2.5,
mark options={solid,green}
]
table[row sep=crcr]{
0.0 0\\
0.1 0\\
};

\addplot [
color=black!70,
solid,
line width=1.0pt,
mark=star,
mark size = 2.5,
mark options={solid,orange}
]
table[row sep=crcr]{
0.0 0\\
0.1 0\\
};

\addplot [
color=black,
solid,
line width=1.0pt,
mark=x,
mark size = 2.5,
mark options={solid,blue}
]
table[row sep=crcr]{
0.0 0\\
0.1 0\\
};

\addplot [
color=black,
dotted,
line width=1.0pt
]
table[row sep=crcr]{
1.0 0\\
1.0 2000\\
};

\addplot [
color=black,
solid,
line width=1.0pt
]
table[row sep=crcr]{
0.840336 237\\
0.847458 232\\
0.854701 226\\
0.862069 222\\
0.869565 219\\
0.877193 215\\
0.884956 215\\
0.892857 207\\
0.900901 200\\
0.909091 201\\
0.917431 201\\
0.925926 204\\
0.934579 205\\
0.943396 196\\
0.952381 194\\
0.961538 188\\
0.970874 174\\
0.980392 174\\
0.990099 175\\
1 174\\
1.0101 172\\
1.02041 173\\
1.03093 168\\
1.04167 169\\
1.05263 168\\
1.06383 168\\
1.07527 165\\
1.08696 165\\
1.0989 166\\
1.11111 174\\
1.1236 173\\
1.13636 176\\
1.14943 174\\
1.16279 178\\
1.17647 192\\
1.19048 201\\
1.20482 201\\
1.21951 201\\
1.23457 201\\
1.25 202\\
};

\addplot [
color=black,
solid,
line width=1.0pt
]
table[row sep=crcr]{
0.840336 445\\
0.847458 446\\
0.854701 444\\
0.862069 443\\
0.869565 440\\
0.877193 435\\
0.884956 430\\
0.892857 418\\
0.900901 413\\
0.909091 409\\
0.917431 402\\
0.925926 398\\
0.934579 395\\
0.943396 392\\
0.952381 388\\
0.961538 383\\
0.970874 373\\
0.980392 361\\
0.990099 351\\
1 360\\
1.0101 356\\
1.02041 353\\
1.03093 349\\
1.04167 348\\
1.05263 348\\
1.06383 347\\
1.07527 345\\
1.08696 345\\
1.0989 348\\
1.11111 349\\
1.1236 350\\
1.13636 349\\
1.14943 351\\
1.16279 373\\
1.17647 384\\
1.19048 394\\
1.20482 412\\
1.21951 406\\
1.23457 405\\
1.25 415\\
};

\addplot [
color=black,
solid,
line width=1.0pt
]
table[row sep=crcr]{
0.840336 897\\
0.847458 907\\
0.854701 892\\
0.862069 905\\
0.869565 890\\
0.877193 902\\
0.884956 896\\
0.892857 910\\
0.900901 876\\
0.909091 846\\
0.917431 810\\
0.925926 808\\
0.934579 812\\
0.943396 808\\
0.952381 798\\
0.961538 805\\
0.970874 818\\
0.980392 804\\
0.990099 804\\
1 813\\
1.0101 783\\
1.02041 772\\
1.03093 775\\
1.04167 763\\
1.05263 755\\
1.06383 718\\
1.07527 722\\
1.08696 708\\
1.0989 727\\
1.11111 762\\
1.1236 764\\
1.13636 759\\
1.14943 774\\
1.16279 783\\
1.17647 794\\
1.19048 797\\
1.20482 822\\
1.21951 866\\
1.23457 881\\
1.25 876\\
};

\addplot [
color=black!70,
solid,
line width=1.0pt
]
table[row sep=crcr]{
0.769231 488\\
0.775194 484\\
0.78125 484\\
0.787402 476\\
0.793651 471\\
0.8 461\\
0.806452 447\\
0.813008 446\\
0.819672 446\\
0.826446 444\\
0.833333 444\\
0.840336 443\\
0.847458 441\\
0.854701 436\\
0.862069 434\\
0.869565 424\\
0.877193 412\\
0.884956 397\\
0.892857 396\\
0.900901 396\\
0.909091 396\\
0.917431 394\\
0.925926 392\\
0.934579 391\\
0.943396 384\\
0.952381 371\\
0.961538 348\\
0.970874 347\\
0.980392 346\\
0.990099 346\\
1 344\\
1.0101 340\\
1.02041 331\\
1.03093 317\\
1.04167 302\\
1.05263 298\\
1.06383 298\\
1.07527 298\\
1.08696 320\\
1.0989 341\\
1.11111 346\\
1.1236 346\\
1.13636 347\\
1.14943 358\\
1.16279 382\\
1.17647 392\\
1.19048 394\\
1.20482 396\\
1.21951 396\\
1.23457 396\\
1.25 422\\
1.26582 437\\
1.28205 443\\
1.2987 446\\
1.31579 446\\
1.33333 449\\
1.35135 471\\
1.36986 481\\
1.38889 489\\
1.40845 491\\
1.42857 493\\
};

\addplot [
color=black,
solid,
line width=1.0pt,
mark=x,
mark size = 2.5,
mark options={solid,blue}
]
table[row sep=crcr]{
0.840336 237\\
0.847458 232\\
0.854701 226\\
0.862069 222\\
0.869565 219\\
0.877193 215\\
0.884956 215\\
0.892857 207\\
0.900901 200\\
0.909091 201\\
0.917431 201\\
0.925926 204\\
0.934579 205\\
0.943396 196\\
0.952381 194\\
0.961538 188\\
0.970874 174\\
0.980392 174\\
0.990099 175\\
1 174\\
1.0101 172\\
1.02041 173\\
1.03093 168\\
1.04167 169\\
1.05263 168\\
1.06383 168\\
1.07527 165\\
1.08696 165\\
1.0989 166\\
1.11111 174\\
1.1236 173\\
1.13636 176\\
1.14943 174\\
1.16279 178\\
1.17647 192\\
1.19048 201\\
1.20482 201\\
1.21951 201\\
1.23457 201\\
1.25 202\\
};

\addplot [
color=black,
solid,
line width=1.0pt,
mark=o,
mark size = 2.5,
mark options={solid,green}
]
table[row sep=crcr]{
0.840336 445\\
0.847458 446\\
0.854701 444\\
0.862069 443\\
0.869565 440\\
0.877193 435\\
0.884956 430\\
0.892857 418\\
0.900901 413\\
0.909091 409\\
0.917431 402\\
0.925926 398\\
0.934579 395\\
0.943396 392\\
0.952381 388\\
0.961538 383\\
0.970874 373\\
0.980392 361\\
0.990099 351\\
1 360\\
1.0101 356\\
1.02041 353\\
1.03093 349\\
1.04167 348\\
1.05263 348\\
1.06383 347\\
1.07527 345\\
1.08696 345\\
1.0989 348\\
1.11111 349\\
1.1236 350\\
1.13636 349\\
1.14943 351\\
1.16279 373\\
1.17647 384\\
1.19048 394\\
1.20482 412\\
1.21951 406\\
1.23457 405\\
1.25 415\\
};

\addplot [
color=black,
solid,
line width=1.0pt,
mark=+,
mark size = 2.5,
mark options={solid,red}
]
table[row sep=crcr]{
0.840336 897\\
0.847458 907\\
0.854701 892\\
0.862069 905\\
0.869565 890\\
0.877193 902\\
0.884956 896\\
0.892857 910\\
0.900901 876\\
0.909091 846\\
0.917431 810\\
0.925926 808\\
0.934579 812\\
0.943396 808\\
0.952381 798\\
0.961538 805\\
0.970874 818\\
0.980392 804\\
0.990099 804\\
1 813\\
1.0101 783\\
1.02041 772\\
1.03093 775\\
1.04167 763\\
1.05263 755\\
1.06383 718\\
1.07527 722\\
1.08696 708\\
1.0989 727\\
1.11111 762\\
1.1236 764\\
1.13636 759\\
1.14943 774\\
1.16279 783\\
1.17647 794\\
1.19048 797\\
1.20482 822\\
1.21951 866\\
1.23457 881\\
1.25 876\\
};

\addplot [
color=black!70,
solid,
line width=1.0pt,
mark=star,
mark size = 2.5,
mark options={solid,orange}
]
table[row sep=crcr]{
0.769231 488\\
0.775194 484\\
0.78125 484\\
0.787402 476\\
0.793651 471\\
0.8 461\\
0.806452 447\\
0.813008 446\\
0.819672 446\\
0.826446 444\\
0.833333 444\\
0.840336 443\\
0.847458 441\\
0.854701 436\\
0.862069 434\\
0.869565 424\\
0.877193 412\\
0.884956 397\\
0.892857 396\\
0.900901 396\\
0.909091 396\\
0.917431 394\\
0.925926 392\\
0.934579 391\\
0.943396 384\\
0.952381 371\\
0.961538 348\\
0.970874 347\\
0.980392 346\\
0.990099 346\\
1 344\\
1.0101 340\\
1.02041 331\\
1.03093 317\\
1.04167 302\\
1.05263 298\\
1.06383 298\\
1.07527 298\\
1.08696 320\\
1.0989 341\\
1.11111 346\\
1.1236 346\\
1.13636 347\\
1.14943 358\\
1.16279 382\\
1.17647 392\\
1.19048 394\\
1.20482 396\\
1.21951 396\\
1.23457 396\\
1.25 422\\
1.26582 437\\
1.28205 443\\
1.2987 446\\
1.31579 446\\
1.33333 449\\
1.35135 471\\
1.36986 481\\
1.38889 489\\
1.40845 491\\
1.42857 493\\
};

\end{axis}
\end{tikzpicture}
\caption{Total fixed-stress iterations for varying step length size $\tau_n$ for dG(1) in time.}
\label{fig:4:dG1:tau}
\end{figure}

\begin{figure}[p]
\centering
%
\begin{tikzpicture}

\begin{axis}[%
width=4.82222222222222in,
height=3.5in,
scale only axis,
/pgf/number format/.cd, 1000 sep={},
xlabel={$\omega$},
ylabel={iterations},
xmin=0.84,
xmax=1.25,
ymin=100,
ymax=1100,
yminorticks=true,
legend style={legend pos=north east},
legend style={draw=black,fill=white,legend cell align=left},
legend entries = {
  {cGP(2), $m=2$, $s=2$, $\tau_n=0.005$},
  {cGP(2), $m=2$, $s=2$, $\tau_n=0.010$},
  {\textcolor{black!70}{cGP(1), $m=2$, $s=2$, $\tau_n=0.010$}},
  {cGP(2), $m=2$, $s=2$, $\tau_n=0.020$}
}
]

\addplot [
color=black,
solid,
line width=1.0pt,
mark=+,
mark size = 2.5,
mark options={solid,red}
]
table[row sep=crcr]{
0.0 0\\
0.1 0\\
};

\addplot [
color=black,
solid,
line width=1.0pt,
mark=o,
mark size = 2.5,
mark options={solid,green}
]
table[row sep=crcr]{
0.0 0\\
0.1 0\\
};

\addplot [
color=black!70,
solid,
line width=1.0pt,
mark=star,
mark size = 2.5,
mark options={solid,orange}
]
table[row sep=crcr]{
0.0 0\\
0.1 0\\
};

\addplot [
color=black,
solid,
line width=1.0pt,
mark=x,
mark size = 2.5,
mark options={solid,blue}
]
table[row sep=crcr]{
0.0 0\\
0.1 0\\
};

\addplot [
color=black,
dotted,
line width=1.0pt
]
table[row sep=crcr]{
1.0 0\\
1.0 2000\\
};

\addplot [
color=black,
solid,
line width=1.0pt
]
table[row sep=crcr]{
0.840336 222\\
0.847458 221\\
0.854701 219\\
0.862069 219\\
0.869565 214\\
0.877193 207\\
0.884956 199\\
0.892857 199\\
0.900901 199\\
0.909091 199\\
0.917431 200\\
0.925926 198\\
0.934579 200\\
0.943396 201\\
0.952381 188\\
0.961538 179\\
0.970874 174\\
0.980392 174\\
0.990099 174\\
1 173\\
1.0101 172\\
1.02041 171\\
1.03093 169\\
1.04167 169\\
1.05263 167\\
1.06383 167\\
1.07527 168\\
1.08696 166\\
1.0989 169\\
1.11111 172\\
1.1236 173\\
1.13636 174\\
1.14943 174\\
1.16279 174\\
1.17647 184\\
1.19048 192\\
1.20482 197\\
1.21951 199\\
1.23457 199\\
1.25 199\\
};

\addplot [
color=black,
solid,
line width=1.0pt
]
table[row sep=crcr]{
0.840336 442\\
0.847458 438\\
0.854701 435\\
0.862069 428\\
0.869565 423\\
0.877193 401\\
0.884956 400\\
0.892857 398\\
0.900901 399\\
0.909091 400\\
0.917431 399\\
0.925926 396\\
0.934579 393\\
0.943396 390\\
0.952381 387\\
0.961538 385\\
0.970874 384\\
0.980392 379\\
0.990099 376\\
1 370\\
1.0101 368\\
1.02041 357\\
1.03093 348\\
1.04167 347\\
1.05263 347\\
1.06383 347\\
1.07527 347\\
1.08696 346\\
1.0989 346\\
1.11111 342\\
1.1236 345\\
1.13636 347\\
1.14943 349\\
1.16279 348\\
1.17647 350\\
1.19048 372\\
1.20482 385\\
1.21951 390\\
1.23457 395\\
1.25 397\\
};

\addplot [
color=black,
solid,
line width=1.0pt
]
table[row sep=crcr]{
0.840336 878\\
0.847458 874\\
0.854701 868\\
0.862069 863\\
0.869565 857\\
0.877193 857\\
0.884956 846\\
0.892857 839\\
0.900901 824\\
0.909091 810\\
0.917431 810\\
0.925926 792\\
0.934579 787\\
0.943396 786\\
0.952381 785\\
0.961538 784\\
0.970874 780\\
0.980392 778\\
0.990099 773\\
1 769\\
1.0101 763\\
1.02041 757\\
1.03093 746\\
1.04167 732\\
1.05263 693\\
1.06383 694\\
1.07527 692\\
1.08696 689\\
1.0989 687\\
1.11111 685\\
1.1236 685\\
1.13636 685\\
1.14943 691\\
1.16279 694\\
1.17647 697\\
1.19048 698\\
1.20482 749\\
1.21951 773\\
1.23457 781\\
1.25 788\\
};

\addplot [
color=black!70,
solid,
line width=1.0pt
]
table[row sep=crcr]{
0.769231 488\\
0.775194 484\\
0.78125 484\\
0.787402 476\\
0.793651 471\\
0.8 461\\
0.806452 447\\
0.813008 446\\
0.819672 446\\
0.826446 444\\
0.833333 444\\
0.840336 443\\
0.847458 441\\
0.854701 436\\
0.862069 434\\
0.869565 424\\
0.877193 412\\
0.884956 397\\
0.892857 396\\
0.900901 396\\
0.909091 396\\
0.917431 394\\
0.925926 392\\
0.934579 391\\
0.943396 384\\
0.952381 371\\
0.961538 348\\
0.970874 347\\
0.980392 346\\
0.990099 346\\
1 344\\
1.0101 340\\
1.02041 331\\
1.03093 317\\
1.04167 302\\
1.05263 298\\
1.06383 298\\
1.07527 298\\
1.08696 320\\
1.0989 341\\
1.11111 346\\
1.1236 346\\
1.13636 347\\
1.14943 358\\
1.16279 382\\
1.17647 392\\
1.19048 394\\
1.20482 396\\
1.21951 396\\
1.23457 396\\
1.25 422\\
1.26582 437\\
1.28205 443\\
1.2987 446\\
1.31579 446\\
1.33333 449\\
1.35135 471\\
1.36986 481\\
1.38889 489\\
1.40845 491\\
1.42857 493\\
};

\addplot [
color=black,
solid,
line width=1.0pt,
mark=x,
mark size = 2.5,
mark options={solid,blue}
]
table[row sep=crcr]{
0.840336 222\\
0.847458 221\\
0.854701 219\\
0.862069 219\\
0.869565 214\\
0.877193 207\\
0.884956 199\\
0.892857 199\\
0.900901 199\\
0.909091 199\\
0.917431 200\\
0.925926 198\\
0.934579 200\\
0.943396 201\\
0.952381 188\\
0.961538 179\\
0.970874 174\\
0.980392 174\\
0.990099 174\\
1 173\\
1.0101 172\\
1.02041 171\\
1.03093 169\\
1.04167 169\\
1.05263 167\\
1.06383 167\\
1.07527 168\\
1.08696 166\\
1.0989 169\\
1.11111 172\\
1.1236 173\\
1.13636 174\\
1.14943 174\\
1.16279 174\\
1.17647 184\\
1.19048 192\\
1.20482 197\\
1.21951 199\\
1.23457 199\\
1.25 199\\
};

\addplot [
color=black,
solid,
line width=1.0pt,
mark=o,
mark size = 2.5,
mark options={solid,green}
]
table[row sep=crcr]{
0.840336 442\\
0.847458 438\\
0.854701 435\\
0.862069 428\\
0.869565 423\\
0.877193 401\\
0.884956 400\\
0.892857 398\\
0.900901 399\\
0.909091 400\\
0.917431 399\\
0.925926 396\\
0.934579 393\\
0.943396 390\\
0.952381 387\\
0.961538 385\\
0.970874 384\\
0.980392 379\\
0.990099 376\\
1 370\\
1.0101 368\\
1.02041 357\\
1.03093 348\\
1.04167 347\\
1.05263 347\\
1.06383 347\\
1.07527 347\\
1.08696 346\\
1.0989 346\\
1.11111 342\\
1.1236 345\\
1.13636 347\\
1.14943 349\\
1.16279 348\\
1.17647 350\\
1.19048 372\\
1.20482 385\\
1.21951 390\\
1.23457 395\\
1.25 397\\
};

\addplot [
color=black,
solid,
line width=1.0pt,
mark=+,
mark size = 2.5,
mark options={solid,red}
]
table[row sep=crcr]{
0.840336 878\\
0.847458 874\\
0.854701 868\\
0.862069 863\\
0.869565 857\\
0.877193 857\\
0.884956 846\\
0.892857 839\\
0.900901 824\\
0.909091 810\\
0.917431 810\\
0.925926 792\\
0.934579 787\\
0.943396 786\\
0.952381 785\\
0.961538 784\\
0.970874 780\\
0.980392 778\\
0.990099 773\\
1 769\\
1.0101 763\\
1.02041 757\\
1.03093 746\\
1.04167 732\\
1.05263 693\\
1.06383 694\\
1.07527 692\\
1.08696 689\\
1.0989 687\\
1.11111 685\\
1.1236 685\\
1.13636 685\\
1.14943 691\\
1.16279 694\\
1.17647 697\\
1.19048 698\\
1.20482 749\\
1.21951 773\\
1.23457 781\\
1.25 788\\
};

\addplot [
color=black!70,
solid,
line width=1.0pt,
mark=star,
mark size = 2.5,
mark options={solid,orange}
]
table[row sep=crcr]{
0.769231 488\\
0.775194 484\\
0.78125 484\\
0.787402 476\\
0.793651 471\\
0.8 461\\
0.806452 447\\
0.813008 446\\
0.819672 446\\
0.826446 444\\
0.833333 444\\
0.840336 443\\
0.847458 441\\
0.854701 436\\
0.862069 434\\
0.869565 424\\
0.877193 412\\
0.884956 397\\
0.892857 396\\
0.900901 396\\
0.909091 396\\
0.917431 394\\
0.925926 392\\
0.934579 391\\
0.943396 384\\
0.952381 371\\
0.961538 348\\
0.970874 347\\
0.980392 346\\
0.990099 346\\
1 344\\
1.0101 340\\
1.02041 331\\
1.03093 317\\
1.04167 302\\
1.05263 298\\
1.06383 298\\
1.07527 298\\
1.08696 320\\
1.0989 341\\
1.11111 346\\
1.1236 346\\
1.13636 347\\
1.14943 358\\
1.16279 382\\
1.17647 392\\
1.19048 394\\
1.20482 396\\
1.21951 396\\
1.23457 396\\
1.25 422\\
1.26582 437\\
1.28205 443\\
1.2987 446\\
1.31579 446\\
1.33333 449\\
1.35135 471\\
1.36986 481\\
1.38889 489\\
1.40845 491\\
1.42857 493\\
};

\end{axis}
\end{tikzpicture}
\caption{Total fixed-stress iterations for varying step length size $\tau_n$ for cGP(2) in 
time.}
\label{fig:5:cG2:tau}
\end{figure}